\documentclass[mathpazo]{cicp}

\usepackage{microtype}
\usepackage{mlmath}
\usepackage{subfigure}
\usepackage{float}

\makeatletter

\usepackage{url}
\usepackage{booktabs}
\usepackage{nicefrac}
\usepackage{caption}
\usepackage{xcolor}

\definecolor{zxcol}{rgb}{0, 1, 0}

\definecolor{zllcol}{rgb}{0, 0, 1}

\@ifundefined{showcaptionsetup}{}{%
 \PassOptionsToPackage{caption=false}{subfig}}
\usepackage{subfig}
\makeatother

\graphicspath{%
    {./figures/}
    {./figures/Poisson_figure/}
    {./figures/deepritz_laplace/}
    {./figures/constructed_solver_case3/}
    {./figures/bolzmann_continuous_fixed/}
    {./figures/boltzmann_sigma/}
    {./figures/boltzmann_continuous_vary/}
    {./figures/burgers_operator_upwind/}
    {./figures/nonlinear_laplace/}
}
\begin{document}

\title{MOD-Net: A Machine Learning Approach via Model-Operator-Data Network for Solving PDEs}



\author[Zhang L L et.~al.]{Lulu Zhang\affil{1}, Tao Luo\affil{2,3}, Yaoyu Zhang\affil{1,3,4}, Weinan E\affil{5,6}, Zhi-Qin John Xu\affil{1,3}\comma\footnotemark, Zheng Ma\affil{2,3}\comma\corrauth}
\address{
    \affilnum{1}\ Institute of Natural Sciences and School of Mathematical Sciences, Shanghai Jiao Tong University, Shanghai, 200240, China. \\
    \affilnum{2}\ School of Mathematical Sciences and Institute of Natural Sciences, Shanghai Jiao Tong University, CMA-Shanghai, Shanghai, 200240, China. \\
    \affilnum{3}\ MOE-LSC and Qing Yuan Research Institute, Shanghai Jiao Tong University, Shanghai, 200240, China. \\
    \affilnum{4}\ Shanghai Center for Brain Science and Brain-Inspired Technology, Shanghai, 200031, China. \\
    \affilnum{5}\ School of Mathematical Sciences, Peking University, Beijing, 100871, China. \\
    \affilnum{6}\ AI for Science Institute, Beijing, 100080, China.
}

\emails{ %
    {\tt xuzhiqin@sjtu.edu.cn} (Z.~Xu), {\tt zhengma@sjtu.edu.cn} (Z.~Ma)
}

\begin{abstract}
    In this paper, we propose a a machine learning approach via model-operator-data network (MOD-Net) for solving PDEs. A MOD-Net is driven by a model to solve PDEs based on operator representation with regularization from data.
    For linear PDEs, we use a DNN to parameterize the Green's function and obtain the neural operator to approximate the solution according to the Green's method. To train the DNN, the empirical risk consists of the mean squared loss with the least square formulation or the variational formulation of the governing equation and
    boundary conditions. For complicated problems, the empirical risk also includes a few labels, which are computed
    on coarse grid points with cheap computation cost and significantly improves the model accuracy.
    Intuitively, the labeled dataset works as a regularization in addition to the model constraints.
    The MOD-Net solves a family of PDEs rather than a specific one and is much more efficient than original neural operator because few expensive labels are required.
    We numerically show MOD-Net is very efficient in solving Poisson equation and one-dimensional  radiative transfer equation. For nonlinear PDEs,
    the nonlinear MOD-Net can be similarly used as an ansatz for solving nonlinear PDEs, exemplified by solving several nonlinear PDE problems, such  as the Burgers equation.
\end{abstract}

\maketitle

\section{Introduction}
Nowadays, using deep neural networks (DNNs) to solve PDEs attracts more and more attention~\cite{weinan2020machine,E2020algorithms,weinan2017deep, weinan2018deep, he2020relu, liao2021deep, raissi2019physics, hamilton2019dnn, strofer2019data,liu2020multi,li2020multi,li2020neural,li2020fourier,lu2019deeponet,markidis2021physics}. Here we review three DNN approaches for solving PDEs.

The first approach is to parameterize the solution by a DNN and use the mean square of the residual of the equation~\cite{dissanayake1994neural,raissi2019physics} or the variational forms~\cite{weinan2018deep,liao2021deep} as risk or loss, by minimizing which the DNN output satisfies PDE. A comprehensive overview can be found in~\cite{E2020algorithms}.
This parameterization approach can solve very high-dimensional PDEs and does not require any labels.
However, it only solves a specific PDE during each training trial, that is, if the PDE setup changes, such as the source terms, the boundary conditions or other parameters in the PDE, we have to train a new DNN.
An important characteristic of the parameterization approach is slow learning of the high frequency part as indicated by the frequency principle~\cite{xu_training_2018,xu2019frequency,markidis2021physics}.
To overcome the curse of high frequency, a series of multiscale approaches are proposed~\cite{cai2020phase,liu2020multi,li2020multi,Matthew2020Fourier,wang2021eigenvector}.
The second approach uses DNN to learn the mapping from the source term to the solution~\cite{fan2019multiscale}.
In this approach, the source function and the solution are sampled at fixed grid points as two vectors.
Then the vector of the source function is fed into the DNN to predict the vector of the solution function.
The advantage of this mapping approach is that the DNN solves the PDE for any source function, thus it can be very convenient in application.
However, the mapping approach can only evaluate the solution at fixed points. DeepOnet~\cite{lu2019deeponet} is proposed that the source function is still fed into the network on fixed grid points but the output can be evaluated on any points by adding on extra inputs of the points to the network.
Such approach requires very large sample points, which is often computational inefficient or intractable, especially when dealing with high-dimensional PDEs or complicated integro-differential equations, such as Boltzmann equation and radiative transfer equation (RTE).
The third approach is called neural operator~\cite{li2020neural,li2020fourier}, which represents the solution based on the form similar to the idea of the Green's function and the DNN is used to parametrize the Green's function.
The neural operator solves a type of PDEs but not a specific PDE and can be evaluated at any time or spatial points.
Training of the neural operator is to minimize the difference between the learned solution and the true solution at randomly sampled points.
Therefore, the neural operator is a data-driven method and requires a large amount of labels, a similar difficulty to the mapping approach.

In this work, we propose a machine learning approach via model-operator-data network (MOD-Net) for solving PDEs.
The MOD-Net has advantages including: (i) obtaining a functional representation of the solution which allows evaluating the solution at any points; (ii) requiring none or few labels numerically computed by a traditional scheme on coarse grid points with cheap computation; (iii) solving a family of PDEs but not a specific PDE.
The three key components of MOD-Net are illustrated as follows.

\textbf{Model driven.} MOD-Net is driven by the physical model to avoid using too much expensive labeled data.
That is, the empirical risk, i.e., training loss requires the solution satisfying the constraints of the governing PDE or equivalent forms and boundary conditions.
To realize the model constraint, one can use various methods, such as minimizing the mean square of the residual of the governing PDE and boundary conditions (e.g., physics-informed neural network~\cite{dissanayake1994neural,raissi2019physics}), or minimizing the variational form of the governing PDE and boundary conditions (e.g., Deep-Ritz method~\cite{weinan2018deep}).

\textbf{Operator representation.} Similar to the neural operator~\cite{li2020neural,li2020fourier}, the MOD-Net represents the solution operator of a PDE, i.e, mapping from source terms, boundary conditions, or parameters to the solution.
In this work, DNN is used to parameterize the Green's function, however, it is not restricted to use Green's function and can be generalized to other architectures.
This operator representation utilizes the invariant characteristic of Green's function in solving PDE, thus, might be more efficient than an end-to-end representation by parameterizing solution with a DNN directly.

\textbf{Data regularization.} In MOD-Net, we find that in complicated problems, with only the model constraints, the solution is often very inaccurate even when the empirical risk is reasonably small.
For example, the radiative/linear transport equations only have a hypercoercive integro-differential operator instead of a nice coercive operator like common elliptic equations.
This degeneracy property makes the velocity space may have bad regularity provided some singular coefficients in the equation thus leads to the existence of many ``weak solutions'' to choose from.
To overcome this problem, we add a regularization term by minimizing the difference between the MOD-Net prediction and a few labels numerically computed by a traditional scheme on coarse grid points with cheap computation cost.
Note that with only the small amount of labeled data, the MOD-Net cannot be well trained either.
Therefore, the effect of the labeled data in MOD-Net is different from supervised learning, in which a DNN training often requires a large amount of accurate labeled data.
Intuitively, the effect of data in MOD-Net works as an regularization similar to various regularization terms in traditional optimization problems.

We first apply MOD-Net to solve simple Poisson equations, in which we show that without labels, MOD-Net with the mean square loss or the variational loss of PDE, i.e., governing equation and boundary conditions, can learn the Poisson equation well.
We further apply MOD-Net to a class of equations controlled by parameters, which can be regarded as PDEs with uncertainty or a simplification of complicated control problem.
For these equations, we can not train the MOD-Net well with only the physical information, however, with the model information and a few labels, we can well train MOD-Net.
The data regularization also significantly improves the model accuracy for the RTE~\cite{chandrasekhar2013radiative,lenoble1985radiative}, which is important in real applications, such as simulation of nuclear reactor, optical tomography and radiation therapy.
Besides for these linear PDEs, we also apply nonlinear MOD-Net to the one-dimensional Burgers equation and two-dimensional nonlinear equations.
In these two cases, we use no labeled data and only utilize the information of the PDE and we can also well train the nonlinear MOD-Net.

The rest of the paper is organized as follows.
In section~\ref{sec:DNN}, we will give a brief introduction of DNNs.
Section~\ref{sec:MOD-Net} will present MOD-Net structures.
Section~\ref{sec:poi} will show the numerical results for Poisson equations.
In section~\ref{sec:toyexample}, we show the numerical experiments for constructed toy equations.
In section~\ref{sec:boltz}, we show numerical experiments for one-dimensional RTE.
In section~\ref{sec:Burgers}, we show numerical results for one-dimensional Burgers equation.
Section~\ref{sec:nonlinear} will show the numerical experiments for two-dimensional nonlinear equation.
Finally, section~\ref{sec:con} gives a conclusion and some discussions for future work.
\section{Preliminary: Deep neural networks}
\label{sec:DNN}
We introduce the following conventional notations for DNNs\footnote{BAAI.2020.\ Suggested Notation for Machine Learning.\ https://github.com/mazhengcn/suggested-notation-for-machine-learning.}.
An $L$-layer neural network is defined recursively as,
\begin{equation}
    \begin{aligned}
         & \vf_{\vtheta}^{[0]}(\vx)=\vx,                                                                                      \\
         & \vf_{\vtheta}^{[l]}(\vx)=\sigma\circ(\mW^{[l-1]} \vf_{\vtheta}^{[l-1]}(\vx) + \vb^{[l-1]}), \quad 1\leq l\leq L-1, \\
         & \vf_{\vtheta}(\vx)=\vf_{\vtheta}^{[L]}(\vx)=\mW^{[L-1]} \vf_{\vtheta}^{[L-1]}(\vx) + \vb^{[L-1]},
    \end{aligned}
\end{equation}
where $\mW^{[l]} \in \sR^{m_{l+1}\times m_{l}}$, $\vb^{[l]}=\sR^{m_{l+1}}$, $m_0=d_{\rm in}=d$ is the input dimension, $m_{L}=d_{\rm o}$ is the output dimension, $\sigma$ is a scalar function and ``$\circ$'' means entry-wise operation.
We denote the set of parameters by $\vtheta$.

A loss function $\ell(f_{\vtheta}(\vx),\vy)$ measures the difference between a prediction and a true label.
The empirical risk, also known as the training loss for a set $S=\{(\vx_i,\vy_i)\}_{i=1}^{n}$ is denoted by $\RS(\vtheta)$,
\begin{equation}
    \RS(\vtheta) = \frac{1}{n}\sum_{i=1}^n\ell(f_{\vtheta}(\vx_i),\vy_i).
\end{equation}
More generally, the empirical risk can be defined without labels. Two common empirical risks used for solving PDEs are least square loss (e.g., physics-informed neural network~\cite{dissanayake1994neural,raissi2019physics}) and variational loss (e.g., Deep-Ritz method~\cite{weinan2018deep}).
For example, if we want to use DNN $f_{\vtheta}$ to learn the solution $u(\vx)$ of an equation such as $\fL[u](\vx)=\vzero$  for $\vx\in\Omega$, the empirical risk with least square loss for this equation can be defined by
\begin{equation}
    \RS(\vtheta) =\frac{1}{n}\sum_{i=1}^n\|\fL[f_{\vtheta}](\vx_i)\|_2^2,
\end{equation}
where the dataset $\{\vx_i\}_{i=1}^n$ is randomly sampled from $\Omega$ at each iteration step.
For some problems, we can use the variational form of the equation to define the loss.
More precisely if the solution is the minimizer of the functional $I[u](x)$ whose density denoted by $\ell[u](x)$, then the variational loss can be defined by
\begin{equation}
    \RS(\vtheta) =\frac{1}{n}\sum_{i=1}^n \ell[f_{\vtheta}](\vx_i),
\end{equation}
where the dataset $\{\vx_i\}_{i=1}^n$ is randomly sampled from $\Omega$ at each iteration step.

The empirical risk can also be defined by the weighted summation of the empirical risk of the labeled data and the empirical risk of the equation.
If more constraints needed in the problem, such as boundary conditions, the empirical risk can be similarly adopted.

\section{Model-operator-data network (MOD-Net)}
\label{sec:MOD-Net}
Our goal is to solve the following PDE efficiently and accurately,
\begin{equation}\label{linearPDE}
    \left\{
    \begin{aligned}
         & \fL [u](\vx) = g(\vx),\quad \vx \in \Omega,        \\
         & u(\vx) = \phi(\vx), \quad \vx \in \partial \Omega,
    \end{aligned}
    \right.
\end{equation}
where $\fL$ is the operator that can be a usual differential operator or even integro-differential operator.
Our basic idea is to use a DNN to learn the operator $\fG: (\phi,g) \mapsto u$, i.e., for each given boundary condition $\phi$ and source term $g$, there is $\fG(\phi,g) = u$.

For a linear PDE, Green's function can help us obtain the solution of PDE due to the superposition principle.
With the Green's function method, we have the following representation of the solution of PDE~\eqref{linearPDE},
\begin{equation}\label{representation}
    u(\vx;\phi,g) = \int_{ \Omega }G_1(\vx,\vx') g(\vx') \diff{\vx'} + \int_{\partial \Omega} G_2(\vx,\vx') \phi(\vx') \diff{\vx'}.
\end{equation}
where for any fixed $\vx'\in \Omega$, $G_1(\vx,\vx')$ is the solution of the following equation,
\begin{equation*}
    \left\{
    \begin{aligned}
         & \fL [G_1](\vx) = \delta(\vx-\vx'),\quad \vx \in \Omega, \\
         & G_1(\vx,\vx') = 0, \quad \vx \in \partial \Omega,
    \end{aligned}
    \right.
\end{equation*}
and for any fixed $\vx'\in \partial \Omega$, $G_2(\vx,\vx')$ is the solution of the following equation,
\begin{equation*}
    \left\{
    \begin{aligned}
         & \fL [G_2](\vx) = 0,\quad \vx \in \Omega,                         \\
         & G_2(\vx,\vx') = \delta(\vx-\vx'), \quad \vx \in \partial \Omega.
    \end{aligned}
    \right.
\end{equation*}
For the nonlinear PDE, we can extend the Green's function method for the nonlinear case and  use the following representation,
\begin{equation*}\label{nonlinear_representation}
    \begin{aligned}
        u(\vx;\phi,g) = F\left(\int_{ \Omega }
        G_1(\vx,\vx') g(\vx') \diff{\vx'} +
        \int_{\partial \Omega }
        G_2(\vx,\vx') \phi(\vx') \diff{\vx'}\right),
    \end{aligned}
\end{equation*}
where $F(\vx)$ is a nonlinear function and is represented by DNN $F_{\vtheta}(\vx)$ in this work.

However, it is difficult to obtain the analytical formula of the operators $G_1$ and $G_2$. In the following, we consider using DNN to represent operators $G_1$ and $G_2$, i.e., a DNN $G_{\vtheta_1}(\vx,\vx')$ is trained to represent
$G_1(\vx,\vx')$, similarly, another DNN  $G_{\vtheta_2}(\vx,\vx')$ is used for $G_2(\vx,\vx')$. By implementing $G_{\vtheta_1}(\vx,\vx')$ and $G_{\vtheta_2}(\vx,\vx')$ into Eq. (\ref{representation}) and Eq. (\ref{nonlinear_representation}), we obtain a DNN representation for $u(\vx;\phi,g)$ as $u_{\vtheta_1,\vtheta_2}(\vx;\phi,g)$. In application, the integration in  Eq. (\ref{representation}) and Eq. (\ref{nonlinear_representation}) is realized by discrete numerical schemes. For example, we consider Monte-Carlo algorithm where we uniformly sample a set $S_{G,\Omega}$ from $\Omega$ and a set  $S_{G,\partial\Omega}$ from $\partial\Omega$, then, for the linear PDE,
\begin{equation}
    u_{\vtheta_1,\vtheta_2}(\vx;\phi,g) = \frac{|\Omega|}{|S_{G,\Omega}|}\sum_{\vx'\in S_{G,\Omega}}
    G_{\vtheta_1}(\vx,\vx') g(\vx')
    +\frac{|\partial\Omega|}{|S_{G,\partial\Omega}|} \sum_{\vx'\in S_{G,\partial\Omega}} G_{\vtheta_2}(\vx,\vx')\phi(\vx'),
\end{equation}
and for the nonlinear PDE,
\begin{equation}
    u_{\vtheta_1,\vtheta_2}(\vx;\phi,g) = F_{\vtheta}\Big( \frac{|\Omega|}{|S_{G,\Omega}|}\sum_{\vx'\in S_{G,\Omega}}
    G_{\vtheta_1}(\vx,\vx') g(\vx')
    +\frac{|\partial\Omega|}{|S_{G,\partial\Omega}|} \sum_{\vx'\in S_{G,\partial\Omega}} G_{\vtheta_2}(\vx,\vx')\phi(\vx')\Big).
\end{equation}

\begin{figure}[H]
    \centering
    \includegraphics[scale=0.06]{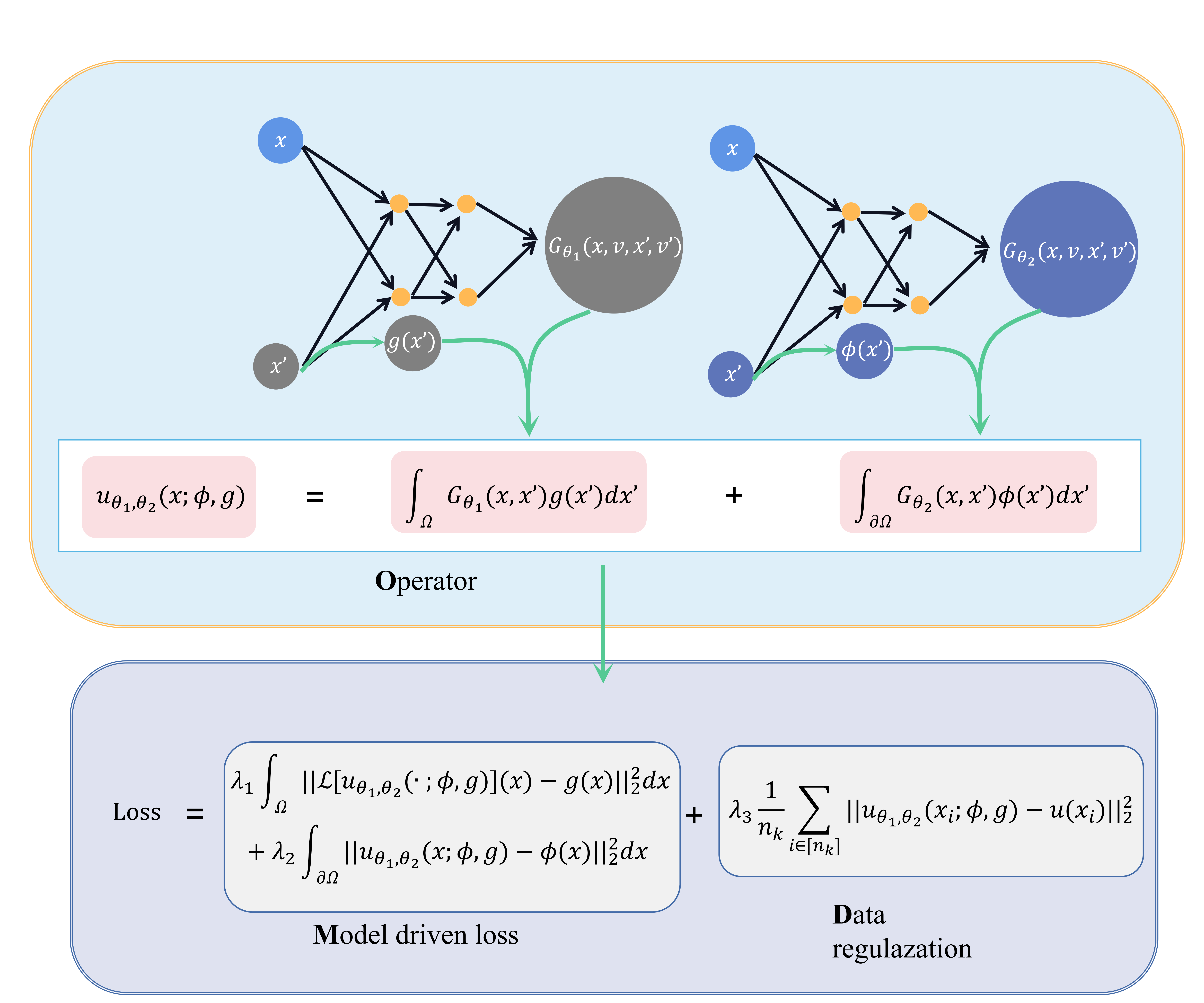}
    \caption{Schematic of MOD-Net approach.  In MOD-Net approach, we use the DNN $G_{\vtheta_1}$ and $G_{\vtheta_2}$ to parameter the Green's function $G_1$ and $G_2$, respectively. And according to the Green's formula, we obtain $u_{\vtheta_1,\vtheta_2}$ which can approximate solution operator. Then we train these DNNs utilizing the information of PDE, i.e., governing equation and boundary conditions, and a few data.
    }
    \label{fig:schematic}
\end{figure}
To train the neural networks, we would utilize the information of PDE, i.e., governing equation and boundary conditions, and a few data $S^{u,k}=\{\vx_i,u^{k}(\vx_i)\}_{i\in [n_k]}$ for each $\{\phi^{k},g^k\}, k=1,2,\cdots,K$, where $K$ is the total number of examples/observations and $u^{k}(\cdot)=u(\cdot;\,\phi^{k},g^k)$.
Note that $S^{u,k}$ can be numerically solved by traditional schemes on coarse grid points, which is not computationally expensive or even obtained from experiment observations.
To utilize the constraint of governing PDE, for each $k$, we uniformly sample a set of data from $\Omega$, i.e.,  $S^{\Omega,k}$.
To utilize the information of boundary constraint, for each $k$, we uniformly sample a set of data from $\partial\Omega$, i.e., $S^{\partial\Omega,k}$.
Then, we train the neural networks by minimizing the empirical risk defined as follows,
\begin{equation}
    \begin{aligned}
        \RS
         & = \frac{1}{K}\sum_{k\in[K]}  \Big(
        \lambda_1 \frac{1}{|S^{\Omega,k}|} \sum_{\vx\in S^{\Omega,k}}
        \| \fL [u_{\vtheta_1,\vtheta_2}(\vx;\phi^{k},g^k)](\vx)-g^k(\vx) \|_2^2
        \\
         & ~~~~+\lambda_2 \frac{1}{|S^{\partial\Omega,k}|} \sum_{\vx\in S^{\partial\Omega,k}}  \| u_{\vtheta_1,\vtheta_2}(\vx;\phi^{k},g^k) - \phi^k(\vx) \|_2^2
        \\
         & ~~~~+ \lambda_3 \frac{1}{n_k} \sum_{i\in [n_k]} \|u_{\vtheta_1,\vtheta_2}(\vx_i;\phi^{k},g^k)-u^{k}(\vx_i)\|_2^2
        \Big),
        \label{lossfunctiondefinition}
    \end{aligned}
\end{equation}
where $\lambda_1,\lambda_2,\lambda_3$ are hyperparameters used to tune the weight of each part in the total risk.

Here we remark that the least square loss is not crucial, we can use other loss such as variational loss, see example~2.
Also the labeled data is not restricted the true solution, we can use other data as regularization term, such as the macroscopic quantities, e.g., the density function $\rho(x_i)$ in the RTE, which is a moment of the solution.
An important advantage of our proposed MOD-Net method is that we take advantage of the PDE constraint and use cheap not-so-accurate labeled data. For convenience, the notations are listed in Table~\ref{table1}.

In real applications, depending on the problems, if the boundary condition is zero, such as the Poisson case in section~\ref{sec:poi}, $G_2$ is ignored, and if the source term is zero, such as the RTE case in section~\ref{sec:boltz}, $G_1$ is ignored. For one-dimensional case, the integration over boundary is the summation at two points, therefore, we can use two DNNs to learn $G_2(\vx,\vx')|_{\vx^\prime = \vx_{L}}$ and $G_2(\vx,\vx')|_{\vx^\prime = \vx_{R}}$, respectively, where $\vx_{L}$ and $\vx_{R}$ are boundary points. The input of the two DNNs are lower dimensional due to the fixation of $\vx^\prime$. The amount of the labeled data required also depends on the problem. For simple problems, such as Poisson problem, we use no  labeled data, but for complicated constructed toy example or RTE problem, we use a few labeled data.

\begin{table}[htbp]
    \centering
    \caption{Notation}
    \label{table1}
    \begin{tabular}{|c|c|}

        \hline

        $\fL$                                                & PDE operator                                                                   \\

        \hline

        $g$                                                  & source term in PDE                                                             \\
        \hline

        $\phi$                                               & boundary value in PDE                                                          \\
        \hline

        $\fG$                                                & operator to be learned, general definition is $\fG: (\phi,\sigma,g) \mapsto u$ \\
        \hline

        $[n]$                                                & index set $\{1,2,...,n\}$                                                      \\
        \hline

        $S_{G_x},S_{G_y},S_{G_v}^+,S_{G_v}^-$                & the set of integration points in solution's representation                     \\
        \hline

        $S_{v}$                                              & the set of integration points in $v$ direction                                 \\
        \hline

        $\omega$                                             & integration coefficients                                                       \\
        \hline
        $S^{\Omega}$                                         & a set of data uniformly sampled from $\Omega$                                  \\
        \hline

        $S^{\partial\Omega}$                                 & a set of data uniformly sampled from $\partial\Omega$                          \\
        \hline

        $S^{u,k}=\{\vx_i,u^k(\vx_i)\}_{i=1}^{n_k}$           & labeled data set for $k$th  PDE, where $u^k(\cdot)$ is solution of $k$th PDE   \\
        \hline
        $F_{\vtheta}(\vx),G_{\vtheta}(\vx),u_{\vtheta}(\vx)$ & DNN                                                                            \\
        \hline
    \end{tabular}
\end{table}

\section{Numerical experiments: $2$D Poisson equation}
\label{sec:poi}
We consider the Poisson equations in $2$D,
\begin{equation}\label{Poisson}
    \begin{aligned}
        -\Delta u(\vx) & =g(\vx), \quad \vx\in \Omega,
        \\
        u(\vx)         & =0, \quad \vx\in \partial \Omega.
    \end{aligned}
\end{equation}
where the source function $g(x,y)=-a(x^2-x+y^2-y)$, i.e.,
\begin{equation}\label{2D_Poisson_case}
    \begin{aligned}
        -(\partial_{xx}u+\partial_{yy}u) & =-a(x^2-x+y^2-y), \quad (x,y)\in \Omega,
        \\
        u                                & = 0, \quad (x,y)\in \partial \Omega,
    \end{aligned}
\end{equation}
where $\Omega=[0,1]^2$ and the constant $a$ controls the source term.
Obviously, the analytical solution is $u(x,y;g)=\frac{a}{2}x(x-1)y(y-1)$.

\subsection{Use DNN to fit Green's function}
For Poisson equation, which is a linear PDE, using the Green's function method, the solution of~\eqref{Poisson} can be represented by
\begin{equation}
    u(\vx;g) =\int_{\Omega} G(\vx,\vx')g(\vx') \diff{\vx'},
\end{equation}
where for any $x'\in \Omega$, the Green's function $G(\vx,\vx')$ is the solution of  following problem,
\begin{equation}
    \begin{aligned}
         & -\Delta G(\vx,\vx')=\delta(\vx-\vx') \quad \vx\in \Omega
        \\
         & G(\vx,\vx')=0, \quad \vx\in \partial \Omega.
    \end{aligned}
\end{equation}
In the considered $2$D case, $\Omega=[0,1]^2$ and $g(x,y)=-a(x^2-x+y^2-y)$, we have
\begin{equation}
    u(x,y;g) = \int_0^1\int_0^1 G(x,y,x',y')g(x',y') \diff{x'}\diff{y'}.
\end{equation}
For demonstration, although we can obtain the analytical form of the Green's function, we use a DNN of hidden layer size $128$-$128$-$128$-$128$ $G_{\vtheta}(x,y,x',y')$ to fit the Green's function $G(x,y,x',y')$.

When we calculate the integral, it is impossible to integrate it analytically.
In practice, we often use the numerical integration.
We use the Gauss-Legendre quadrature.
Then we can represent neural operator $u_{\vtheta}(x,y;g)$ with Green's function DNN $G_{\vtheta}(x,y,x',y')$, that is,
\begin{equation}
    \begin{aligned}
        u_{\vtheta}(x,y;g) = \sum_{x^\prime \in S_{G_x}}\sum_{y^\prime \in S_{G_y}} \omega_{x^\prime} \omega_{y^\prime} G_{\vtheta}(x,y,x^\prime,y^\prime)g(x^\prime,y^\prime),
    \end{aligned}
    \label{Poiss_Gausianint}
\end{equation}
where  $S_{G_x}\subset[0,1]$,$S_{G_y}\subset[0,1]$ consist fixed integration points, determined by $1$D Gauss-Legendre quadrature and $\omega_{x^\prime}$, $\omega_{y^\prime}$ are corresponding coefficients.

\subsection{Empirical risk function}
For this toy example, to train the neural networks, we use no labeled data and only utilize the information of PDE, i.e., governing equation and boundary condition, for each $g^k,k=1,2,\cdots,K$.
To utilize the  constraint of governing equation of PDE, we uniformly sample a set of data from $\Omega=[0,1] \times [0,1]$, i.e., $S^{\Omega,k}$.
To utilize the information of boundary constraint, for each $k$, we uniformly sample a set of data from $\partial\Omega$.
Since the boundary $\partial\Omega$ consists of four line segments, i.e., $\partial\Omega=\bigcup_{i=1}^{4}\partial\Omega_i$, where $\partial\Omega_1=\{(0,y)\}_{y\in[0,1]}$,$\partial\Omega_2= \{(1,y)\}_{y\in[0,1]}$,$\partial\Omega_3=\{(x,0)\}_{x\in[0,1]}$,$\partial\Omega_4=\{(x,1)\}_{x\in[0,1]}$, we uniformly sample a set of data from $\partial\Omega_i$ respectively, i.e., $S^{\partial\Omega_i,k}$, i=1,2,3,4.

Since we use no labeled data, $\lambda_3$ in the general definition (\ref{lossfunctiondefinition}) is set as zero.
The empirical risk for this example is as follows,
\begin{equation}
    \begin{aligned}
        \RS
         & =  \frac{1}{K}\sum_{k\in[K]}  \bigg(\lambda_1
        \frac{1}{|S^{\Omega,k}|}\sum_{(x,y) \in S^{\Omega,k}}\Big(
        -\big(
        \partial_{xx}u_{\vtheta}(x,y;g^k)+\partial_{yy}u_{\vtheta}(x,y;g^k)\big)
        -g^k(x,y)
        \Big)^{2}                                                                                                                  \\
         & \quad +  \lambda_2 \frac{1}{|S^{\partial\Omega_1,k}|}\sum_{(x,y) \in S^{\partial\Omega_1,k}} {u_{\vtheta}(x,y;g^k)}^{2}
        \\
         & \quad +  \lambda_2 \frac{1}{|S^{\partial\Omega_2,k}|}\sum_{(x,y) \in S^{\partial\Omega_2,k}} {u_{\vtheta}(x,y;g^k)}^{2}
        \\
         & \quad +  \lambda_2\frac{1}{|S^{\partial\Omega_3,k}|}\sum_{(x,y) \in S^{\partial\Omega_3,k}} {u_{\vtheta}(x,y;g^k)}^{2}
        \\
         & \quad + \lambda_2\frac{1}{|S^{\partial\Omega_4,k}|}\sum_{(x,y) \in S^{\partial\Omega_4,k}} {u_{\vtheta}(x,y;g^k)}^{2}
        \bigg).
        \label{laplace_leastsquareloss}
    \end{aligned}
\end{equation}
Note that, the least square loss in (\ref{lossfunctiondefinition}) is not crucial. To support this point, we also use the variational loss used in Deep Ritz Method and the empirical risk is as follows,
\begin{equation}
    \begin{aligned}
        \RS
         & =  \frac{1}{K}\sum_{k\in[K]}  \bigg(\lambda_1
        \frac{1}{|S^{\Omega,k}|}\sum_{(x,y) \in  S^{\Omega,k}}\Big(\frac{1}{2}
        \big(
        |\partial_{x}u_{\vtheta}(x,y;g^k)|^2+|\partial_{y}u_{\vtheta}(x,y;g^k)|^2\big)
        -g^k(x,y)u_{\vtheta}(x,y;g^k)
        \Big)                                                                                                                      \\
         & \quad +  \lambda_2 \frac{1}{|S^{\partial\Omega_1,k}|}\sum_{(x,y) \in S^{\partial\Omega_1,k}} {u_{\vtheta}(x,y;g^k)}^{2}
        \\
         & \quad +  \lambda_2 \frac{1}{|S^{\partial\Omega_2,k}|}\sum_{(x,y) \in S^{\partial\Omega_2,k}} {u_{\vtheta}(x,y;g^k)}^{2}
        \\
         & \quad +  \lambda_2\frac{1}{|S^{\partial\Omega_3,k}|}\sum_{(x,y) \in S^{\partial\Omega_3,k}} {u_{\vtheta}(x,y;g^k)}^{2}
        \\
         & \quad + \lambda_2\frac{1}{|S^{\partial\Omega_4,k}|}\sum_{(x,y) \in S^{\partial\Omega_4,k}} {u_{\vtheta}(x,y;g^k)}^{2}
        \bigg).
        \label{laplace_deepritzloss}
    \end{aligned}
\end{equation}

\subsection{Learning process}
For each training epoch, we first randomly choose source functions $\{g^k\}_{k=1}^{K}$ and calculate their values on fixed  quadrature points$(x',y')$, where $x' \in S_{G_x}$ and $ y' \in S_{G_y}$.
Second, we randomly sample data and obtain data set $S^{\Omega,k},S^{\partial \Omega_i,k},i=1,2,3,4$.
We obtain the dataset $D= \{(x,y,x',y',g^k(x',y'))| (x,y)\in S^{\Omega,k} \cup (\bigcup_{i=1}^{4} S^{\partial\Omega_i,k}),x'\in S_{G_x},y' \in S_{G_y}\}$.
In the following, we feed the data into the neural network $G_{\vtheta}(x,y,x',y')$ and calculate the total risk~\eqref{laplace_leastsquareloss} or~\eqref{laplace_deepritzloss} with neural operator $u_{\vtheta}(x,y;g)$, see Eq. (\ref{Poiss_Gausianint}).
We train neural network $G_{\vtheta}$ with Adam to minimize the total risk, and finally we obtain a well-trained Green's function DNN $G_{\vtheta}$, furthermore, according to Eq.~\ref{Poiss_Gausianint}, we obtain a neural operator $u_{\vtheta}(x,y;g)$.

\subsection{Results}
The source function $g(x,y)=-a(x^2-x+y^2-y)$ is determined by $a$. We then denote source function $g(x,y)=-a(x^2-x+y^2-y)$ by $g_a$.
To illustrate our approach,
we train a neural operator mapping from $g$ to the solution of the Poisson equation~\eqref{2D_Poisson_case}.

\paragraph{Example 1: MOD-Net for a family of Poisson equations using least square loss}

For illustration, we would train the MOD-Net by various source functions and test the MOD-Net with several source functions that are not used for training. In this example, we calculate the empirical risk with least square loss.

\begin{figure}[!htb]
    \centering

    \subfigure[$a=15,x=0.0$]{
        \includegraphics[scale=0.11]{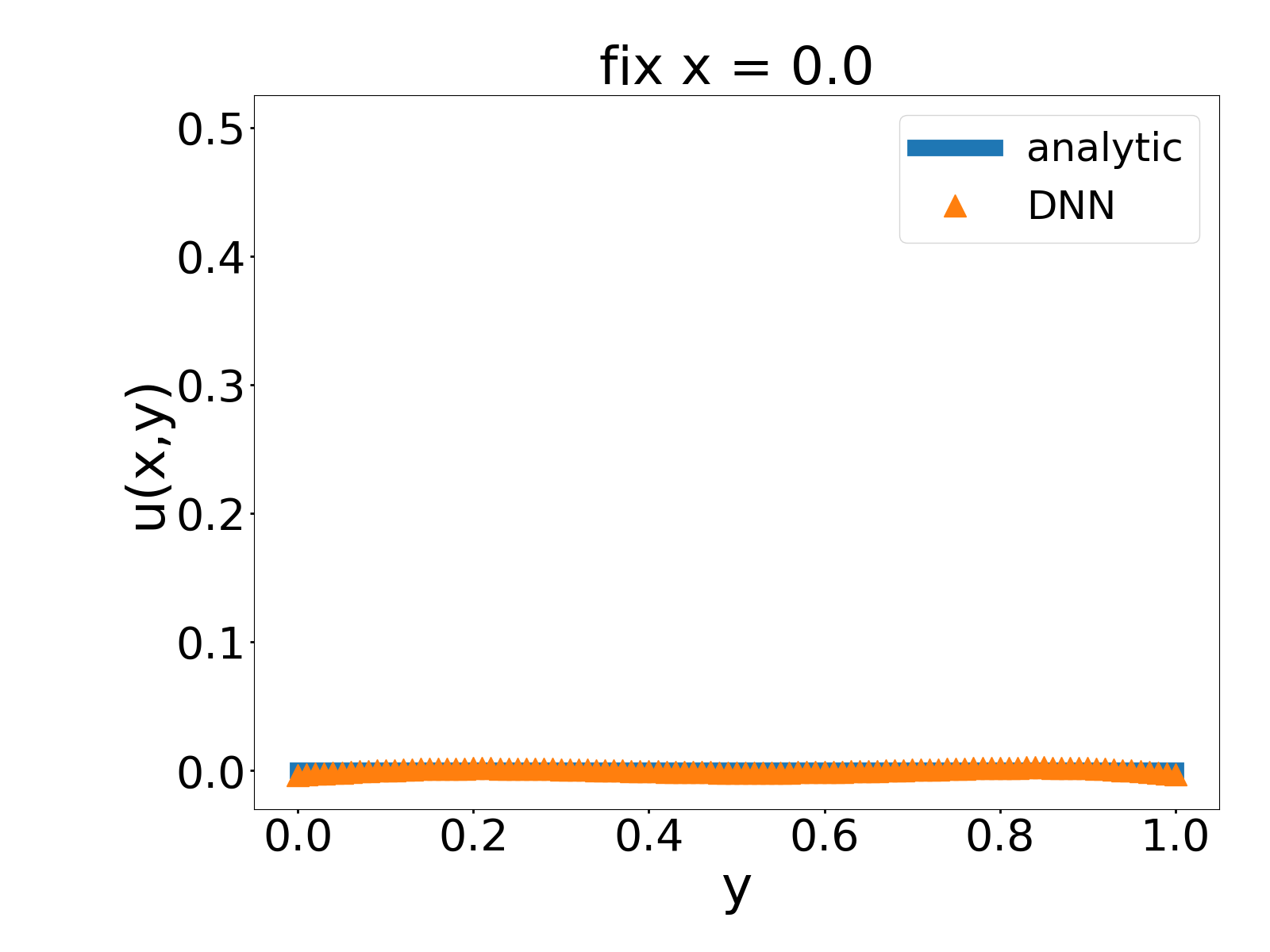}
    }
    \subfigure[$a=15,x=0.5$]{
        \includegraphics[scale=0.11]{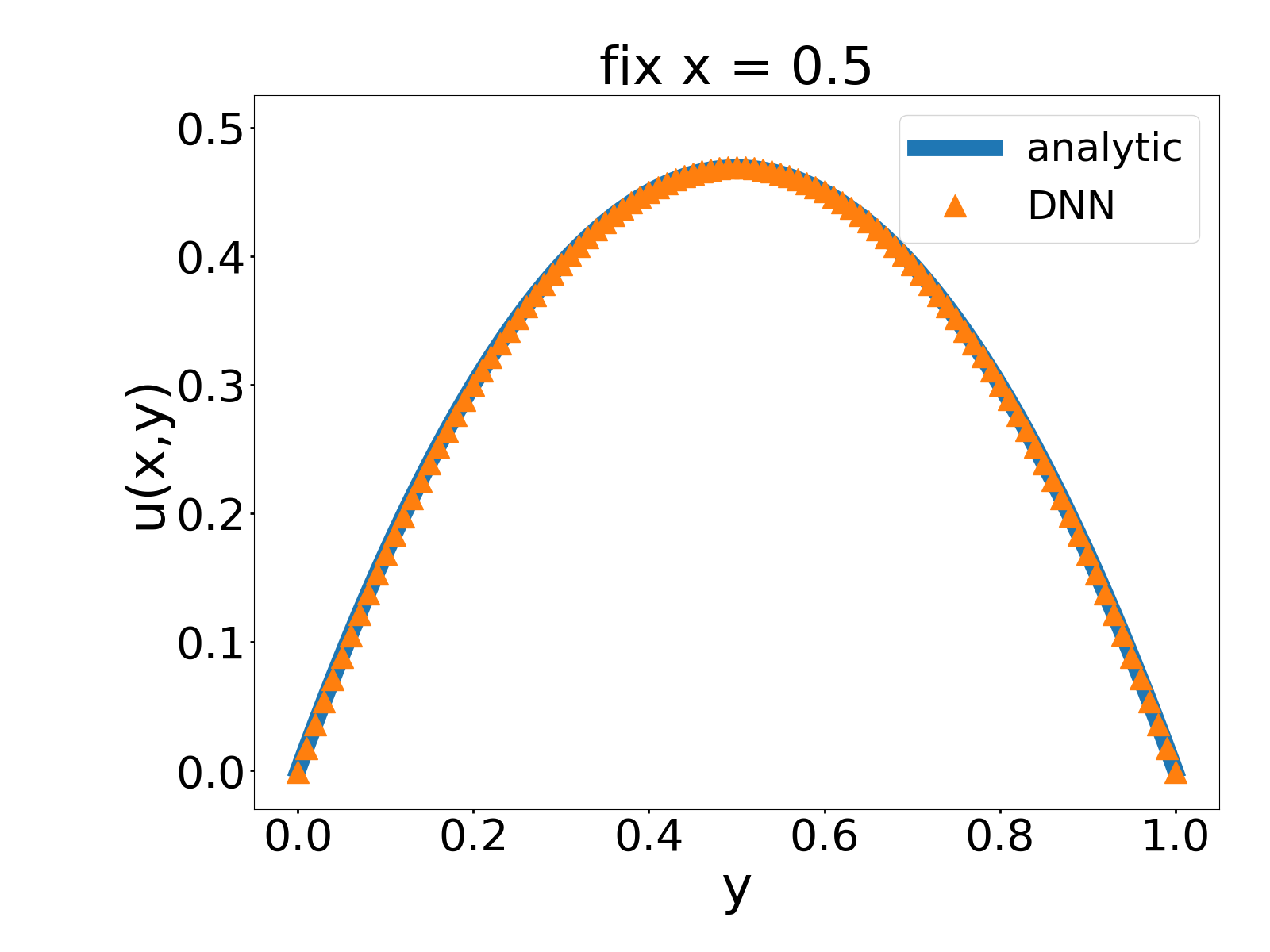}
    }
    \subfigure[$a=15,x=1.0$]{
        \includegraphics[scale=0.11]{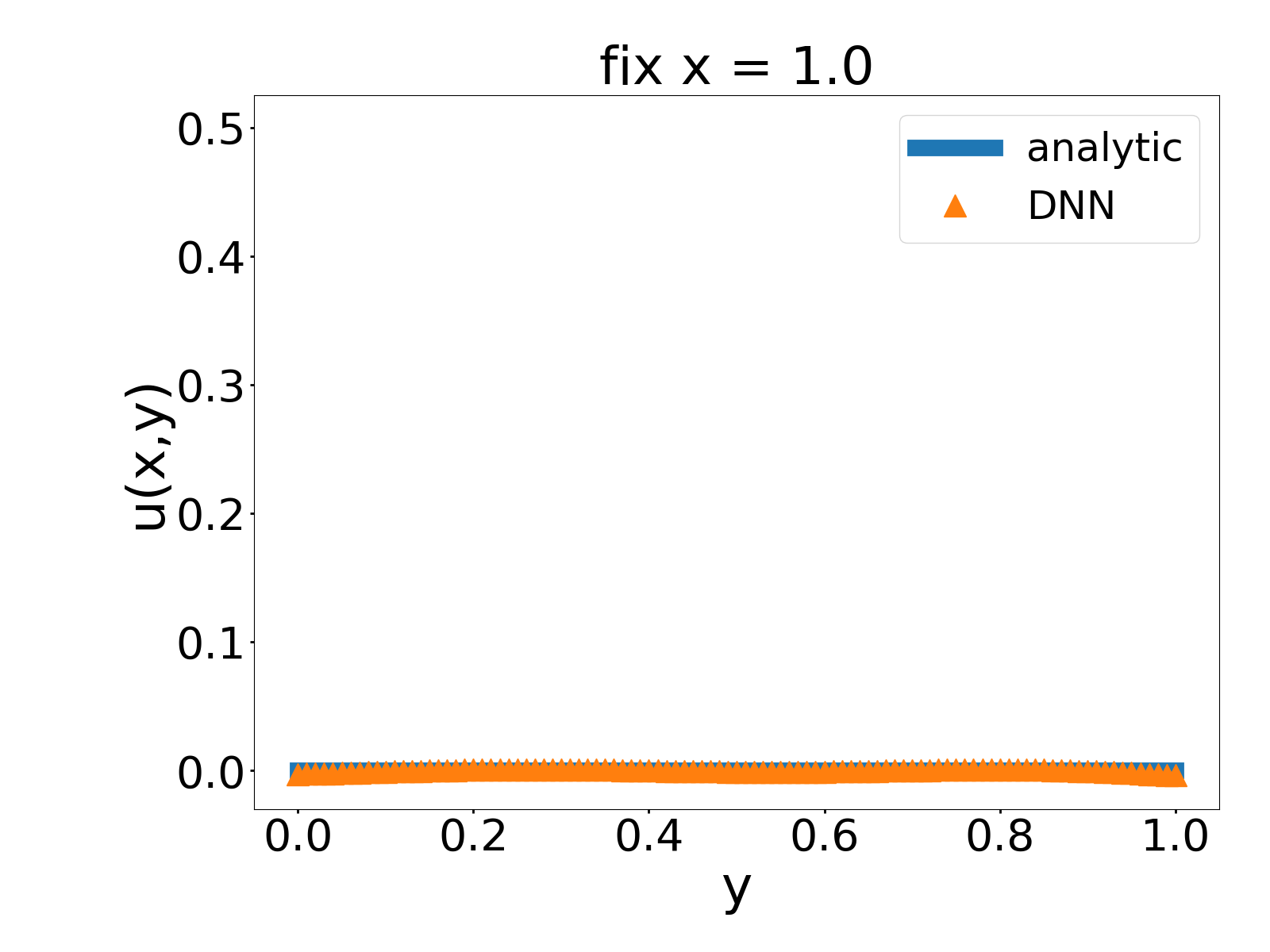}
    }
    \caption{Example 1. Comparison between analytic solution and  MOD-Net solution   on  $x=0,0.5,1$ corresponding to  source function determined by test $a=15$.}
    \label{fig:Poisson_op_contrast_leastsquare_section}
\end{figure}

During training, for each epoch, we choose a family of source functions, for example, we sample $K=10$ source functions $g_a$ from selected region, that is, we sample control parameter $a$ uniformly from $\{10k\}_{k=1}^{20}$.
We set the number of integration points $|S_{G_x}|=10$ and $|S_{G_y}|=10$, and for each epoch, we randomly sample $200$ points in $[0,1]^2$ and sample  $200$ points in each line of boundary, respectively. We set $\lambda_1=1,\lambda_2=1$ in the empirical risk using least square loss~\eqref{laplace_leastsquareloss}.
\begin{figure}
    \centering
    \subfigure[$a=15$]{
        \includegraphics[scale=0.11]{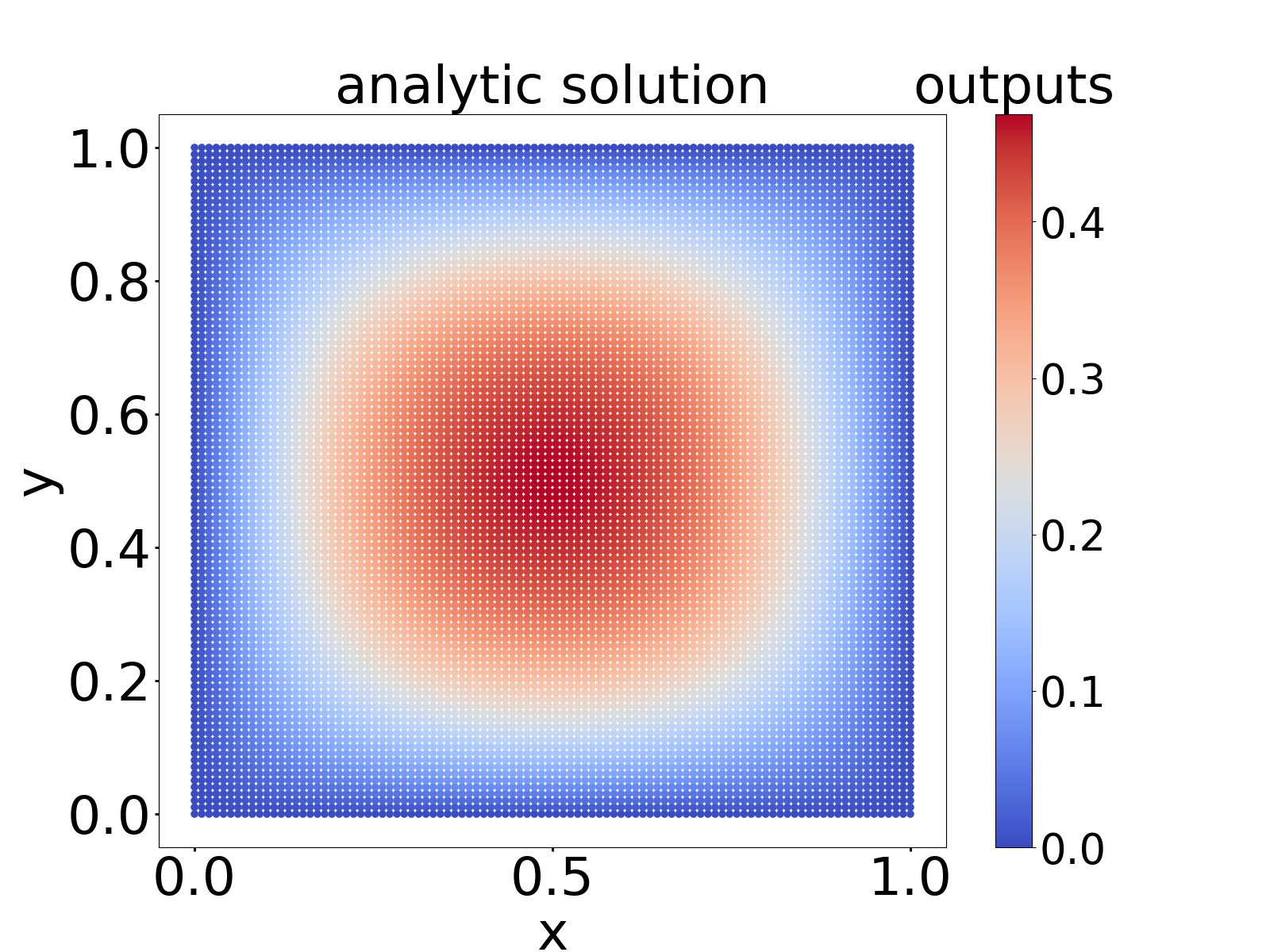}
    }
    \subfigure[$a=15$]{
        \includegraphics[scale=0.11]{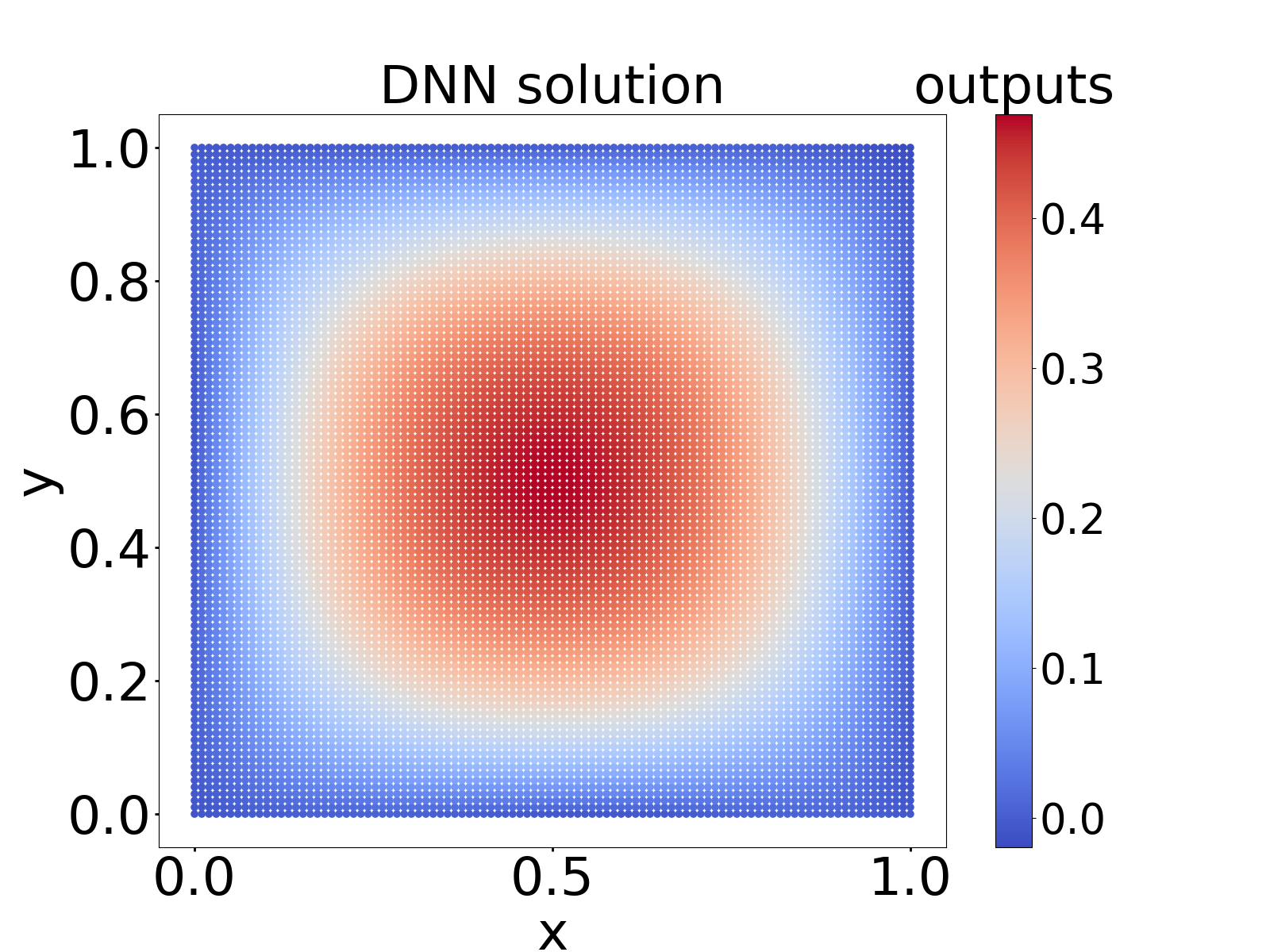}
    }
    \subfigure[$a=15$]{
        \includegraphics[scale=0.11]{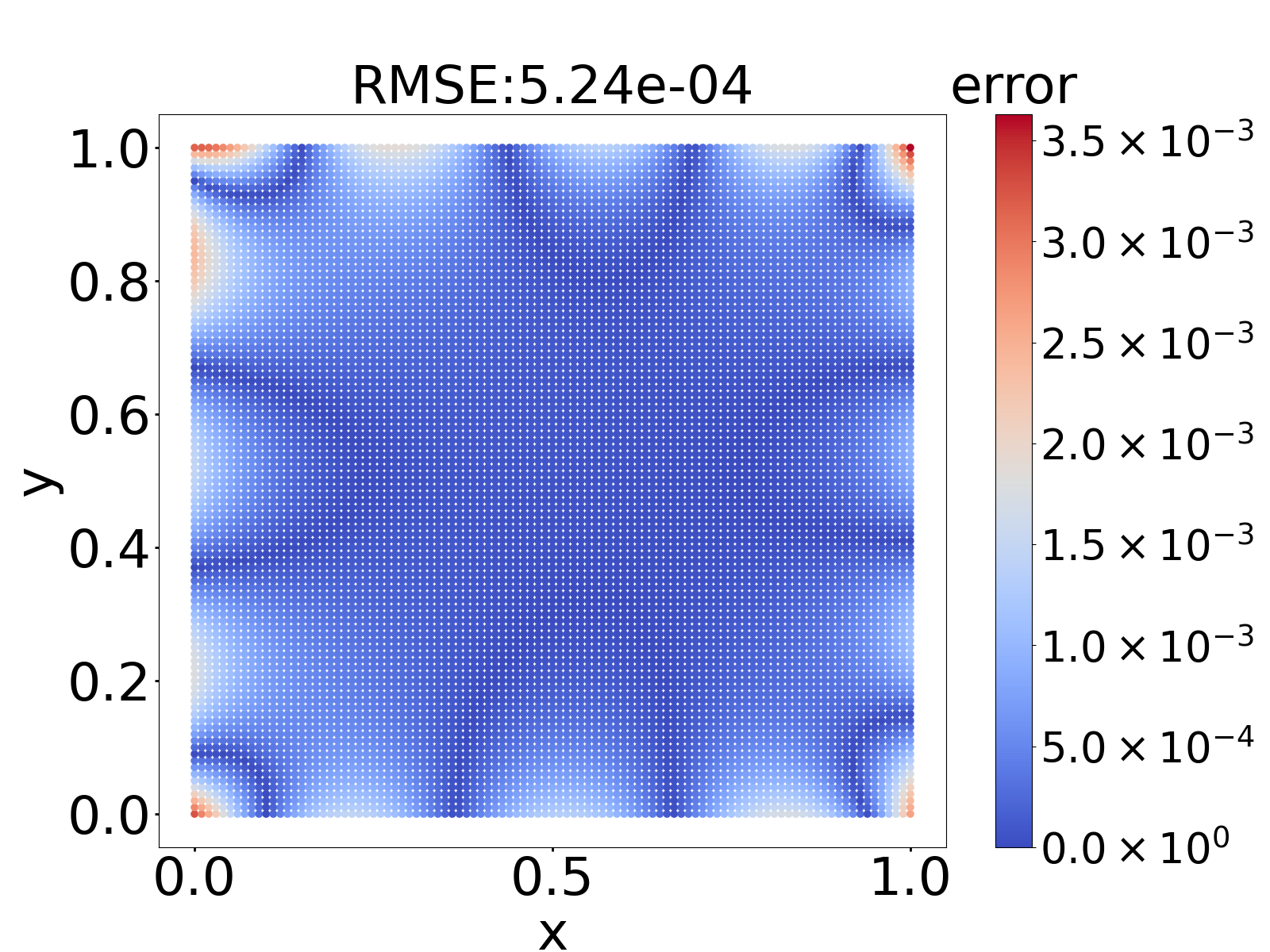}
    }
    \caption{Example 1. Comparison between analytic solution and  MOD-Net solution on $101\times101$ grid points corresponding to source terms determined by $a=15$. (a) Analytic solution. (b) MOD-Net solution. (c) The difference between the analytic solution and MOD-Net solution. Since the error at most points are very small, here we reduce the upper limit of the Color bar to see more information. }
    \label{fig:Poisson_op_contrast_leastsquare_error}
\end{figure}
We test the performance of this well-trained MOD-Net on $a=15$.
For each fixed $a$, the corresponding exact true solution is $u(x,y)=\frac{a}{2}x(x-1)y(y-1)$.
For visualizing the performance of our well-trained MOD-Net, we show the analytic solution and  MOD-Net solution on $101\times 101$ equidistributed isometric grid points by color in Fig. \ref{fig:Poisson_op_contrast_leastsquare_error}.
To compare these two solutions more intuitively, we calculate the difference of the two solutions at each point and show the error by color ,i.e., as the error increases, the color changes from blue to red, in Fig. \ref{fig:Poisson_op_contrast_leastsquare_error}(c).
The root of mean square error (RMSE) is $\sim 5.2\times 10^{-4}$.
We also show the solution obtained by MOD-Net and the corresponding analytical solution on fixed $x=0,0.5,1$ for  considered test source function.
As shown in  Fig. \ref{fig:Poisson_op_contrast_leastsquare_section}, the MOD-Net can well predict the analytic solutions.

\paragraph{Example 2: MOD-Net for a family of Poisson equations using variational loss}
In this experiment, we also train the MOD-Net by various source functions. But, note that, in this example, we calculate the empirical risk with variational loss used in Deep Ritz Method.

In training, for each epoch, similarly we set $K=10$ and we sample control parameter $a$ uniformly from $\{10k\}_{k=1}^{20}$.
We set the number of integration points $|S_{G_x}|=10$ and $|S_{G_y}|=10$, and for each epoch, we randomly sample $600$ points in $[0,1]^2$ and sample  $600$ points in each line of boundary, respectively. We set $\lambda_1=1,\lambda_2=250$ in the empirical risk using least square loss (\ref{laplace_deepritzloss}).
\begin{figure}
    \centering

    \subfigure[$a=15$]{
        \includegraphics[scale=0.11]{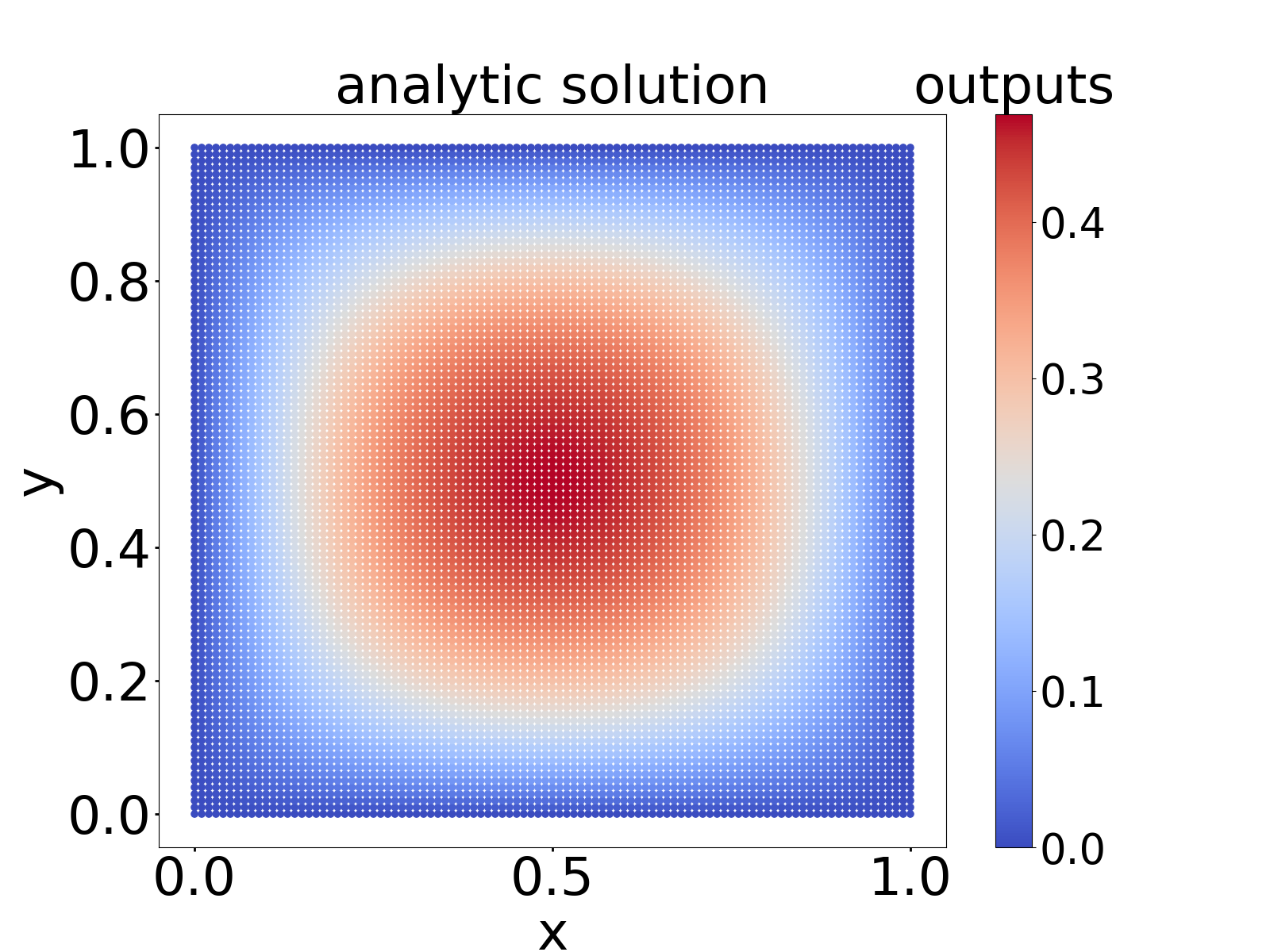}
    }
    \subfigure[$a=15$]{
        \includegraphics[scale=0.11]{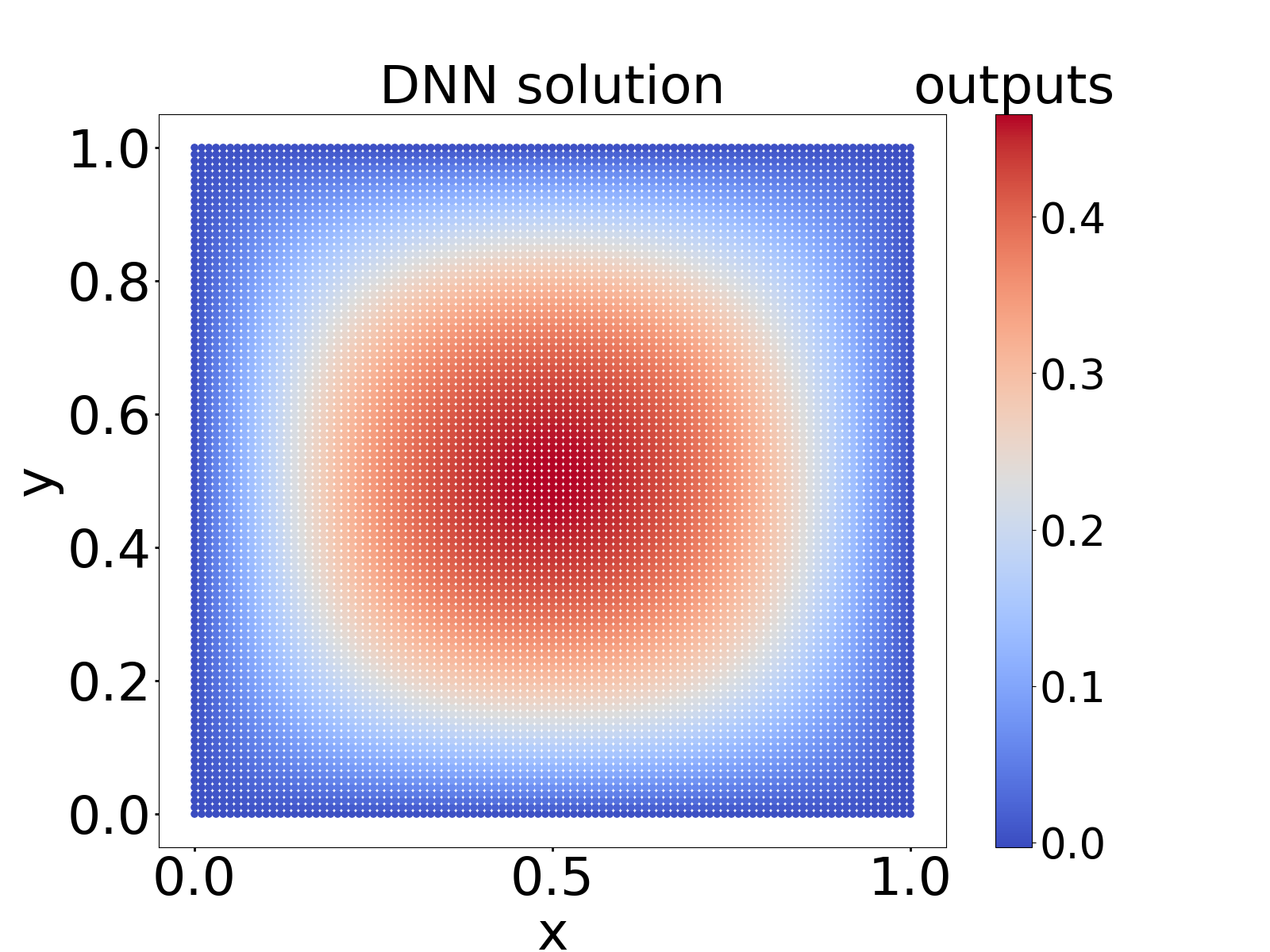}
    }
    \subfigure[$a=15$]{
        \includegraphics[scale=0.11]{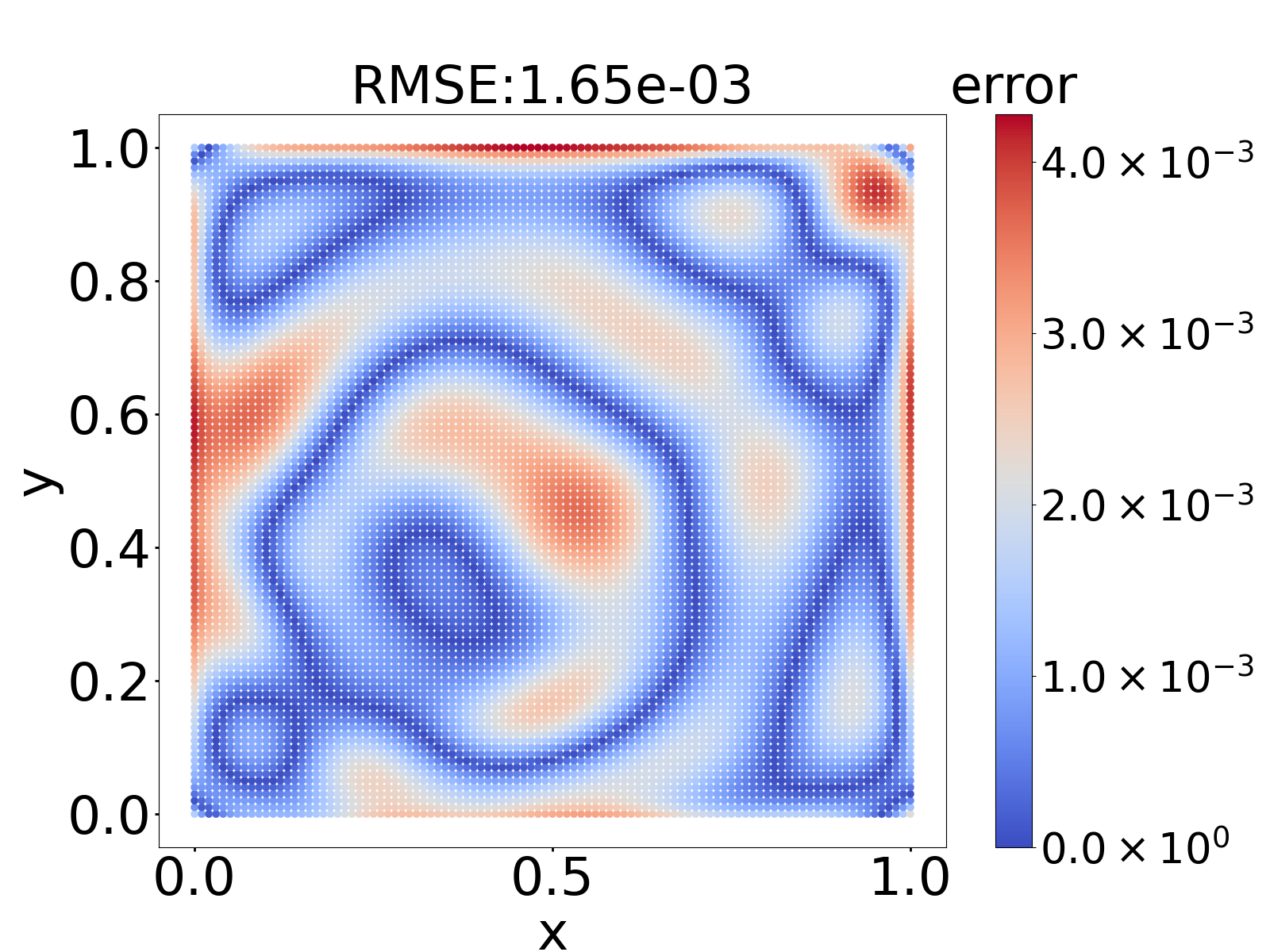}
    }
    \caption{Example 2. Comparison between analytic solution and  MOD-Net solution on $101\times101$ grid points corresponding to source terms determined by $a=15$. (a) Analytic solution. (b) MOD-Net solution. (c) The difference between the analytic solution and MOD-Net solution.}
    \label{fig:Poisson_op_contrast_variationalloss_error}
\end{figure}
We test the performance of this well-trained MOD-Net on $a=15$.
For each fixed $a$, the corresponding exact true solution is $u(x,y)=\frac{a}{2}x(x-1)y(y-1)$. For visualizing the performance of our well-trained MOD-Net, similarly we show the analytic solution and  MOD-Net solution on $101\times101$ equidistributed isometric grid points by color. To compare these two solutions more intuitively, similarly we calculate the difference of the two solutions at each point and show the error by color in Fig. \ref{fig:Poisson_op_contrast_variationalloss_error}. The root mean square error (RMSE) is $\sim1.6\times10^{-3}$.
\begin{figure}[!htb]
    \centering
    \subfigure[$a=15,x=0.0$]{
        \includegraphics[scale=0.11]{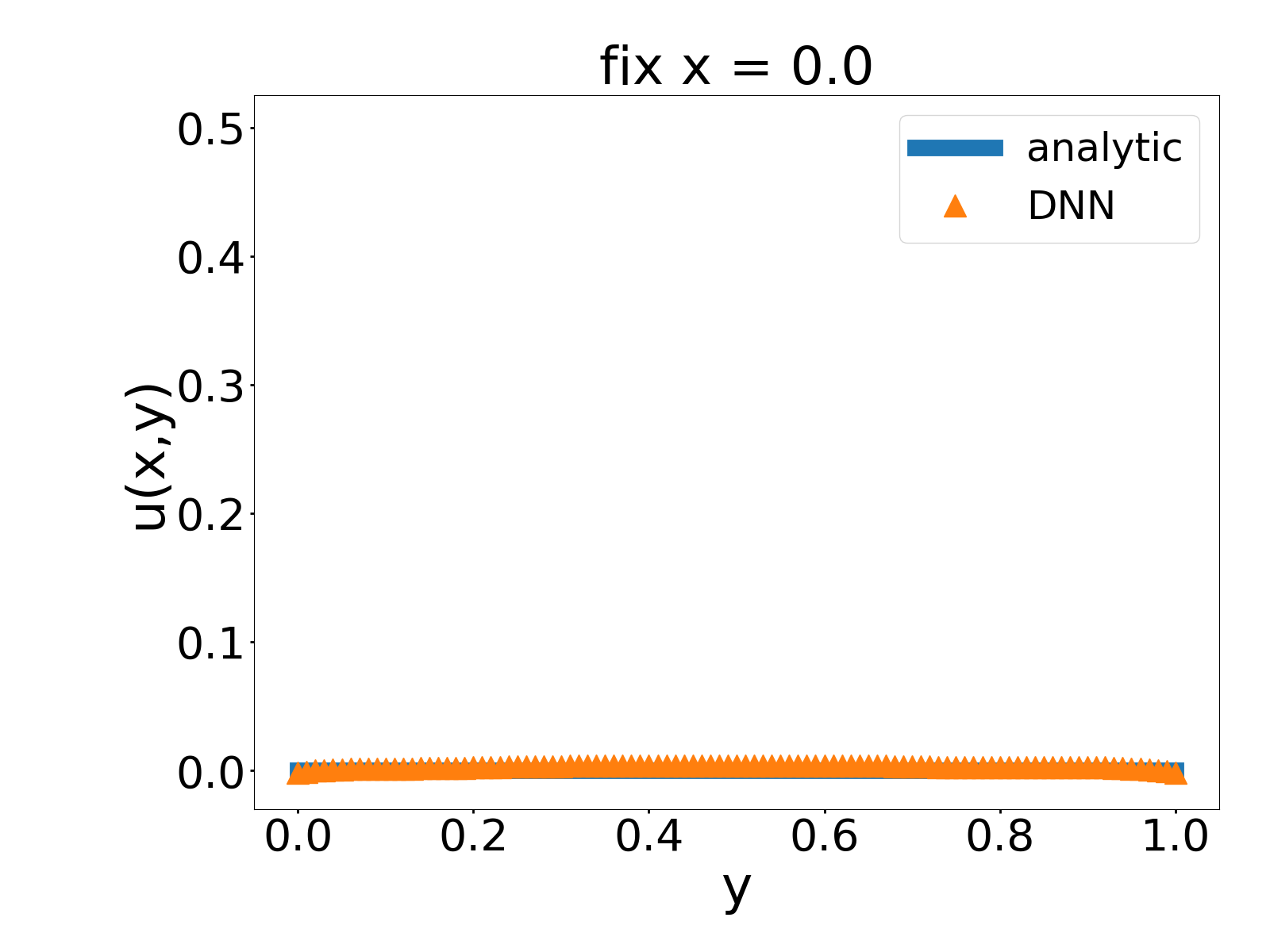}
    }
    \subfigure[$a=15,x=0.5$]{
        \includegraphics[scale=0.11]{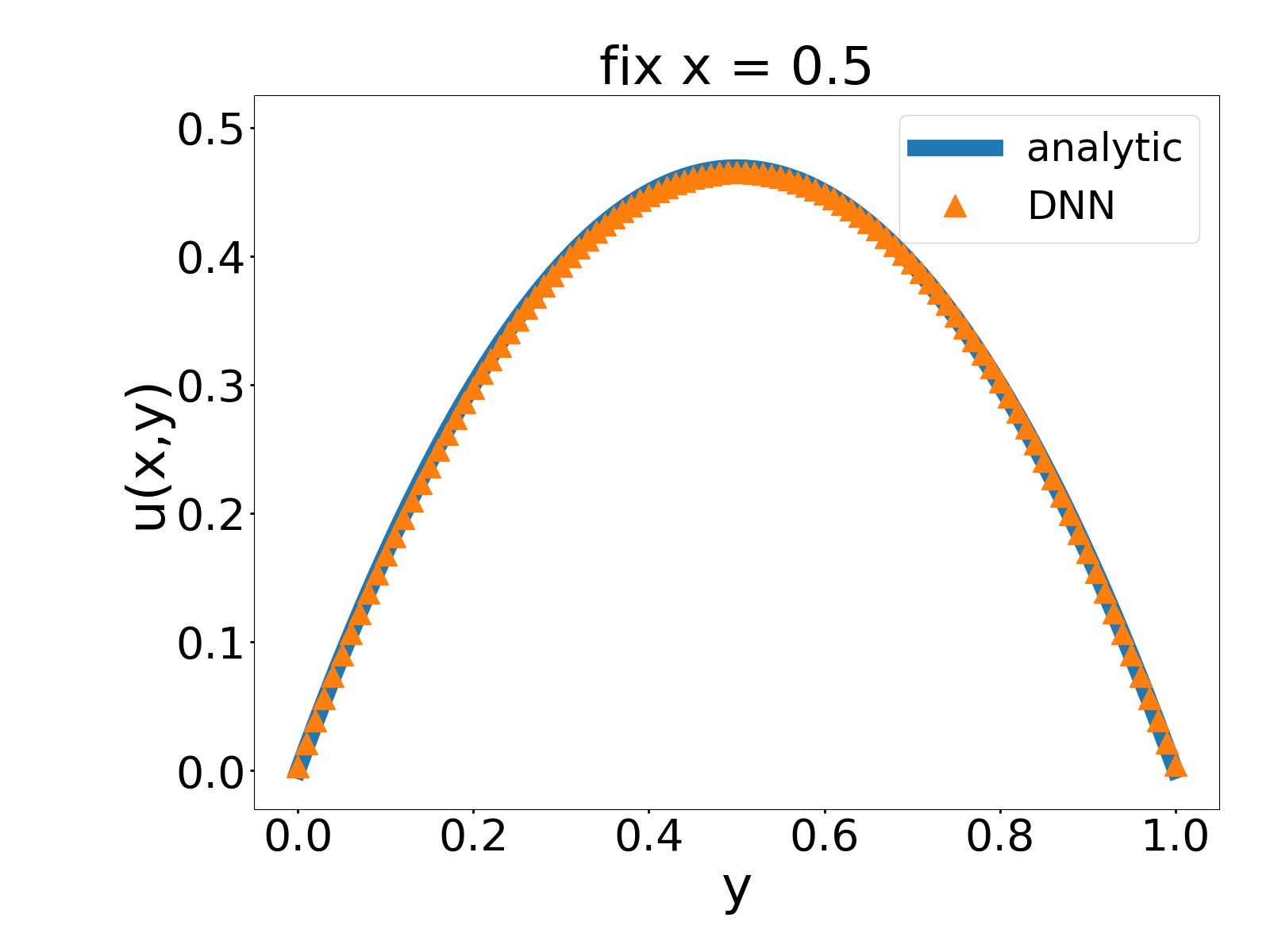}
    }
    \subfigure[$a=15,x=1.0$]{
        \includegraphics[scale=0.11]{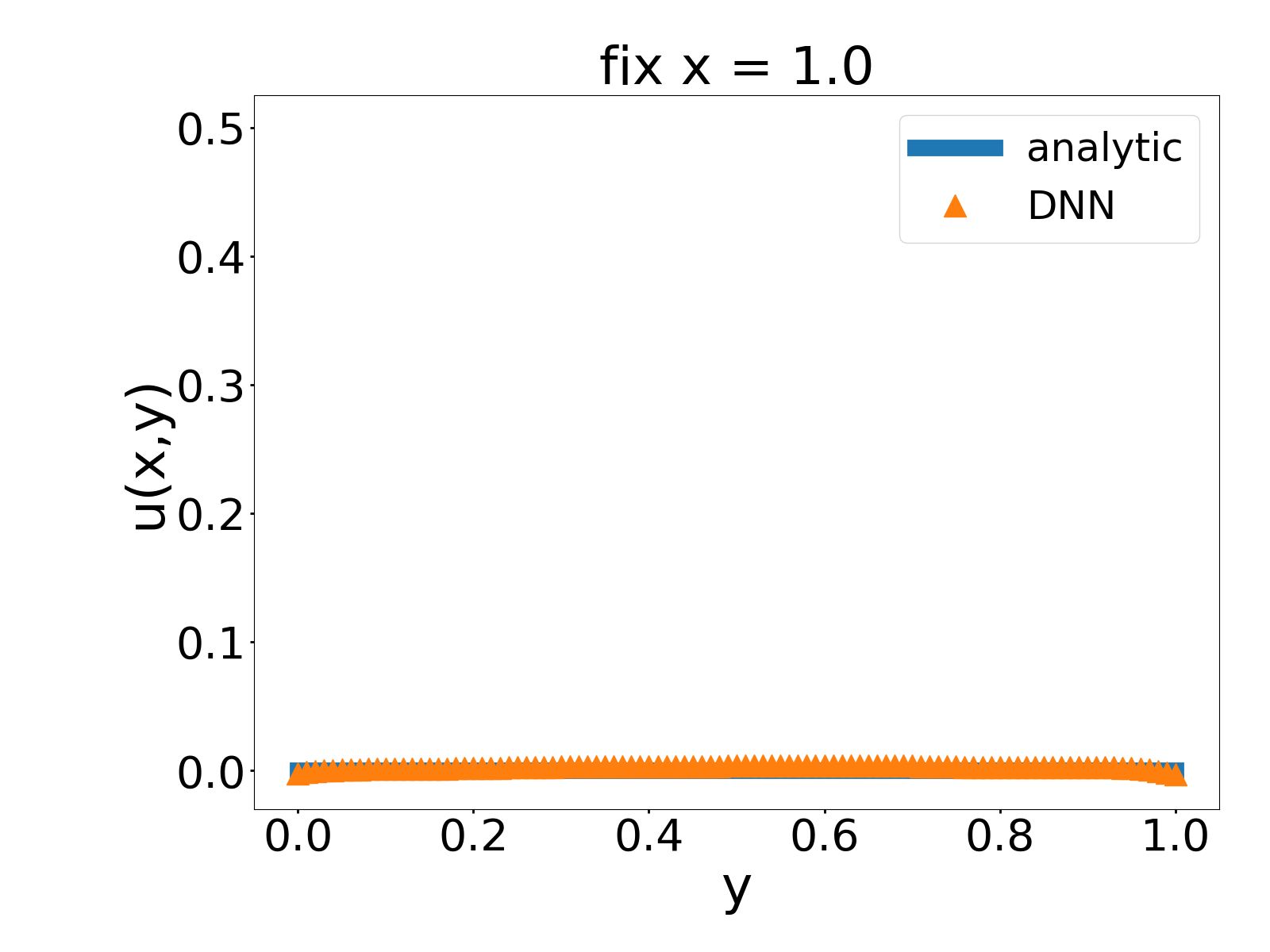}
    }
    \caption{Example 2. Comparison between analytic solution and  MOD-Net solution  on  $x=0,0.5,1$ corresponding to  source function determined by test $a=15$.}
    \label{fig:Poisson_op_contrast_variationalloss_section}
\end{figure}
For visualizing the performance of our well-trained MOD-Net, we show the solution obtained by MOD-Net and the corresponding analytical solution on fixed $x=0,0.5,1$ for  considered test source function. As shown in  Fig. \ref{fig:Poisson_op_contrast_variationalloss_section}, the MOD-Net can well predict the analytic solutions.

\section{Numerical experiments: equation with uncertainty or control variable} \label{sec:toyexample}
In this section, we would use a toy model to show the positive effect of the regularization from a few labels.
We consider the following equation,
\begin{equation}
    \begin{aligned}
         & \partial_{x}u(x,y) + a(x,y)u(x,y) = g(x,y), \quad (x,y)\in \Omega,
        \\
         & u(x,y)=0, \quad (x,y)\in \partial \Omega.
    \end{aligned}
    \label{constucted_equation}
\end{equation}
where $\Omega=[a,b]\times[c,d]$, the variable $y$ is the auxiliary variable that can be uncertainty variable or the control variable. $a(x,y)$ and $g(x,y)$ are coefficient functions that contain the randomness or the control functions in the control problem.
This equation can be regarded as a simplified toy model of a linear ODE with uncertainty or a simplified version of the equation with degeneracy such as the kinetic equations. We use the toy model here to show that data regularization is a crucial ingredients when dealing with this kind of equations.

\subsection{Use DNN to fit Green's function}
Similarly to Poisson equation, using the Green's function method, the solution of (\ref{constucted_equation})
can be represented by
\begin{equation}
    u(x,y;g) = \int_c^d\int_a^b G(x,y,x',y')g(x',y') \diff{x'}\diff{y'}.
\end{equation}
In this experiment, we use a  DNN of hidden layer size $128$-$128$-$128$-$128$ $G_{\vtheta}(x,y,x',y')$ equipped with activation function tanh  to fit the Green's function $G(x,y,x',y')$.

When we calculate the integral, we use the Gauss-Legendre quadrature. Then we can represent neural operator $u_{\vtheta}(x,y;g)$ with Green's function DNN $G_{\vtheta}(x,y,x',y')$, that is,
\begin{equation}
    \begin{aligned}
        u_{\vtheta}(x,y;g) = \sum_{x^\prime \in S_{G_x}}\sum_{y^\prime \in S_{G_y}} \omega_{x^\prime} \omega_{y^\prime} G_{\vtheta}(x,y,x^\prime,y^\prime)g(x^\prime,y^\prime),
    \end{aligned}
    \label{constructed_Gausianint}
\end{equation}
where  $S_{G_x}\subset[a,b]$,$S_{G_y}\subset[c,d]$ consist fixed integration points, determined by $1$D Gauss-Legendre quadrature and $\omega_{x^\prime}$, $\omega_{y^\prime}$ are corresponding coefficients.

\subsection{Empirical risk function}
For this example, to train the neural networks, we utilize the information of PDE, i.e., governing equation and boundary condition and a few data ,
$S^{u,k}=\{x_i,y_i,$ $u^{k}(x_i,y_i)\}_{i\in [n_k]}$ for each $g^k,k=1,2,\cdots,K$, where $u^{k}(\cdot)=u(\cdot;g^k)$.
Note that $S^{u,k}$ can be analytically solved on  grid points.
To utilize the  constraint of governing equation of PDE, we uniformly sample a set of data from $\Omega=[a,b] \times [c,d]$, i.e., $S^{\Omega,k}$. To utilize the information of boundary constraint, for each $k$, we uniformly sample a set of data from $\partial\Omega$.
Since the boundary $\partial\Omega$ consists of four line segments, i.e., $\partial\Omega=\bigcup_{i=1}^{4}\partial\Omega_i$, where $\partial\Omega_1=\{(a,y)\}_{y\in[c,d]}$,$\partial\Omega_2= \{(b,y)\}_{y\in[c,d]}$,$\partial\Omega_3=\{(x,c)\}_{x\in[a,b]}$,$\partial\Omega_4=\{(x,d)\}_{x\in[a,b]}$, we uniformly sample a set of data from $\partial\Omega_i$ respectively, i.e., $S^{\partial\Omega_i,k}$, i=1,2,3,4.

The empirical risk for this example is as follows,
\begin{equation}
    \begin{aligned}
        \RS
         & =  \frac{1}{K}\sum_{k\in[K]}  \bigg(\lambda_1
        \frac{1}{|S^{\Omega,k}|}\sum_{(x,y) \in S^{\Omega,k}}\Big(
        \partial_{x}u_{\vtheta}(x,y;g^k) + a(x,y)u_{\vtheta}(x,y;g^k)  -g^k(x,y)
        \Big)^{2}                                                                                                                  \\
         & \quad +  \lambda_2 \frac{1}{|S^{\partial\Omega_1,k}|}\sum_{(x,y) \in S^{\partial\Omega_1,k}} {u_{\vtheta}(x,y;g^k)}^{2}
        \\
         & \quad +  \lambda_2 \frac{1}{|S^{\partial\Omega_2,k}|}\sum_{(x,y) \in S^{\partial\Omega_2,k}} {u_{\vtheta}(x,y;g^k)}^{2}
        \\
         & \quad +  \lambda_2\frac{1}{|S^{\partial\Omega_3,k}|}\sum_{(x,y) \in S^{\partial\Omega_3,k}} {u_{\vtheta}(x,y;g^k)}^{2}
        \\
         & \quad + \lambda_2\frac{1}{|S^{\partial\Omega_4,k}|}\sum_{(x,y) \in S^{\partial\Omega_4,k}} {u_{\vtheta}(x,y;g^k)}^{2}
        \\
         & \quad + \lambda_3 \frac{1}{n_k} \sum_{i\in [n_k]} |u_{\vtheta}(x_i,y_i;g^k)-u^{k}(x_i,y_i)|^2
        \bigg).
        \label{constructed_leastsquareloss}
    \end{aligned}
\end{equation}

\subsection{Learning process}
Similarly to Poisson equation, for each training epoch, we first randomly choose source functions $\{g^k\}_{k=1}^{K}$ and calculate their values on fixed  integration points$(x',y')$, where $x' \in S_{G_x}$ and $ y' \in S_{G_y}$. Second, we randomly sample data and obtain data set $S^{\Omega,k},S^{\partial \Omega_i,k},i=1,2,3,4$.
We obtain the data set $D= \{(x,y,x',y',g^k(x',y'))| (x,y)\in S^{\Omega,k} \cup (\bigcup_{i=1}^{4} S^{\partial\Omega_i,k}),x'\in S_{G_x},y' \in S_{G_y}\}$. In the following, we feed the data into the neural network $G_{\vtheta}(x,y,x',y')$ and calculate the total risk (\ref{constructed_leastsquareloss})  with neural operator $u_{\vtheta}(x,y;g)$, see Eq. (\ref{Poiss_Gausianint}).
Train neural network $G_{\vtheta}$ with Adam to minimize the total risk, and finally we  obtain a well-trained Green's function DNN $G_{\vtheta}$, furthermore, according to Eq. (\ref{Poiss_Gausianint}), we obtain a neural operator $u_{\vtheta}(x,y;g)$.

\subsection{Results}
To demonstrate data regularization, we apply MOD-Net method in the following example.
\paragraph{Example 3: MOD-Net for equation with auxiliary variable}
\begin{figure}[!htb]
    \centering
    \includegraphics[scale=0.11]{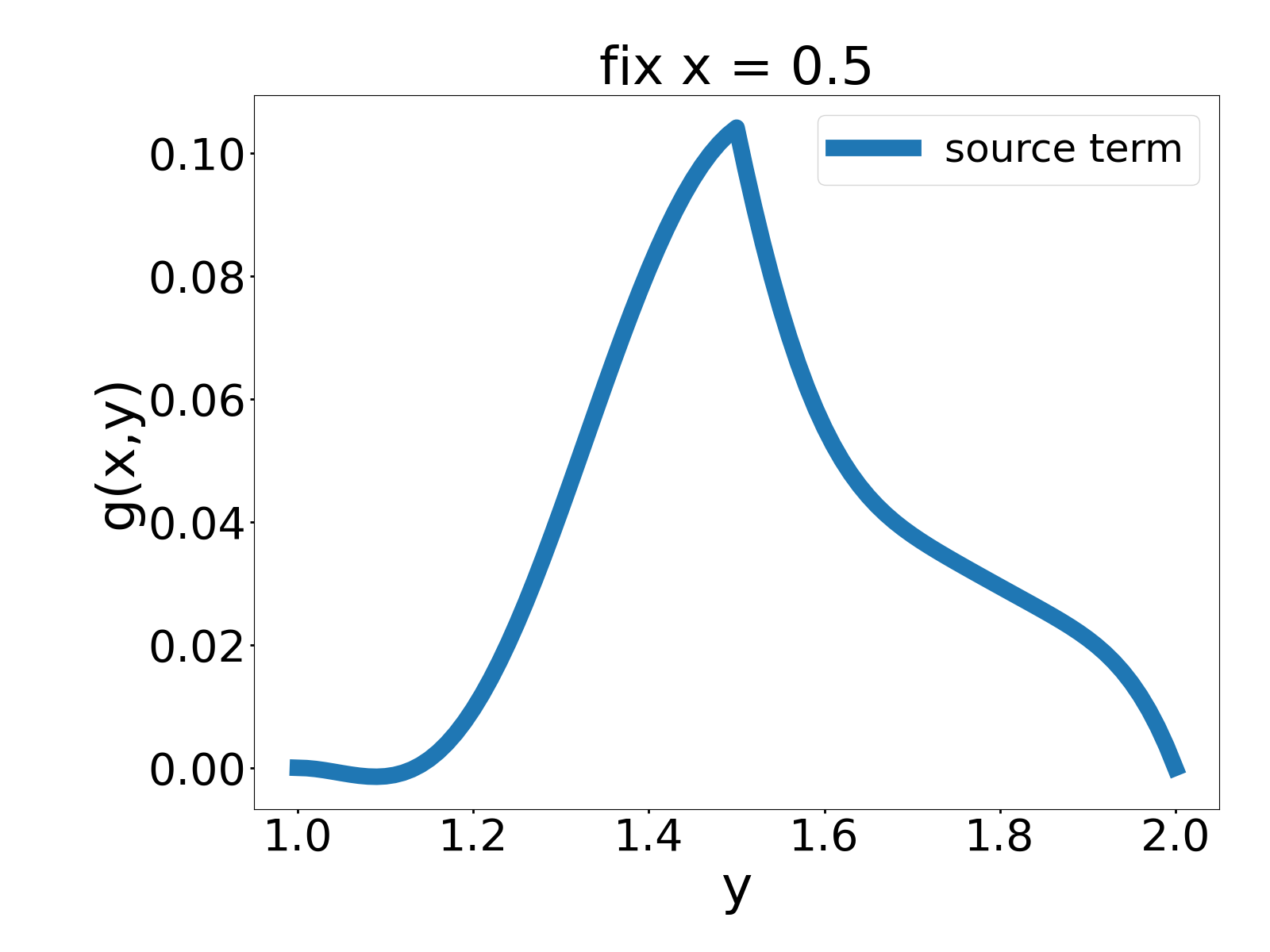}
    \caption{Example 3. Source term $g(x,y)$ on $x=0.5$, i.e., $g(0.5,y)$ vs. $y$. Remark that $\partial_{y} g$ is discontinuous.}
    \label{fig:example5_sourceterm}
\end{figure}
In this example, we consider the equation (\ref{constucted_equation}) in which $\Omega=[0,1] \times [1,2]$, $a(x,y)=\frac{\sin(4\pi xy)+2}{y}$ and the source term
\begin{equation}
    g(x, y)=\left\{
    \begin{aligned}
        (y^2-3y+2) \Bigg(x^2-x +\Big(2x-1 + \frac{x(x-1)(2+\sin(4\pi xy))}{y}\Big)(x+\cos(\pi y))   \Bigg),
        \quad 1\leq y<1.5, \\
        (y^2-3y+2) \Bigg(x^2-x +\Big(2x-1 + \frac{x(x-1)(2+\sin(4\pi xy))}{y}\Big)(x+\sin(2\pi y))   \Bigg), \quad  1.5\leq x\leq 2.
    \end{aligned}\right.
\end{equation}
The graph of source term is shown in Fig. \ref{fig:example5_sourceterm}.
For this equation, the analytical solution
\begin{equation}
    u(x,y) =
    \left\{
    \begin{aligned}
        x(x-1)(y-1)(y-2) \Big(x+\cos(\pi y)\Big), \quad 1\leq y<1.5, \\
        x(x-1)(y-1)(y-2) \Big( x+\sin(2\pi y)\Big), \quad  1.5\leq x\leq 2.
    \end{aligned}\right.
\end{equation}
Note that the first derivative of solution $u(x,y)$ with respect to $x$, i.e., $\partial_{x}u$ is continuous, but  $\partial_{y}u$ is discontinuous.
For this example, we train the MOD-Net to approximate the solution of a specific equation with three different loss functions.

\begin{figure}[!htb]
    \centering
    \subfigure[Training loss]{
        \includegraphics[scale=0.11]{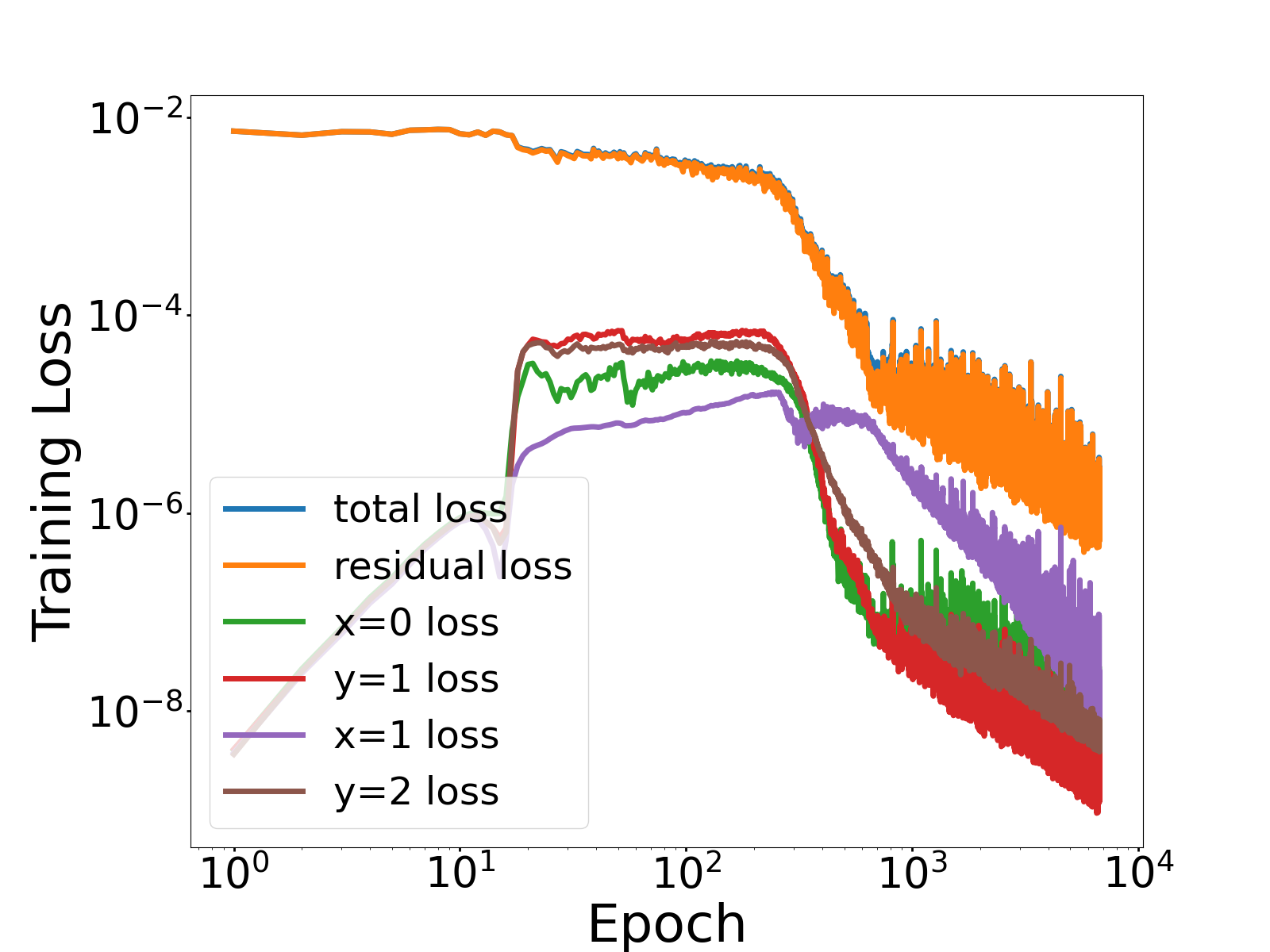}
    }
    \subfigure[$x=0.5$]{
        \includegraphics[scale=0.11]{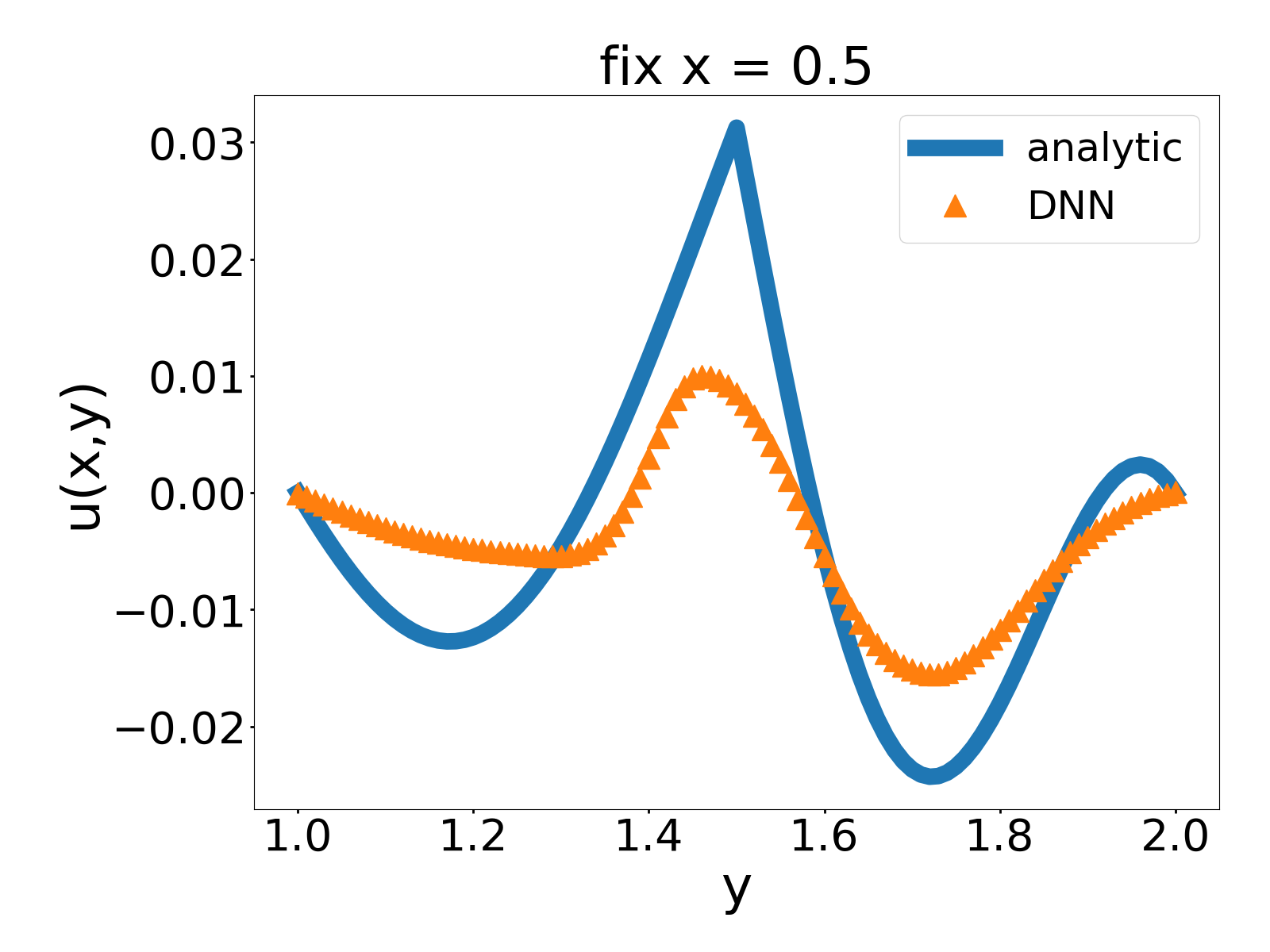}
    }
    \subfigure[$y=1.55$]{
        \includegraphics[scale=0.11]{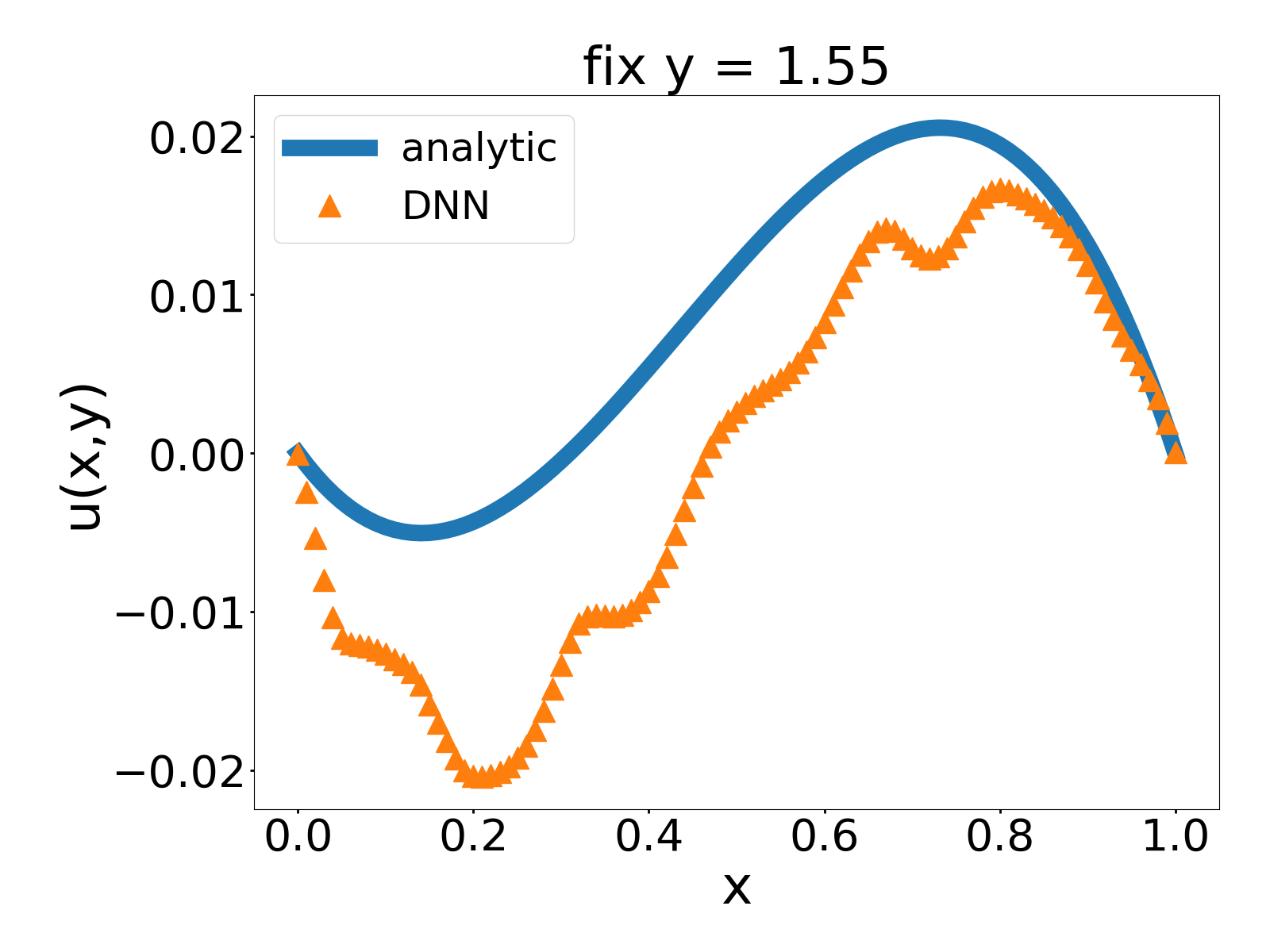}
    }
    \subfigure[Analytic solution]{
        \includegraphics[scale=0.11]{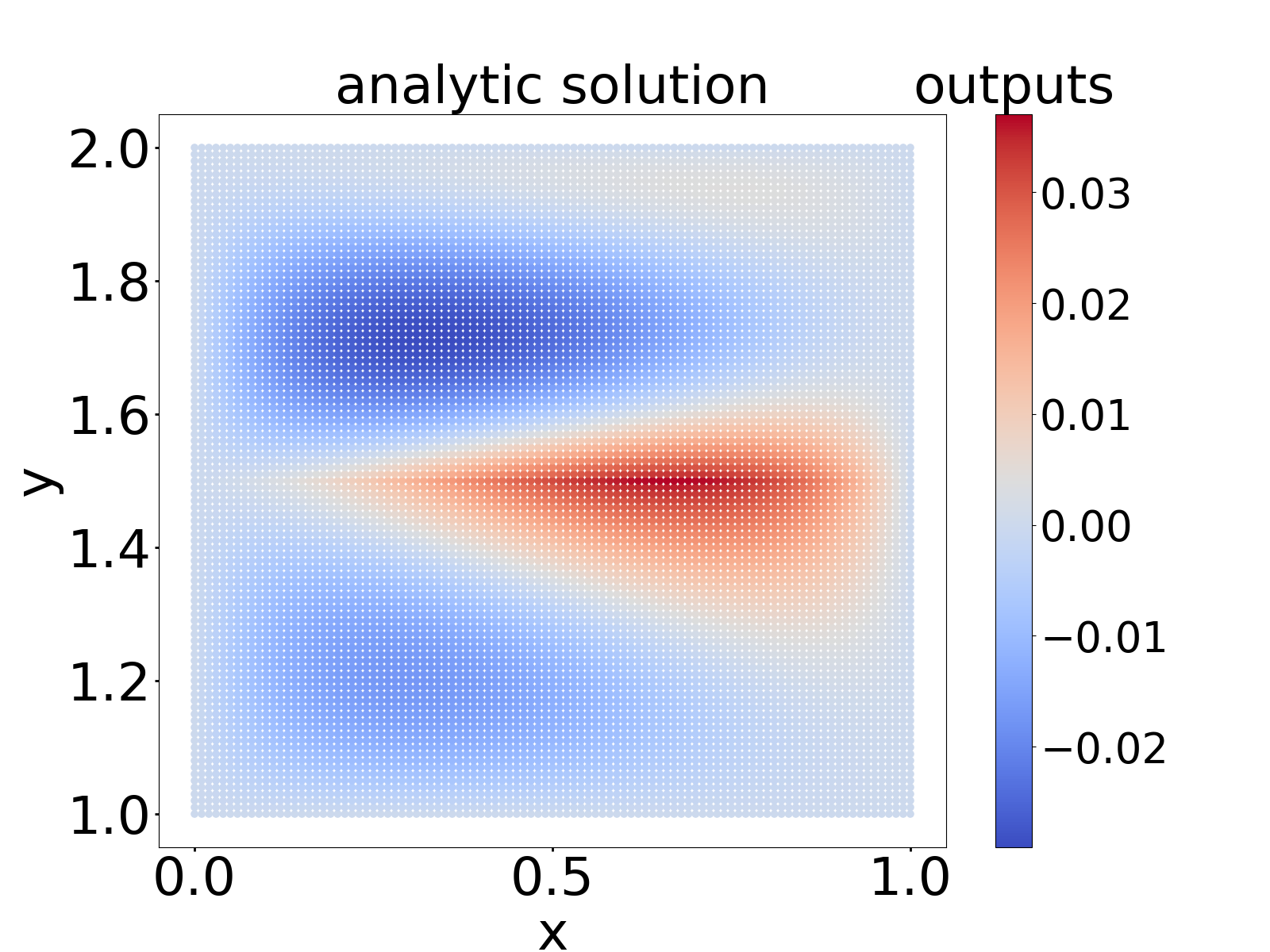}
    }
    \subfigure[MOD-Net solution]{
        \includegraphics[scale=0.11]{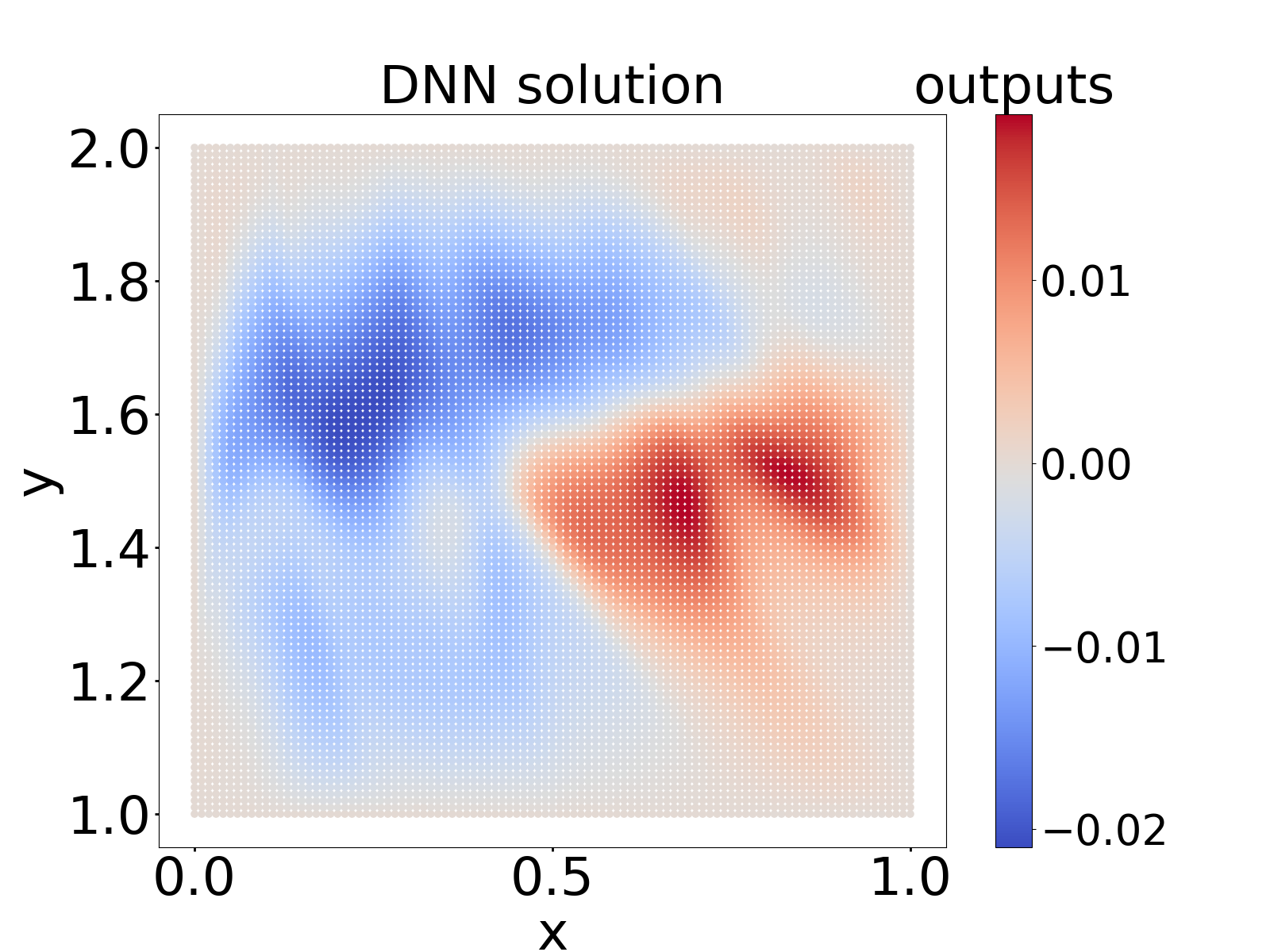}
    }
    \subfigure[Difference between two solutions]{
        \includegraphics[scale=0.11]{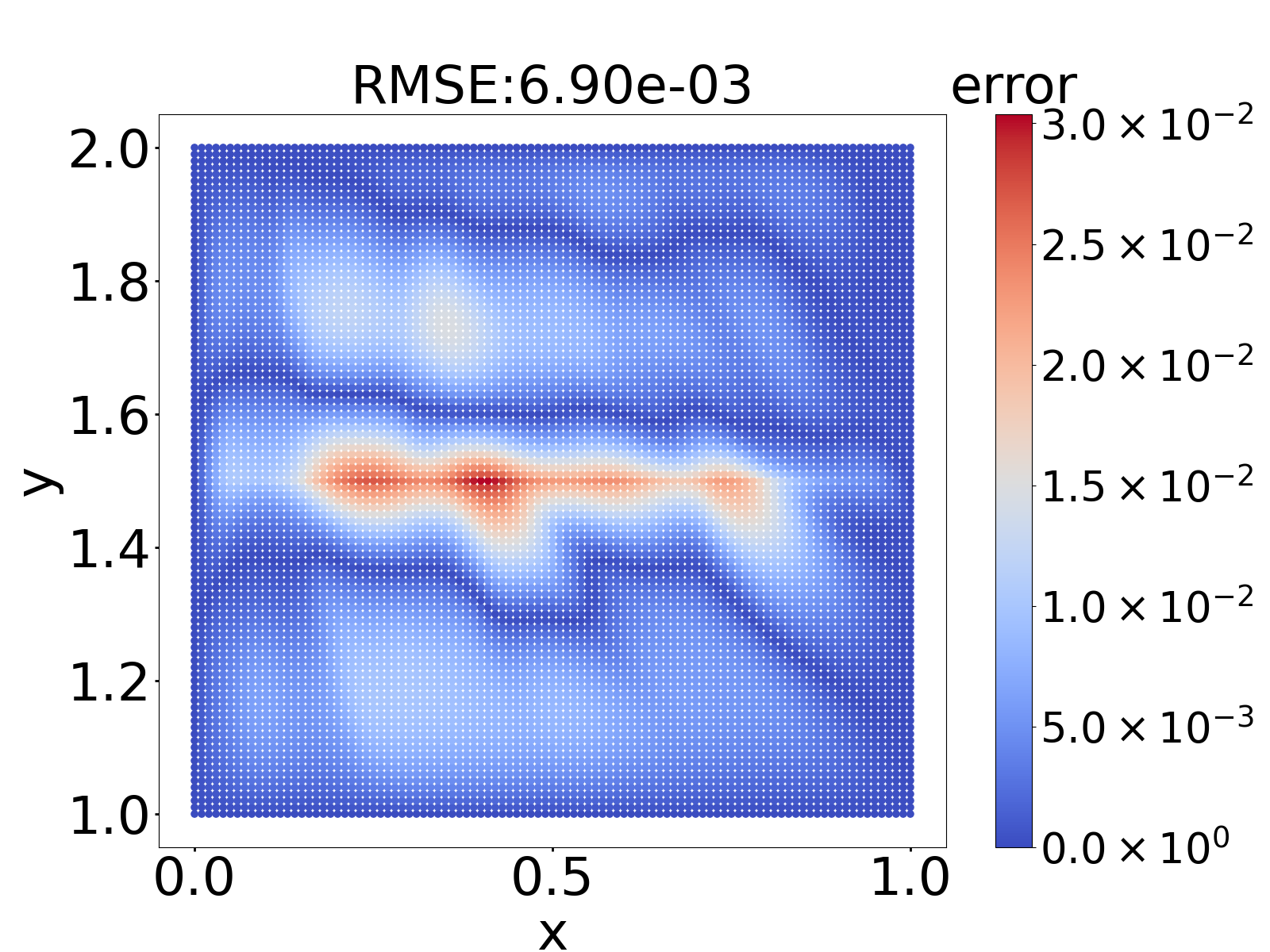}
    }
    \caption{Example 3. Trained by only the PDE,i.e., the governing equation and boundary condition. (a) Training loss. The  blue curve is the total risk Eq. (\ref{constructed_leastsquareloss}). The other five curves represent five terms in the total risk Eq.  (\ref{constructed_leastsquareloss}) except the last term, respectively.
        (b,c)
        Comparison between analytic solution and MOD-Net solution on $x=0.5$ and $y=1.55$.
        (d,e,f)
        Comparison between analytic solution and  MOD-Net solution on $101\times101$ grid points. }
    \label{fig:example5_constructed_onlyPDE}
\end{figure}

First, we set the $\lambda_3$ be zero in (\ref{constructed_leastsquareloss}),
and train the MOD-Net
only by the governing equation and boundary conditions. The empirical risk, i.e., training loss is shown in Fig. \ref{fig:example5_constructed_onlyPDE} (a).
For visualizing the performance of our well-trained MOD-Net $u_{\vtheta}$ intuitively, we plot the MOD-Net solution and the corresponding analytic solution on $x=0.5$ and $y=1.55$. As shown in Fig. \ref{fig:example5_constructed_onlyPDE} (b,c), the MOD-Net solution is not matching with the analytic solution.
We also show the analytic solution and  MOD-Net solution on $101\times101$ equidistributed isometric grid points by color. To compare these two solutions more intuitively, we calculate the difference of the two solutions at each point and show the error by color ,i.e., as the error increases, the color changes from blue to red, in Fig. \ref{fig:example5_constructed_onlyPDE} (d,e,f). The root mean square error (RMSE) is $\sim 6.90\times 10^{-3}$ and relative error is $\sim 57.61\%$.

\begin{figure}[!htb]
    \centering

    \subfigure[Training loss]{
        \includegraphics[scale=0.11]{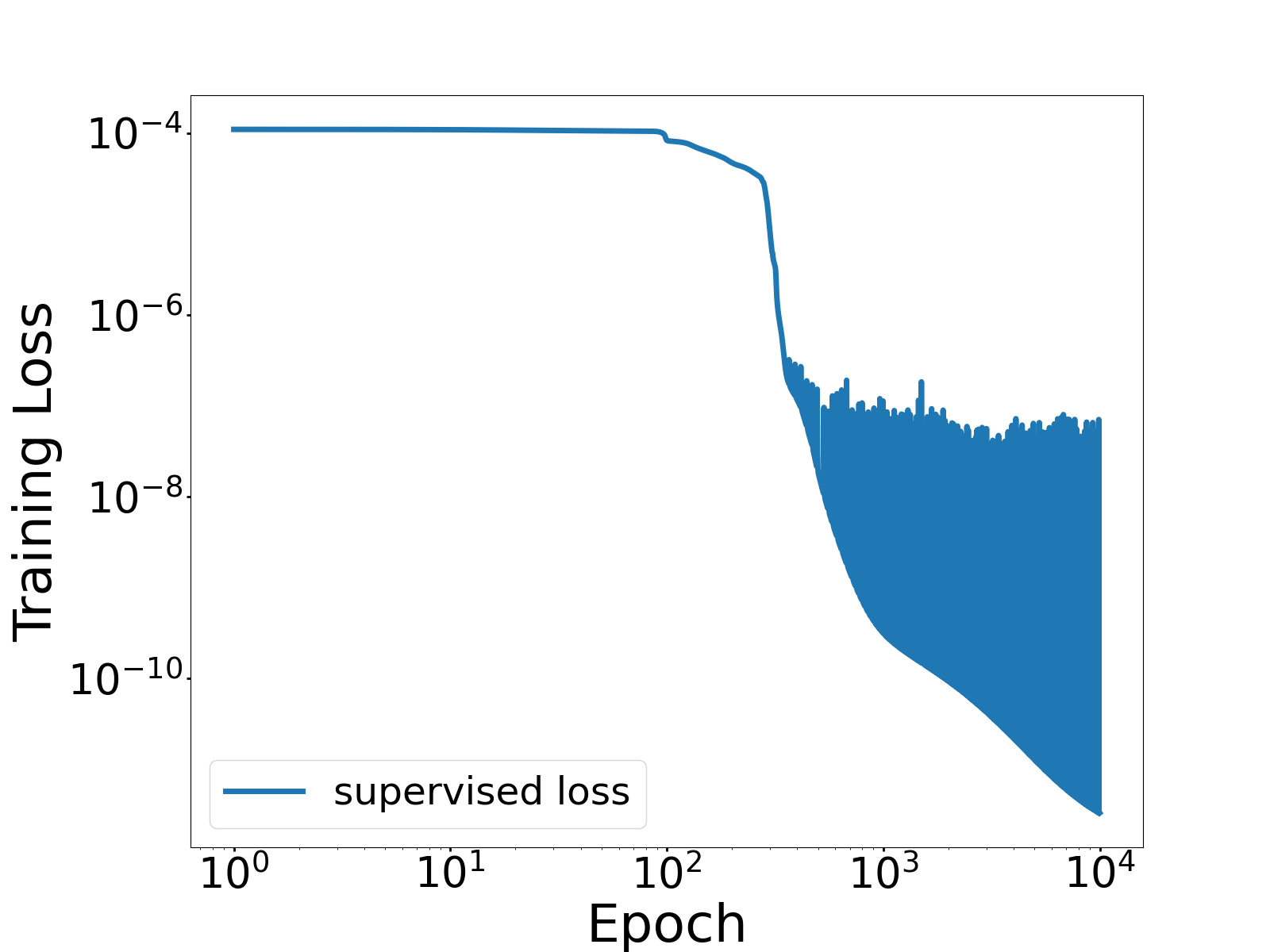}
    }
    \subfigure[$x=0.5$]{
        \includegraphics[scale=0.11]{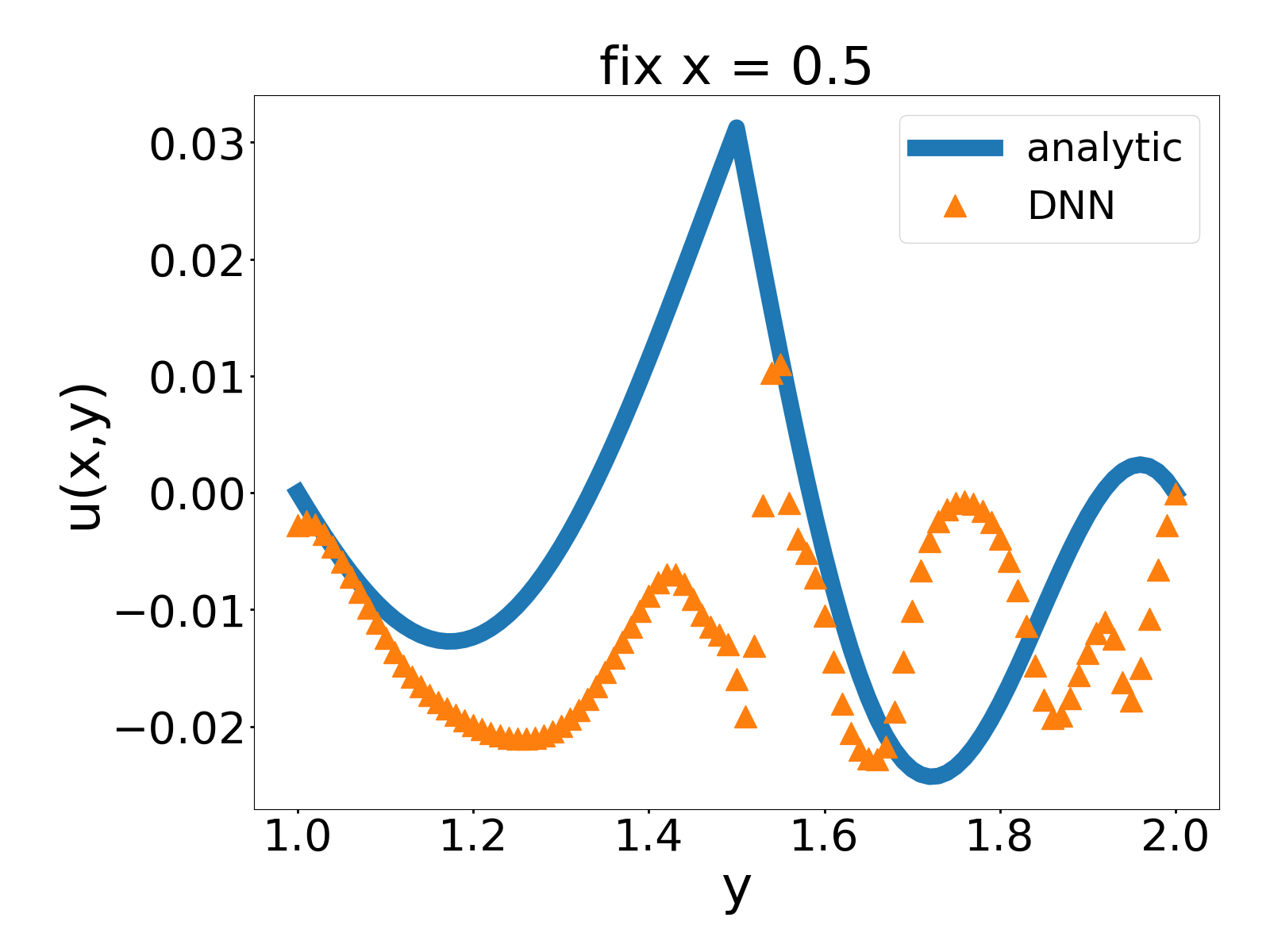}
    }
    \subfigure[$y=1.55$]{
        \includegraphics[scale=0.11]{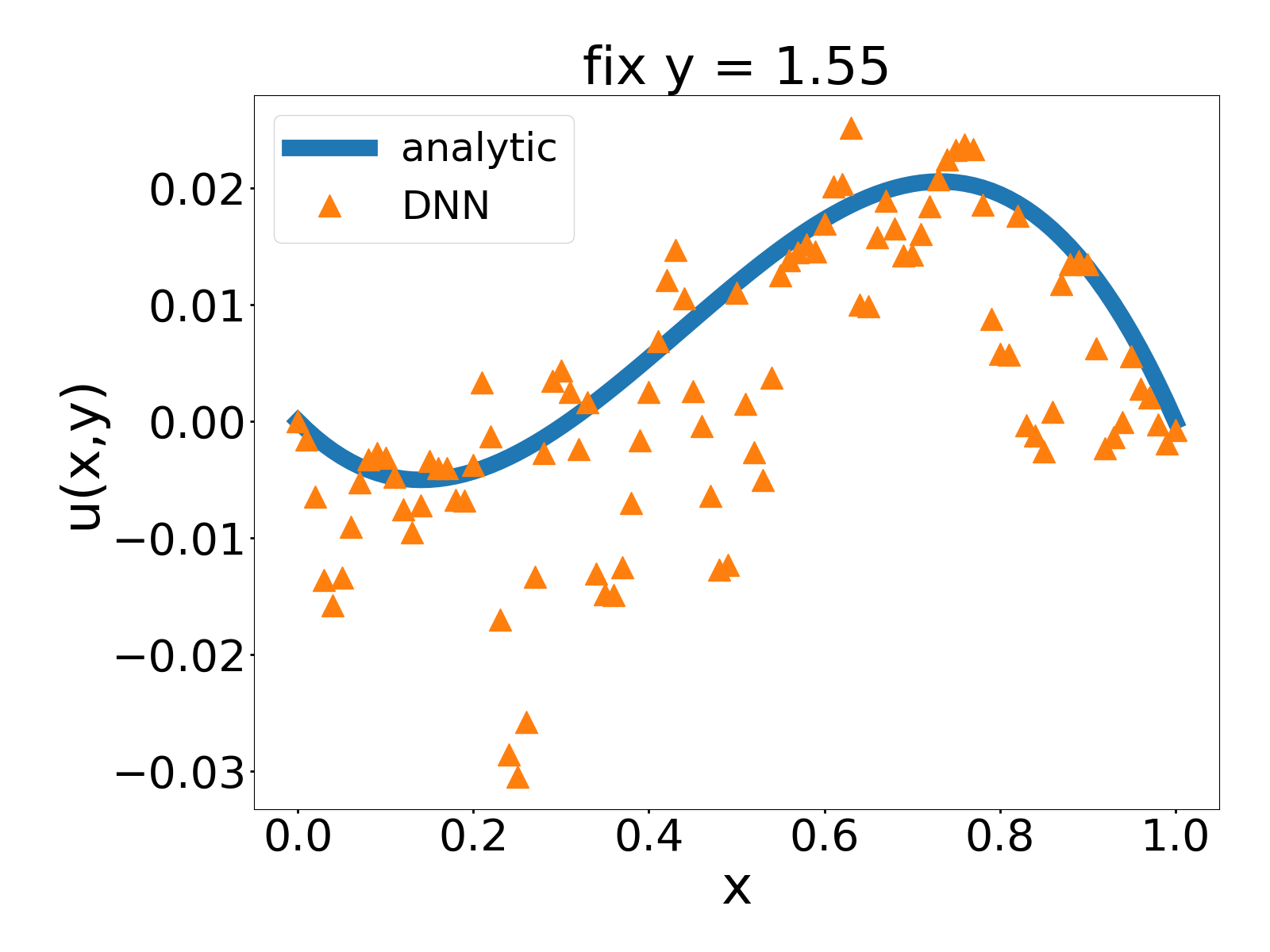}
    }
    \subfigure[Analytic solution]{
        \includegraphics[scale=0.11]{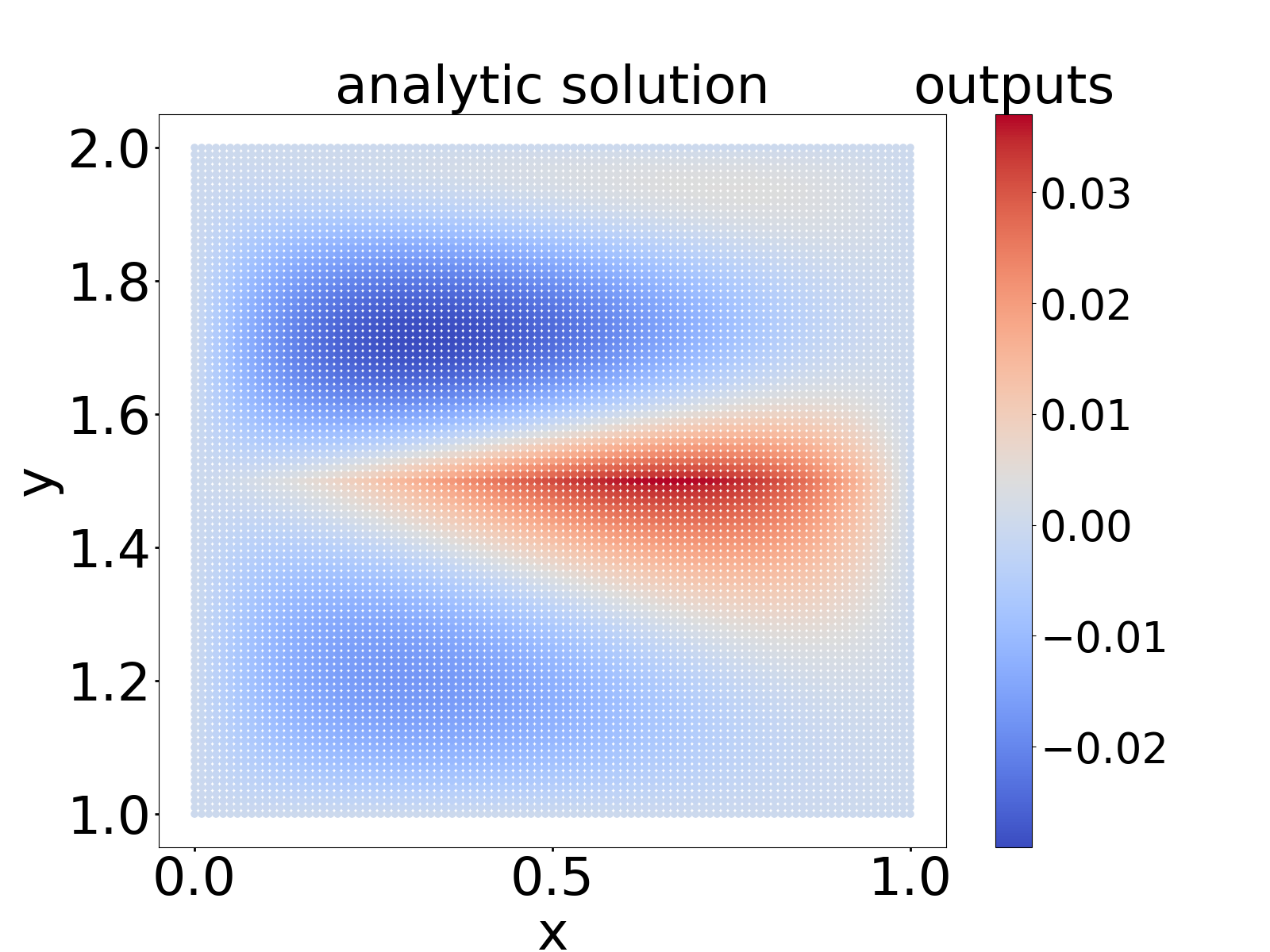}
    }
    \subfigure[MOD-Net solution]{
        \includegraphics[scale=0.11]{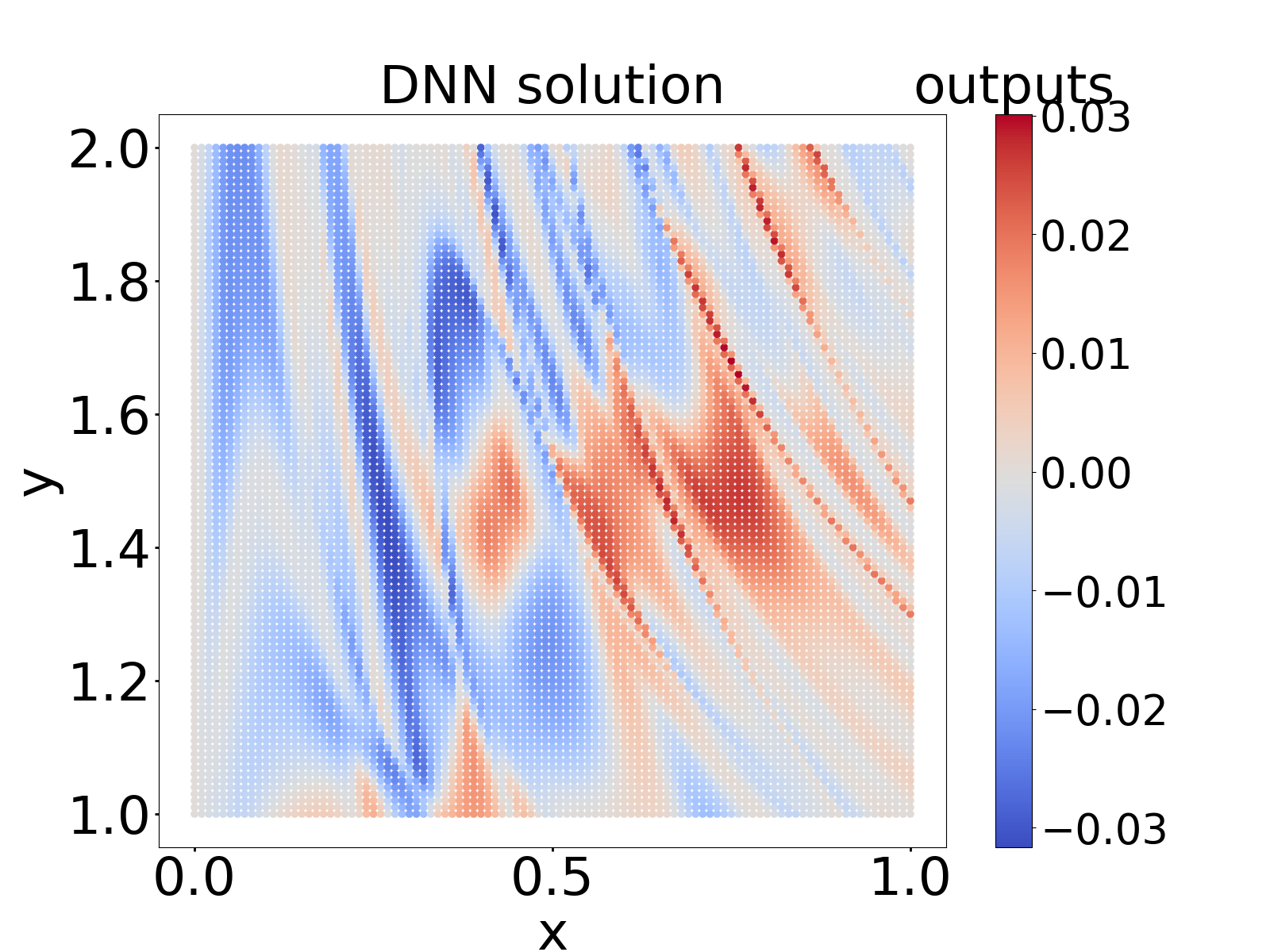}
    }
    \subfigure[Difference between two solutions]{
        \includegraphics[scale=0.11]{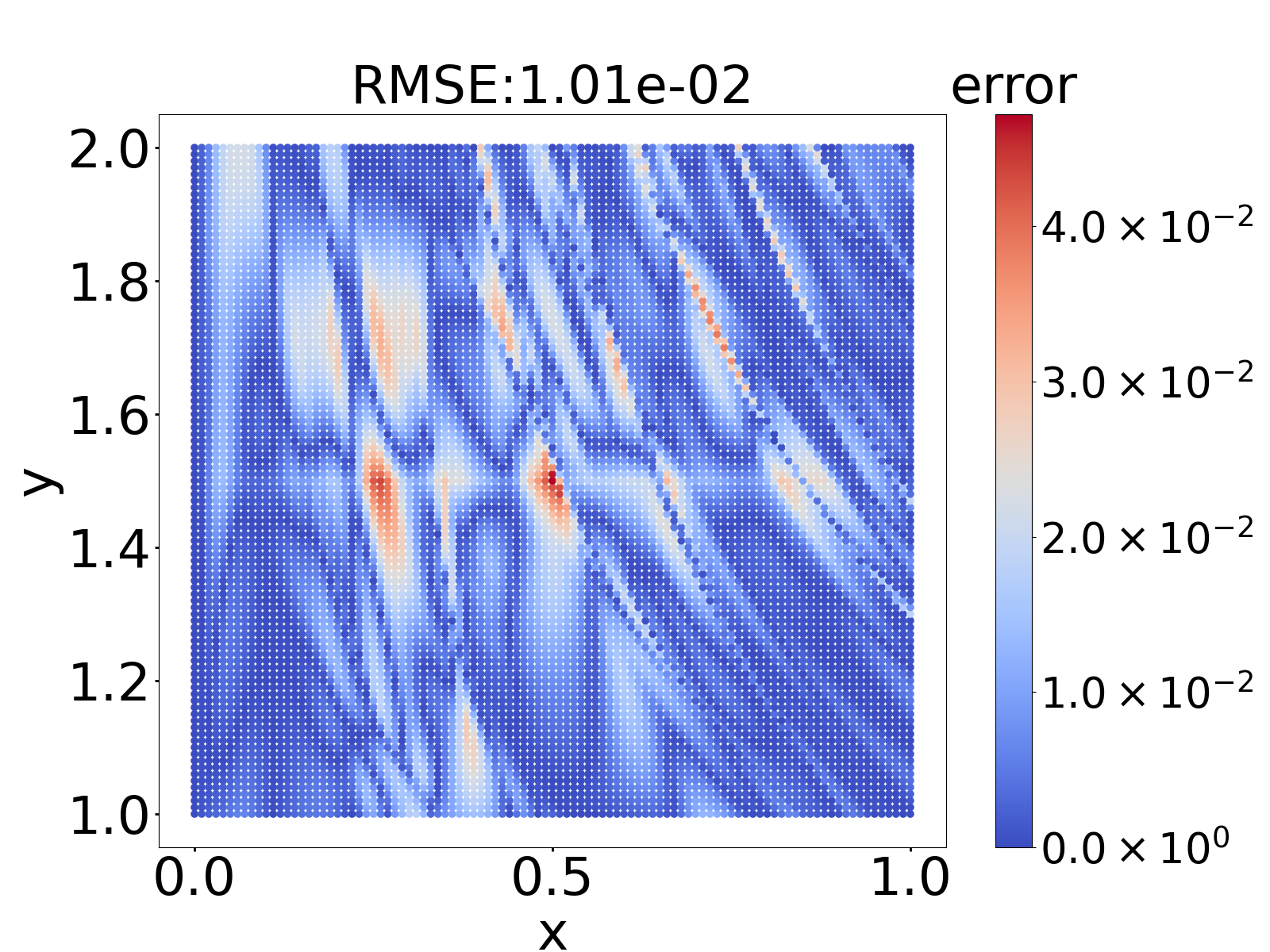}
    }
    \caption{Example 3. Trained by only the $10 \times 10$ labeled data. (a) Training loss. The  blue curve is the supervised risk Eq. (\ref{constructed_leastsquareloss}).
        (b,c)
        Comparison between analytic solution and MOD-Net solution on $x=0.5$ and $y=1.55$.
        (d,e,f)
        Comparison between analytic solution and  MOD-Net solution on $101\times101$ grid points.
    }
    \label{fig:example5_constructed_onlydata}
\end{figure}

Then, we set the $\lambda_1$, $\lambda_2$ be zero  and train the MOD-Net only by the $10 \times 10$ equidistributed isometric labeled data, i.e., $10$ equidistributed isometric points in $x$ direction and $10$ equidistributed isometric points in $y$ direction.  The empirical risk, i.e., training loss is shown in Fig. \ref{fig:example5_constructed_onlydata} (a). Similarly, to visualize  the performance of the obtained solution $u_{\vtheta}$, we plot the MOD-Net solution and the corresponding analytic solution on $x=0.5$ and $y=1.55$. As shown in Fig. \ref{fig:example5_constructed_onlydata} (b,c), the MOD-Net solution deviates from the analytic solution.
We also show the analytic solution and  MOD-Net solution on $101\times101$ equidistributed isometric grid points by color in Fig. \ref{fig:example5_constructed_onlydata} (d,e,f). The root mean square error (RMSE) is $\sim 1.01\times 10^{-2}$ and relative error is $\sim 84.24\%$.

\begin{figure}
    \centering

    \subfigure[Training loss]{
        \includegraphics[scale=0.11]{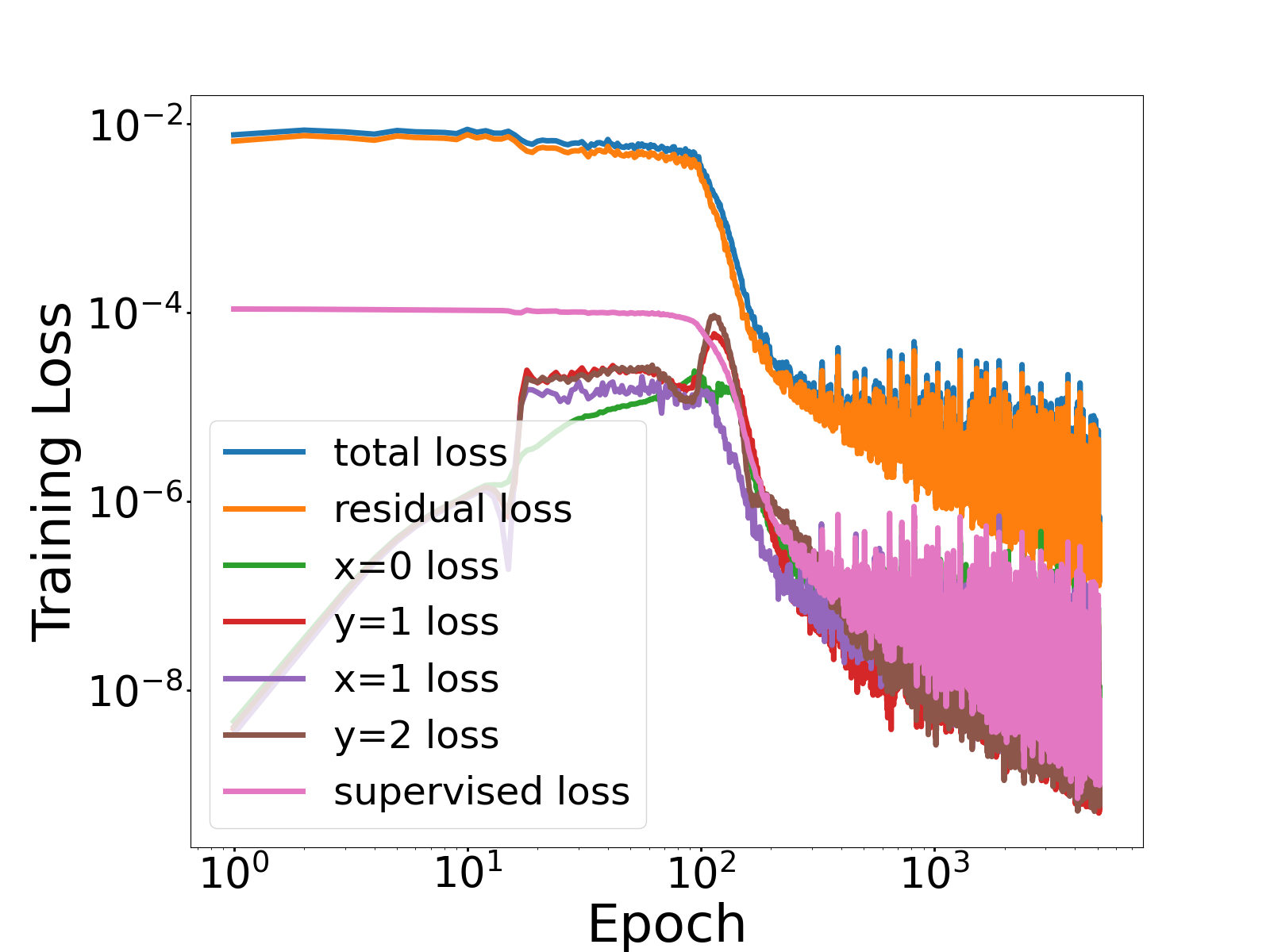}
    }
    \subfigure[$x=0.5$]{
        \includegraphics[scale=0.11]{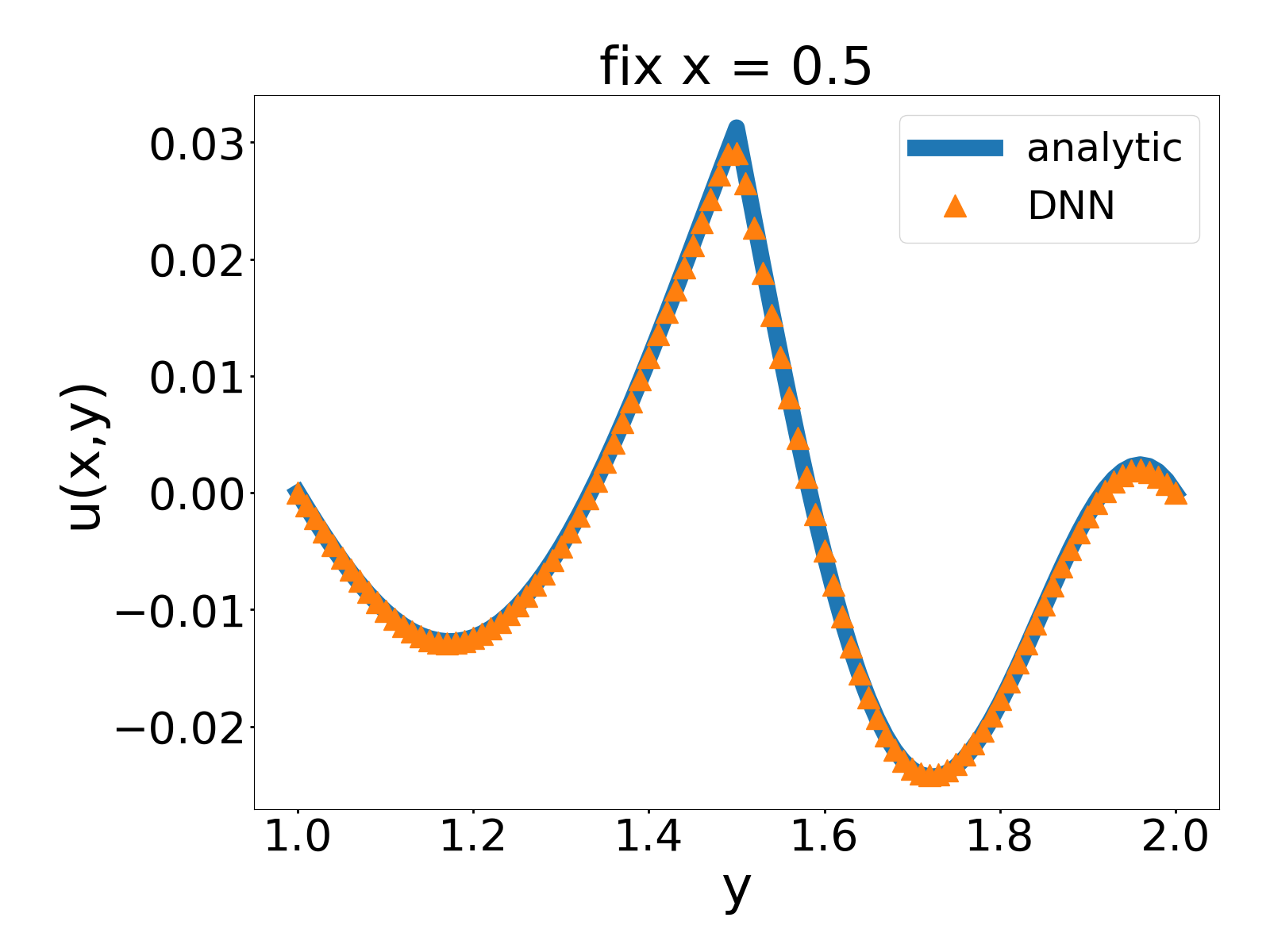}
    }
    \subfigure[$y=1.55$]{
        \includegraphics[scale=0.11]{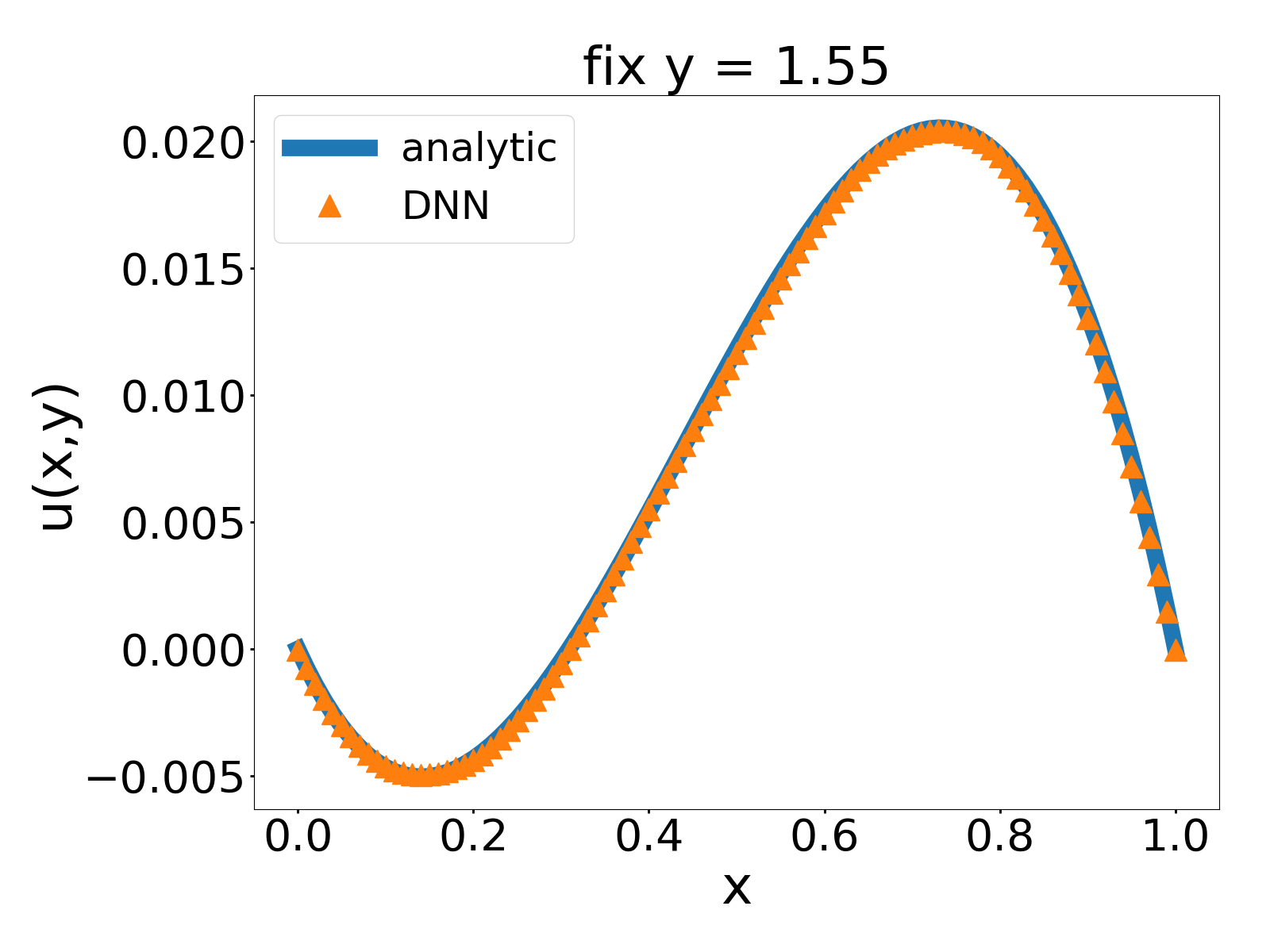}
    }
    \subfigure[Analytic solution]{
        \includegraphics[scale=0.11]{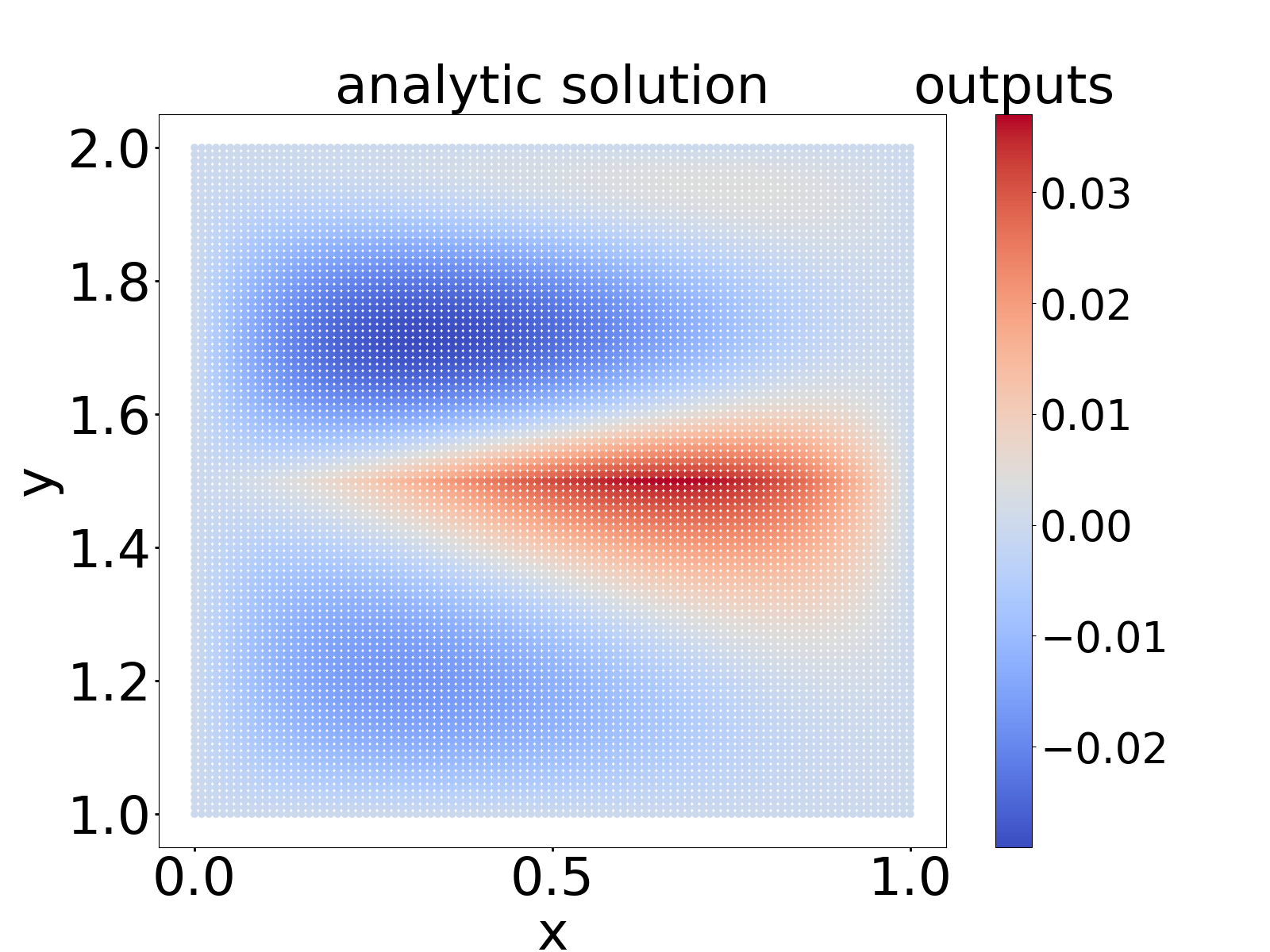}
    }
    \subfigure[MOD-Net solution]{
        \includegraphics[scale=0.11]{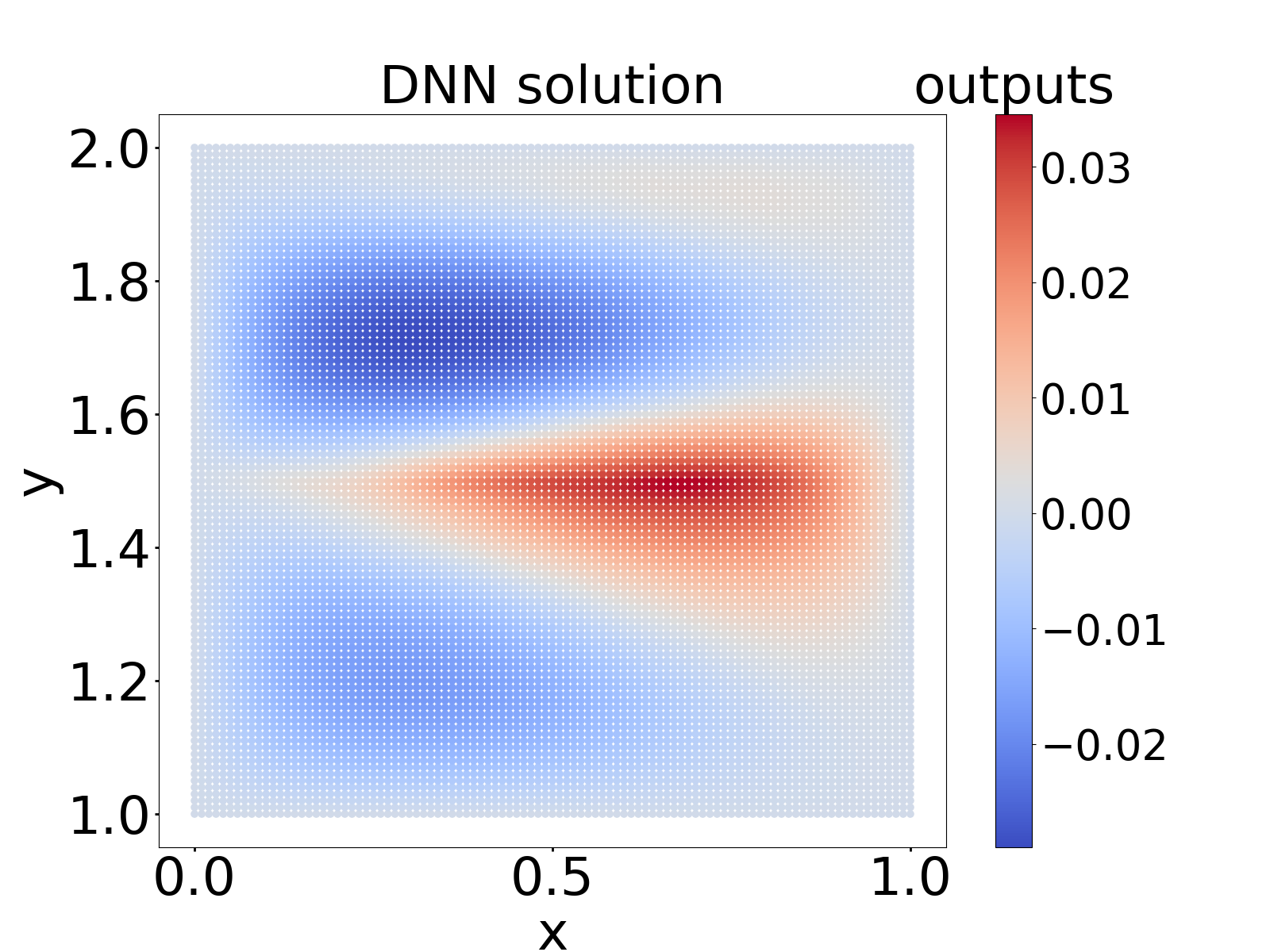}
    }
    \subfigure[Difference between two solutions]{
        \includegraphics[scale=0.11]{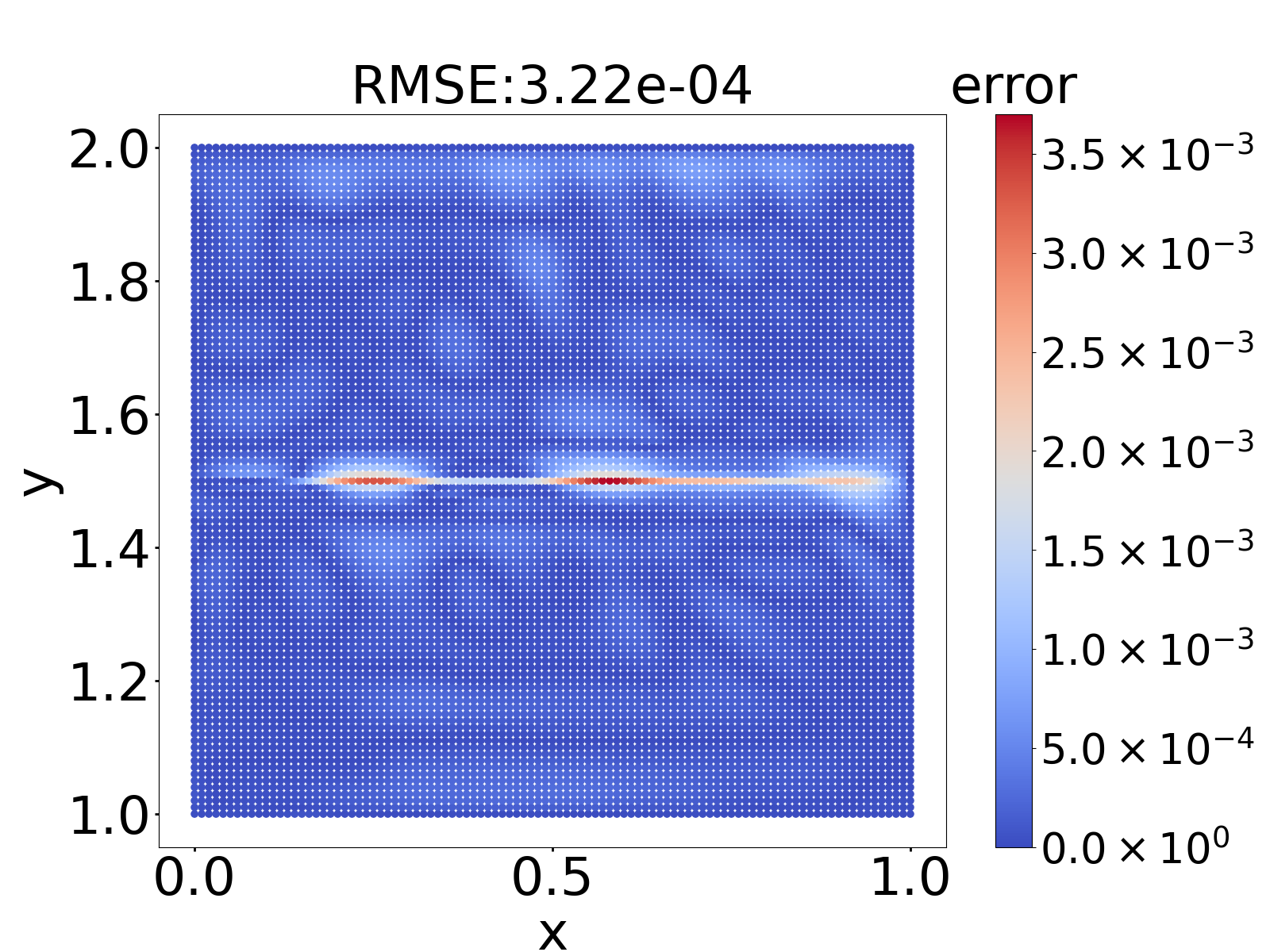}
    }
    \caption{Example 3. Trained by the governing PDE, boundary condition,and coarse grid data together. (a) Training loss. The  blue curve is the total risk Eq. (\ref{constructed_leastsquareloss}). The other six curves represent six terms in the total risk Eq.  (\ref{constructed_leastsquareloss}), respectively.
        (b,c)
        Comparison between analytic solution and MOD-Net solution on $x=0.5$ and $y=1.55$.
        (d,e,f)
        Comparison between analytic solution and  MOD-Net solution on $101\times101$ grid points.
    }
    \label{fig:example5_constructed_PDEanddata}
\end{figure}

To sum up, with the information of only PDE or only a few labeled data, the MOD-Net cannot be trained well, however, combining these two information, we can train the MOD-Net very well. See results in Fig. \ref{fig:example5_constructed_PDEanddata}. The root mean square error (RMSE) is $\sim 3.22\times 10^{-4}$ and relative error is $\sim 2.68\%$.

\section{Numerical experiments: One-dimensional radiative transfer equation} \label{sec:boltz}
In this section, we would apply MOD-Net to solve  one-dimensional steady radiative transfer equation (RTE).
Consider the density of particles in a bounded domain that interact with a background through absorption and scattering processes. The density function $u(\vx, \vvv)$ follows the RTE
\begin{equation}
    \begin{aligned}
         & \vvv \cdot \nabla u(\vx, \vvv)+\frac{\sigma_{T}(\vx)}{\varepsilon(\vx)} u(\vx, \vvv)=\frac{1}{|S|} \left(\frac{\sigma_{T}(\vx)}{\varepsilon(\vx)}-\varepsilon(\vx) \sigma_{a}(\vx)\right) \int_{S} u(\vx, \vxi) \diff{\vxi},
        \\ & u(\vx, \vvv)=\phi(\vx, \vvv), \quad \vx \in \Gamma=\partial \Omega, \quad \vvv \in S, \quad \vvv \cdot \vn_{\vx}<0,
        \label{RTE}
    \end{aligned}
\end{equation}
where $\vx \in \Omega\subset\sR^{d}$ is the $d$-dimensional space variable,
$\vvv$ is the angular variable on unit ball $S^{d-1}\subset\sR^{d}$, $\vn_{\vx}$ is the outward normal vector at $\vx$ on the boundary,
$\sigma_{a}$ is the absorption coefficient, $\sigma_{T}$ is the total scattering coefficient and $\varepsilon$ is Knudsen number. For simplifying the PDE and well describing the real situation,  we use the isotropic hypothesis, then $\sigma_{a}$ and $\sigma_{T}$ are only related to the space variable $\vx$.

High-accuracy numerical methods are developed to solve RTE, e.g., Tailored Finite Point Method (TFPM). However, these numerical methods are usually based on discrete grids. To improve the accuracy, the mesh grid has to be finer and the solving time exponentially increases. Another disadvantage is that these numerical methods usually only solve one particular PDE, i.e., if we change
anyone of $\sigma_{T}$, $\sigma_{a}$, $\varepsilon$ or
$\phi$, we need to reuse this algorithm to solve it. To  alleviate these problems,  MOD-Net approach learns the operator $\fG: (\phi,\sigma_{T},\sigma_{a},\varepsilon) \mapsto u$.


We consider one-dimensional RTEs. Set $\Omega=[x_L,x_R]$, $S=[-1,1]$, normal vector ${n_{x_L}}=-1$, ${n_{x_R}}=1$, we obtain
\begin{equation}
    \begin{aligned}
         & v \partial_{x} u(x,v)+\frac{\sigma_{T}(x)}{\varepsilon(x)} u(x,v)=\left(\frac{\sigma_{T}(x)}{\varepsilon(x)}-\varepsilon(x) \sigma_{a}(x)\right) \frac{1}{2} \int_{-1}^{1} u\left(x, \xi \right) \, \diff{\xi}, \\
         & u\left(x_{L}, v\right)=\phi_{L}(v), \quad v>0,                                                                                                                                                                  \\
         & u\left(x_{R}, v\right)=\phi_{R}(v), \quad v<0.
        \label{onedimensionRTE}
    \end{aligned}
\end{equation}
For convenience, in this paper, we fixed the $\sigma=(\sigma_T,\sigma_a,\varepsilon)$, and our goal is only to learn a operator mapping from $\phi_{L},\phi_{R}$ to the solution $u$ of the RTE.

\subsection{Use DNN to fit Green's function}

Denote $$\fL [u]=\vvv \cdot \nabla u(\vx, \vvv)+\frac{\sigma_{T}(\vx)}{\varepsilon(\vx)} u(\vx, \vvv)-\frac{1}{|S|} \left(\frac{\sigma_{T}(\vx)}{\varepsilon(\vx)}-\varepsilon(\vx) \sigma_{a}(\vx)\right) \int_{S} u(\vx, \vxi) \diff{\vxi},$$
the RTE (\ref{RTE}) can be rewritten as the following linear PDE,
\begin{equation}
    \left\{
    \begin{aligned}
         & \fL [u](\vx, \vvv) = 0,\quad \vx \in \Omega,\quad \vvv \in S,                               \\
         & u(\vx, \vvv) = \phi(\vx, \vvv), \quad \vx \in \partial \Omega,\quad \vvv \cdot \vn_{\vx}<0,
    \end{aligned}
    \right.
    \label{LinearhomogeneousPDE}
\end{equation}
which is a special case of linear PDE (\ref{linearPDE}).

With Green's function method, the solution of (\ref{LinearhomogeneousPDE}),
can be represented by the following formula,
\begin{equation}
    \begin{aligned}
        u(\vx,\vvv;\phi,\sigma) =
        \int_{\partial \Omega } \int_{S \cap \{\vvv'|\vvv' \cdot \vn_{\vx'}<0 \}}
        G(\vx,\vx',\vvv,\vvv') \phi(\vx, \vvv)\, \diff{\vvv'}\diff{\vx'},
    \end{aligned}
    \label{GeneralRTErepr}
\end{equation}
where $G(\vx,\vx',\vvv,\vvv')$ is the solution of the following PDE,
$$
    \left\{
    \begin{aligned}
         & \fL [G](\vx, \vvv) = 0,\quad \vx \in \Omega ,\quad \vvv \in S,                                                                                                \\
         & G(\vx,\vx',\vvv,\vvv') = \delta(\vvv-\vvv')\delta(\vx-\vx'), \quad \vx,\vx' \in \partial \Omega ,\quad \vvv \cdot \vn_{\vx}<0,\quad \vvv' \cdot \vn_{\vx'}<0.
    \end{aligned}
    \right.
$$

Since in the one-dimensional example (\ref{onedimensionRTE}), we consider $\Omega=[x_L,x_R]$, then $x'$ in $G(x,x',$ $v,v')$ only have two values. We can rewrite the integral
in (\ref{GeneralRTErepr}) by the following formula,
\begin{equation}
    \begin{aligned}
        u(x,v;\phi_L,\phi_{R}) = \int^1_0 G_L(x,v,v')\phi_L(v')\diff{v'} + \int^0_{-1} G_R(x,v,v')\phi_{R}(v')\diff{v'},
    \end{aligned}
    \label{1D_repre}
\end{equation}
where $G_L$, $G_R$ are the solution of the following two PDEs respectively,
$$
    \left\{
    \begin{aligned}
        v \partial_{x} G_L+\frac{\sigma_{T}}{\varepsilon} G_L & =\left(\frac{\sigma_{T}}{\varepsilon}-\varepsilon \sigma_{a}\right) \frac{1}{2} \int_{-1}^{1} G_L\left(x, \xi,v'\right)\diff{\xi}, \\
        G_L\left(x_{L}, v, v'\right)                          & =\delta(v-v'), \quad v>0,                                                                                                          \\ \quad G_L\left(x_{R}, v, v'\right)&=0, \quad v<0,
    \end{aligned}
    \right.
$$
$$
    \left\{
    \begin{aligned}
        v \partial_{x} G_R+\frac{\sigma_{T}}{\varepsilon} G_R & =\left(\frac{\sigma_{T}}{\varepsilon}-\varepsilon \sigma_{a}\right) \frac{1}{2} \int_{-1}^{1} G_R\left(x, \xi,v';\sigma\right) \diff{\xi}, \\
        G_R\left(x_{L}, v, v';\sigma\right)                   & =0, \quad v>0,                                                                                                                             \\ \quad G_R\left(x_{R}, v, v';\sigma\right)&=\delta(v-v'), \quad v<0.
    \end{aligned}
    \right.
$$

With the help of Green's function method, to achieve our goal of  fitting the operator $\fG: (\phi_{L},\phi_{R}) \mapsto u$, where $u$ is the solution of (\ref{onedimensionRTE}) with given $\phi_{L},\phi_{R}$,
We  fit $\fG_L$,$\fG_R$ with DNNs of hidden layer size $128$-$256$-$256$-$128$ equipped with activation function tanh, i.e., a DNN  $G_{\vtheta_L}(x,v,v')$ is trained to represent
$G_L(x,v,v')$, similarly, another DNN  $G_{\vtheta_R}(x,v,v')$ is for $G_R(x,v,v')$.
Since the solution reflects the distribution function, we use the exponential function to make sure that $G_L$ and $G_R$ are positive.


In Eq. (\ref{1D_repre}),
when we calculate the $1$D integral,
similarly
we use the Gauss-Legendre quadrature. Then we obtain the following representation of MOD-Net solution,
\begin{equation}
    \begin{aligned}
        u_{\vtheta_1,\vtheta_2}(x,v;\phi_L,\phi_R)
         & = \sum_{v' \in S_{G_v}^+} \omega_{v_+'}G_{\vtheta_L}(x,v,v_+')\phi_L({v_+'})
        \\
         & + \sum_{v' \in S_{G_v}^-} \omega_{v_-'}G_{\vtheta_R}(x,v,v_-')\phi_R({v_-'}),
    \end{aligned}
    \label{RTE_Gausianint}
\end{equation}
where $S_{G_v}^+ \subset [0,1]$, $S_{G_v}^- \subset  [-1,0]$ consist of fixed points determined by Gauss-Legendre quadrature and $\omega_{v_+'}$, $\omega_{v_-'}$ are corresponding coefficients.
\subsection{Empirical risk function}
For one-dimensional case, to train the neural networks, we would utilize the information of PDE, i.e., governing equation and boundary condition, and a few data $S^{u,k}=\{x_i,v_i,$ $u^{k}(x_i,v_i)\}_{i\in [n_k]}$ for each $\{\phi_{L}^k,\phi_{R}^k\},k=1,2,\cdots,K$, where $u^{k}(\cdot)=u(\cdot;\phi_{L}^k,\phi_{R}^k)$. Note that  $S^{u,k}$ can be numerically solved by TFPM on coarse grid points, which is not computationally expensive. For each $k$, we uniformly sample a set of $(x,v)$ from $\Omega \times S=[x_L,x_R] \times [-1,1]$, i.e., $S^{\Omega,S,k}$ and uniformly sample a set of data from boundaries $\partial \Omega_L=\{(x_L,v)\}_{v\in[0,1]}$, $\partial \Omega_R=\{(x_R,v)\}_{v\in[-1,0]}$, respectively, i.e., $S^{\partial\Omega_L,k}$,$S^{\partial\Omega_R,k}$.

We use the general definition of empirical risk (\ref{lossfunctiondefinition}) for one-dimensional RTE.  The empirical risk of solving RTE is as follows,
\begin{align*} 
    \RS & = \frac{1}{K}\sum_{k\in[K]}  \Big(
    \lambda_1 \frac{1}{|S^{\Omega,S,k}|} \sum_{(x,v)\in S^{\Omega,S,k}}
    \| v \partial_{x} u_{\vtheta_1,\vtheta_2}(x,v;\phi_{L}^k,\phi_{R}^k)+\frac{\sigma_{T}(x)}{\varepsilon(x)} u_{\vtheta_1,\vtheta_2}(x,v;\phi_{L}^k,\phi_{R}^k)
    \\
        & \quad -\left(\frac{\sigma_{T}(x)}{\varepsilon}-\varepsilon \sigma_{a}(x)\right) \frac{1}{2} \int_{-1}^{1} u_{\vtheta_1,\vtheta_2}(x,\xi;\phi_{L}^k,\phi_{R}^k) \diff{\xi}
    \|_2^2
    \\
        & \quad +\lambda_{21} \frac{1}{|S^{\partial\Omega_L,k}|} \sum_{(x,v)\in S^{\partial\Omega_L,k}}  \| u_{\vtheta_1,\vtheta_2}(x, v;\phi_{L}^k,\phi_{R}^k)  -\phi_{L}^k(v) \|_2^2
    \\
        & \quad +\lambda_{22} \frac{1}{|S^{\partial\Omega,k}|} \sum_{(x,v)\in S^{\partial\Omega_R,k}}  \| u_{\vtheta_1,\vtheta_2}(x, v;\phi_{L}^k,\phi_{R}^k)  -\phi_{R}^k(v) \|_2^2
    \\
        & \quad + \lambda_3 \frac{1}{n_k} \sum_{i\in [n_k]} \|u_{\vtheta_1,\vtheta_2}(x_i,v_i;\phi_{L}^k,\phi_{R}^k)-u^{k}(x_i,v_i)\|_2^2
    \Big).
\end{align*}

For integral term in the first risk term related to the governing equation of PDE, we use the Gauss-Legendre numerical integral method. We obtain
\begin{align}
    \RS & = \frac{1}{K}\sum_{k\in[K]}  \Big(
    \lambda_1 \frac{1}{|S^{\Omega,S,k}|} \sum_{(x,v)\in S^{\Omega,S,k}}
    \| v \partial_{x} u_{\vtheta_1,\vtheta_2}(x,v;\phi_{L}^k,\phi_{R}^k)+\frac{\sigma_{T}(x)}{\varepsilon(x)} u_{\vtheta_1,\vtheta_2}(x,v;\phi_{L}^k,\phi_{R}^k) \nonumber
    \\
        & \quad -\left(\frac{\sigma_{T}(x)}{\varepsilon}-\varepsilon \sigma_{a}(x)\right) \frac{1}{2} \sum_{\xi\in S_v} \omega_{\xi} u_{\vtheta_1,\vtheta_2}(x, \xi;\phi_{L}^k,\phi_{R}^k)
    \|_2^2  \nonumber
    \\
        & \quad +\lambda_{21} \frac{1}{|S^{\partial\Omega_L,k}|} \sum_{(x,v)\in S^{\partial\Omega_L,k}}  \| u_{\vtheta_1,\vtheta_2}(x, v;\phi_{L}^k,\phi_{R}^k)  -\phi_{L}^k(v) \|_2^2 \label{RTEloss}
    \\
        & \quad +\lambda_{22} \frac{1}{|S^{\partial\Omega,k}|} \sum_{(x,v)\in S^{\partial\Omega_R,k}}  \| u_{\vtheta_1,\vtheta_2}(x, v;\phi_{L}^k,\phi_{R}^k)  -\phi_{R}^k(v) \|_2^2 \nonumber
    \\
        & \quad + \lambda_3 \frac{1}{n_k} \sum_{i\in [n_k]} \|u_{\vtheta_1,\vtheta_2}(x_i,v_i;\phi_{L}^k,\phi_{R}^k)-u^{k}(x_i,v_i)\|_2^2 \Big), \nonumber
\end{align}
where $S_v \subset [-1,1]$ consists of fixed integration points, determined by Gauss-Legendre quadrature and $\omega_{\xi}$'s are corresponding coefficients.

In practical applications,  it is usually not easy to measure $u(x,v)$, but density function $\rho(x)=\frac{1}{2} \int_{-1}^{1} u\left(x, \xi \right) \diff{\xi} $ can be measured. Therefore, we would utilize the information of PDE and a few data $S^{\rho,k}=\{x_i,$ $\rho^{k}(x_i)\}_{i\in [n_k]}$ for each $\{\phi_{L}^k,\phi_{R}^k\},k=1,2,\cdots,K$, where $\rho^{k}(\cdot)=\frac{1}{2} \int_{-1}^{1} u\left(\cdot, \xi;\phi_{L}^k,\phi_{R}^k \right) \diff{\xi} $.
And  we obtain another empirical risk,
\begin{align}
    \RS & = \frac{1}{K}\sum_{k\in[K]}  \Big(
    \lambda_1 \frac{1}{|S^{\Omega,S,k}|} \sum_{(x,v)\in S^{\Omega,S,k}}
    \| v \partial_{x} u_{\vtheta_1,\vtheta_2}(x,v;\phi_{L}^k,\phi_{R}^k)+\frac{\sigma_{T}(x)}{\varepsilon(x)} u_{\vtheta_1,\vtheta_2}(x,v;\phi_{L}^k,\phi_{R}^k) \nonumber
    \\
        & \quad -\left(\frac{\sigma_{T}(x)}{\varepsilon}-\varepsilon \sigma_{a}(x)\right) \frac{1}{2} \sum_{\xi\in S_v} \omega_{\xi} u_{\vtheta_1,\vtheta_2}(x, \xi;\phi_{L}^k,\phi_{R}^k)
    \|_2^2  \nonumber
    \\
        & \quad +\lambda_{21} \frac{1}{|S^{\partial\Omega_L,k}|} \sum_{(x,v)\in S^{\partial\Omega_L,k}}  \| u_{\vtheta_1,\vtheta_2}(x, v;\phi_{L}^k,\phi_{R}^k)  -\phi_{L}^k(v) \|_2^2 \label{RTElossrho}
    \\
        & \quad +\lambda_{22} \frac{1}{|S^{\partial\Omega,k}|} \sum_{(x,v)\in S^{\partial\Omega_R,k}}  \| u_{\vtheta_1,\vtheta_2}(x, v;\phi_{L}^k,\phi_{R}^k)  -\phi_{R}^k(v) \|_2^2 \nonumber
    \\
        & \quad + \lambda_3 \frac{1}{n_k} \sum_{i\in [n_k]} \|\frac{1}{2} \sum_{\xi'\in S_v'} \omega_{\xi'} u_{\vtheta_1,\vtheta_2}(x_i, \xi';\phi_{L}^k,\phi_{R}^k)-\rho^{k}(x_i)\|_2^2 \Big), \nonumber
\end{align}
where $S_v,S_v' \subset [-1,1]$ consists of fixed integration points, determined by Gauss-Legendre quadrature and $\omega_{\xi},\omega_{\xi'}$'s are corresponding coefficients.

\subsection{Learning process}
During training, for each epoch, we first randomly choose $\{(\phi_{L}^k,\phi_{R}^k)\}_{k\in[K]}$ and calculate the values of boundary condition $\phi_{L}^k$, $\phi_{R}^k$ on fixed integration points $v_+'\in S_{G_v}^+$, $v_-' \in S_{G_v}^-$ respectively.
Second, we randomly sample points  and obtain the data set  $S^{\Omega,S,k}$, $S^{\partial\Omega_L,k}$,$S^{\partial\Omega_R,k}$.
We obtain the data set $D_L= \{(x,v,v_+',\phi_{L}^k(v_+'))| (x,v)\in S^{\Omega,S,k} \cup S^{\partial\Omega_L,k} \cup S^{\partial\Omega_R,k},$ $v_+'\in S_{G_v}^+\}$ and $D_R= \{(x,v,v_-',\phi_{R}^k(v_-'))| (x,v)\in S^{\Omega,S,k} \cup S^{\partial\Omega_L,k} \cup S^{\partial\Omega_R,k},v_-'\in S_{G_v}^-\}$. In the following, for each $k$, we feed the data $D_L$, $D_R$ into the neural network $G_{\vtheta_L}(x,v,v')$, $G_{\vtheta_R}(x,v,v')$ respectively,
and calculate the empirical risk (\ref{RTEloss}) or (\ref{RTElossrho}) defined utilizing the governing equation of PDE, boundary condition and a few labeled data $S^{u,k}$.
We train neural network $G_{\vtheta_L}$ and $G_{\vtheta_R}$ with Adam to minimize the  empirical risk, and finally  obtain  well-trained Green's function DNNs $G_{\vtheta_L}$ and $G_{\vtheta_R}$, furthermore, according to (\ref{RTE_Gausianint}), we obtain a neural operator $u_{\vtheta_1,\vtheta_2}(x,v;\phi_{L},\phi_{R})$.

\subsection{Results}
To see the performance of obtained operator, we need the reference solution. Since it is difficult to obtain the analytical solution of RTE, we use the TFPM method to obtain the numerical solution as reference.

\paragraph{Example 4. MOD-Net with data-regularization for radiative transfer equation with  $\sigma$ varying with $x$}

\begin{figure}[!htb]
    \centering
    \subfigure[$\sigma(x)$]{
        \includegraphics[scale=0.11]{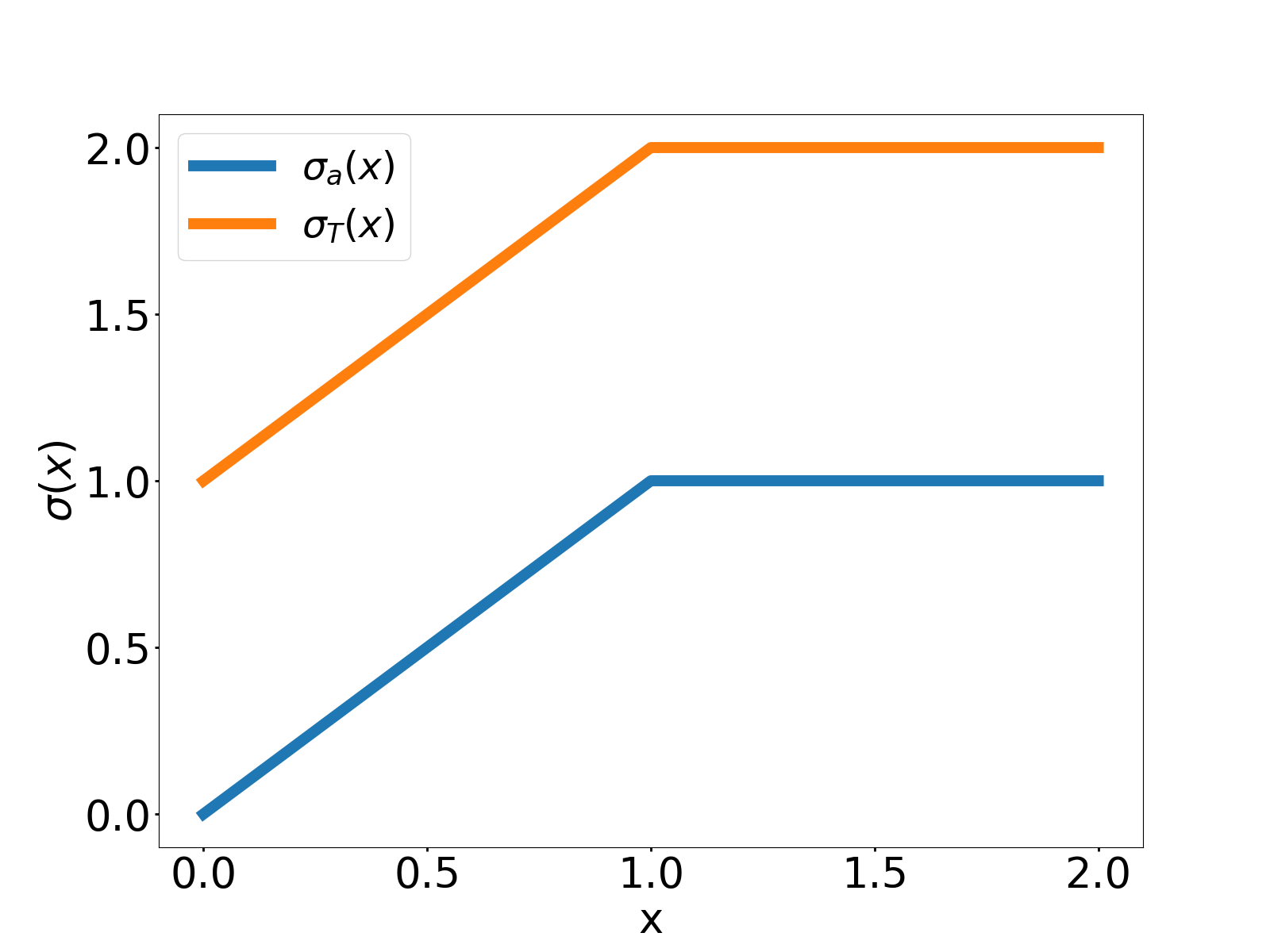}
    }
    \subfigure[Training loss]{
        \includegraphics[scale=0.11]{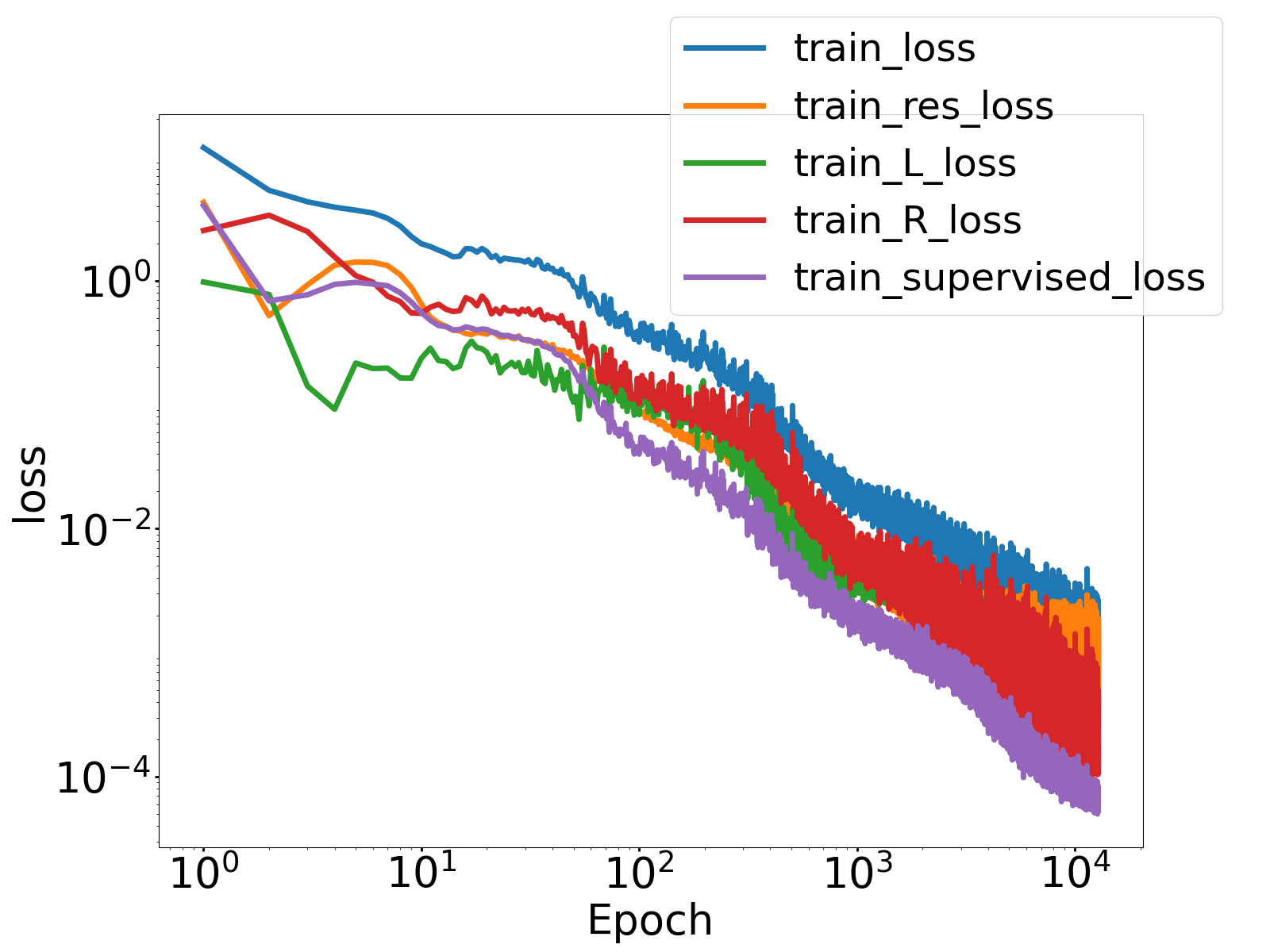}
    }
    \subfigure[$a_2=0.02$]{
        \includegraphics[scale=0.11]{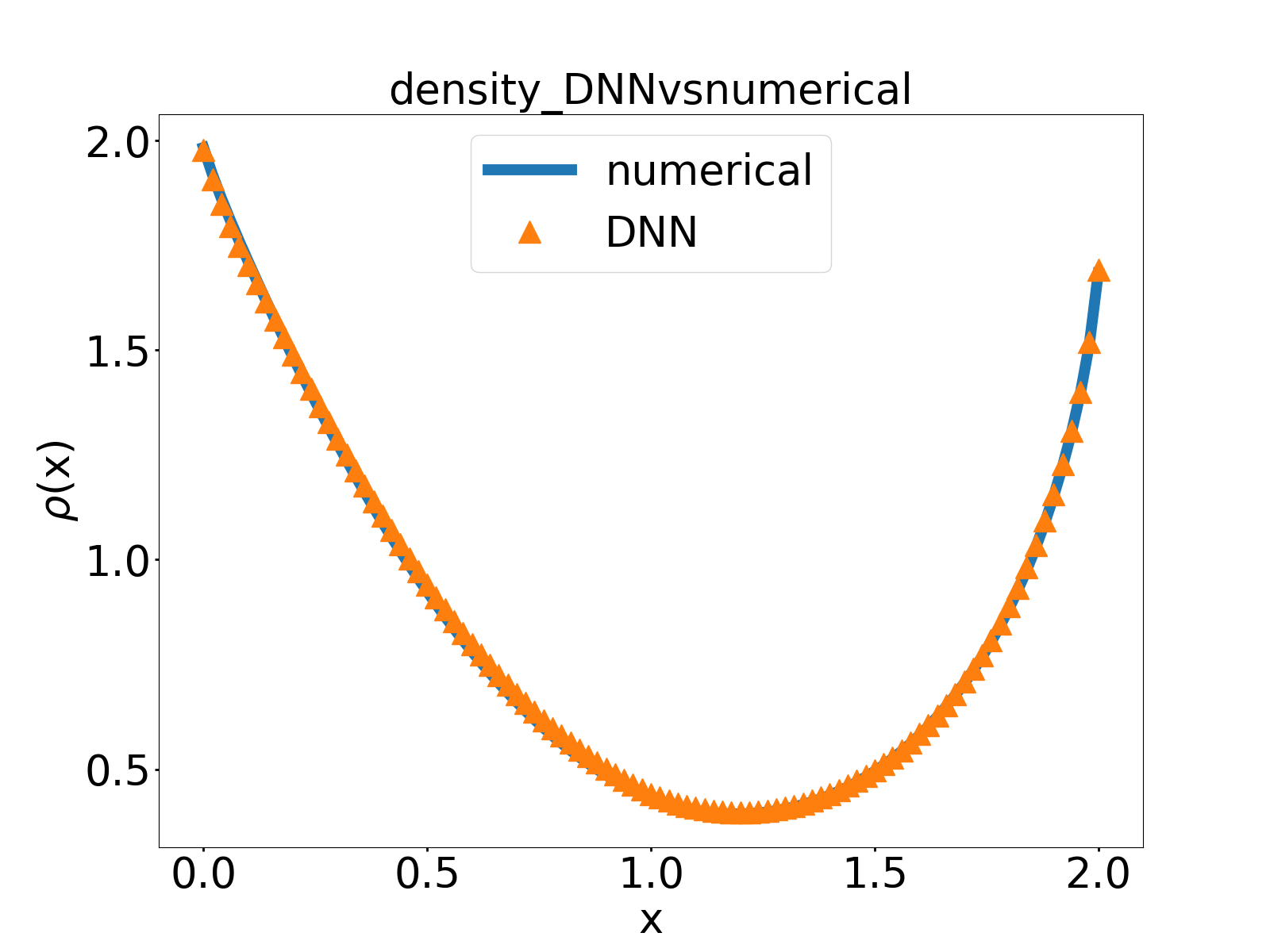}
    }
    \subfigure[$a_2=0.11$]{
        \includegraphics[scale=0.11]{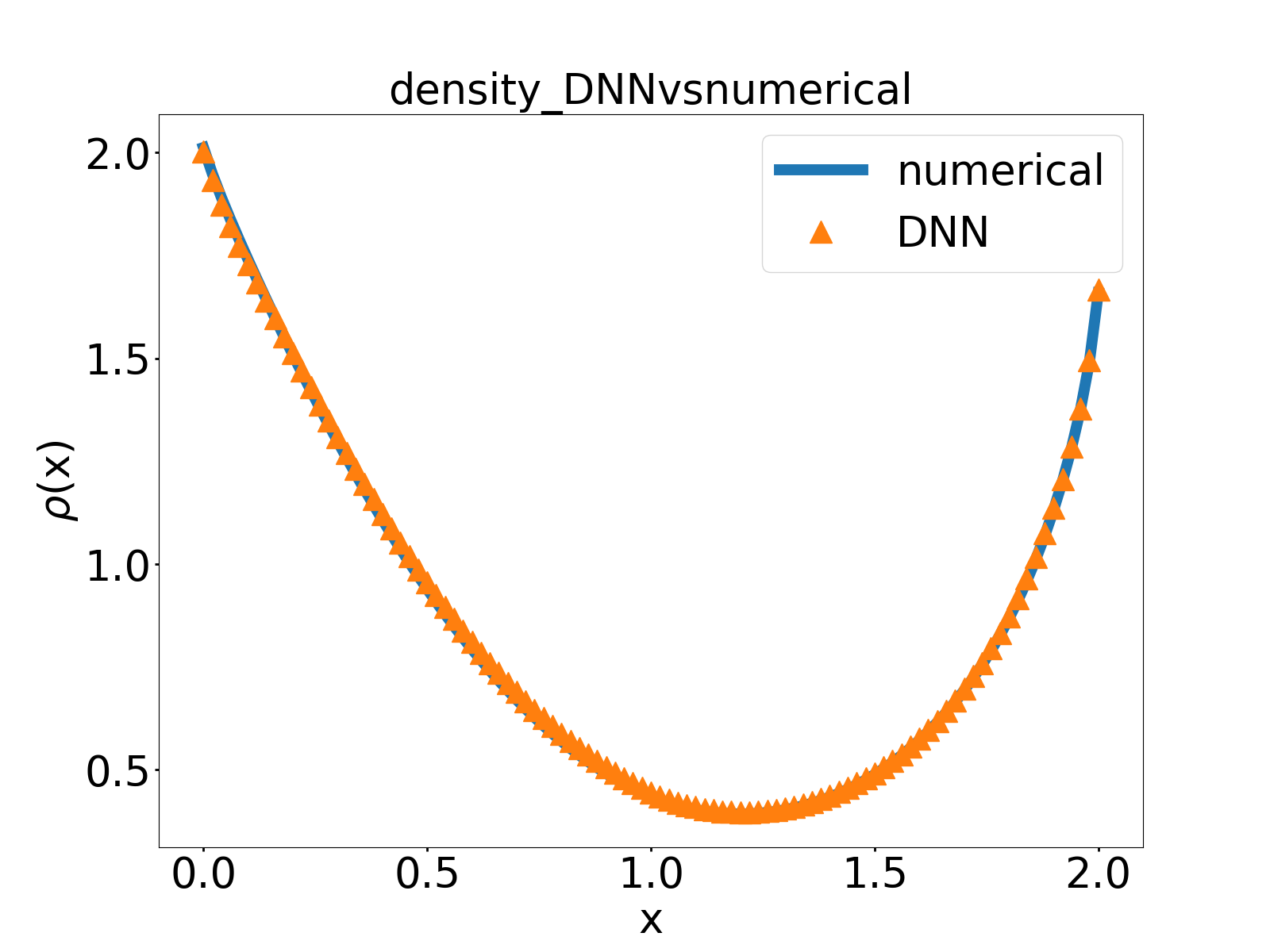}
    }
    \subfigure[$a_2=0.51$]{
        \includegraphics[scale=0.11]{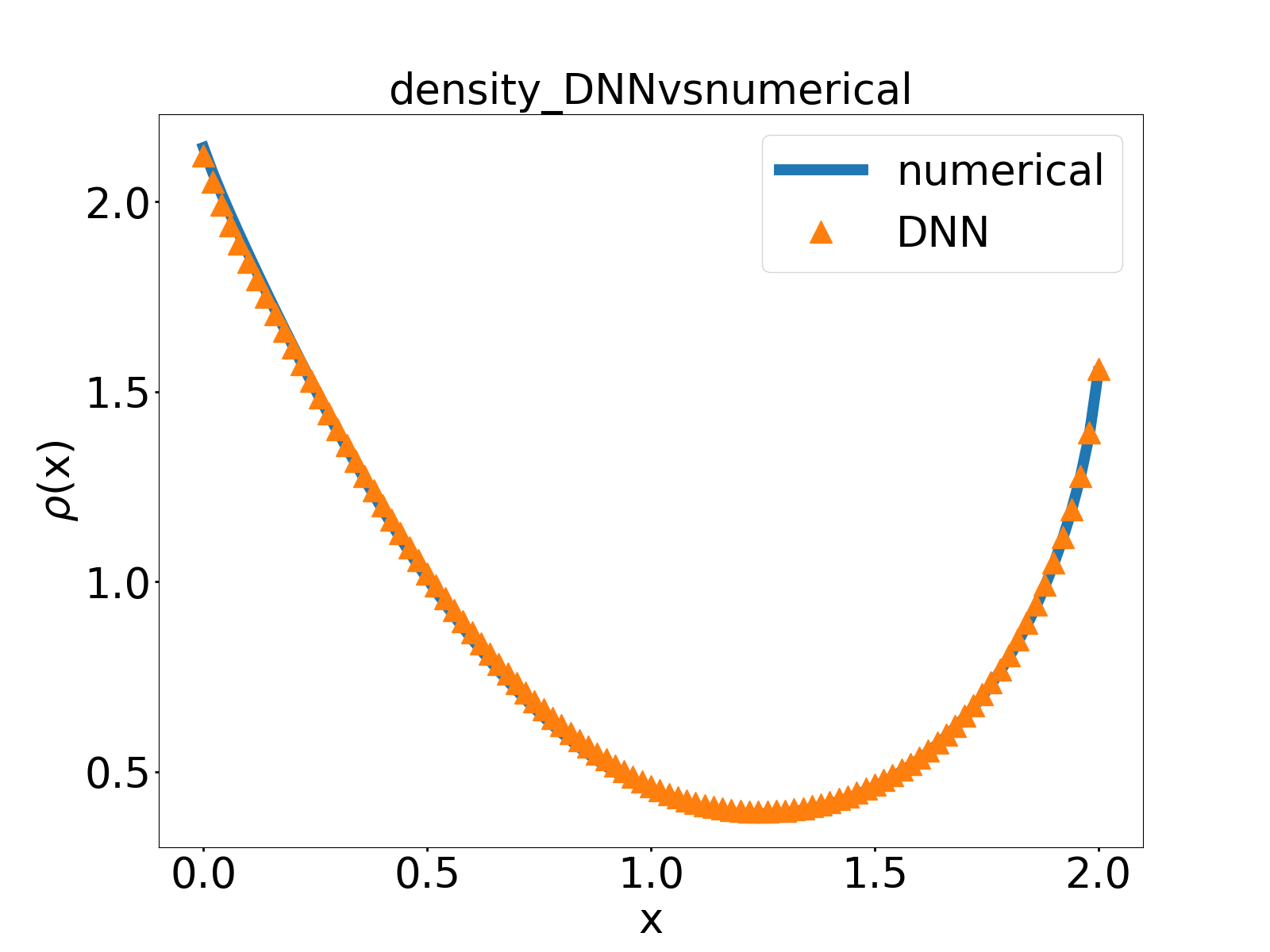}
    }
    \subfigure[$a_2=0.99$]{
        \includegraphics[scale=0.11]{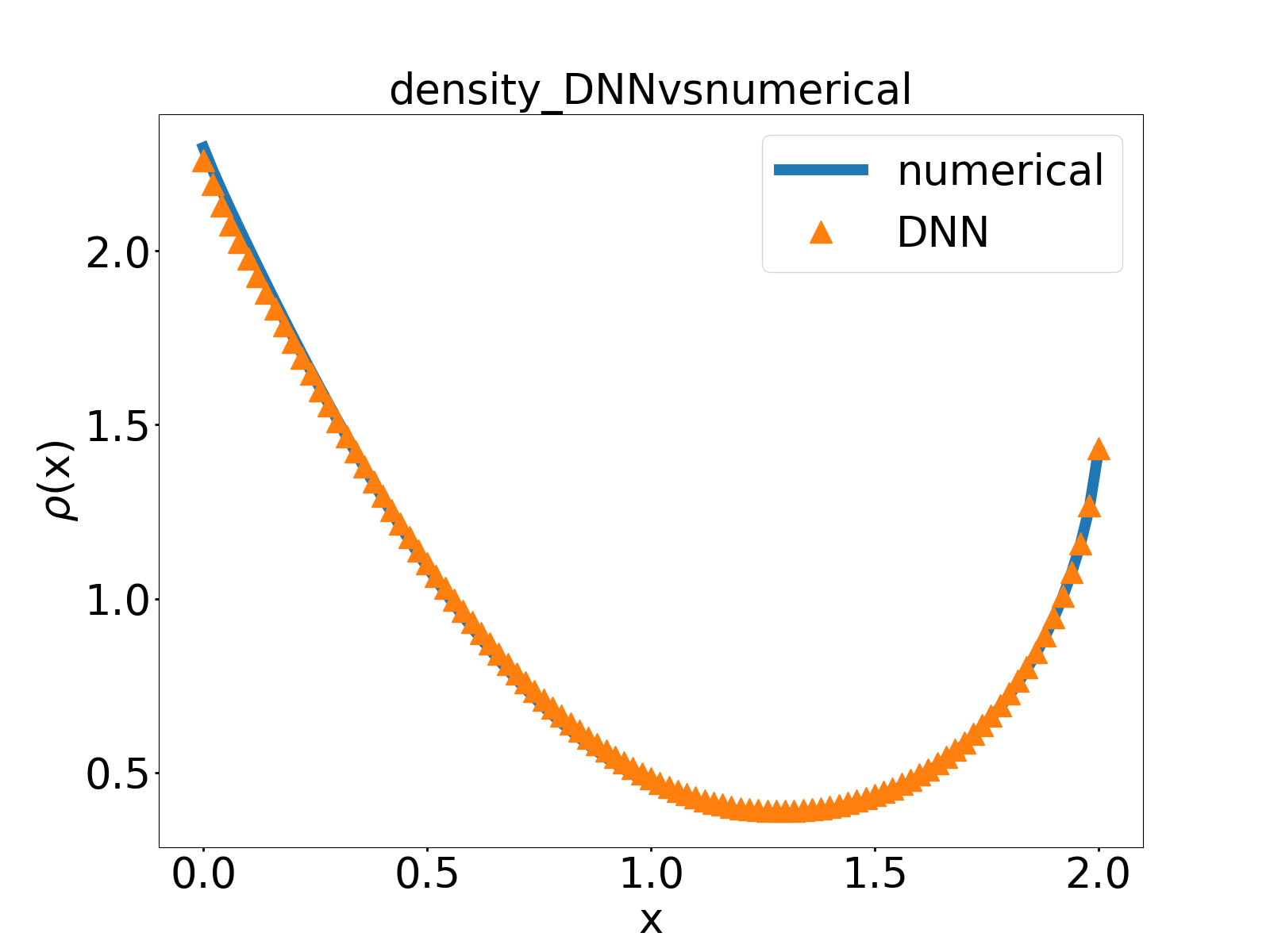}
    }
    \caption{Example 4. Training MOD-Net by the governing PDE, boundary condition, and  coarse grid data together. (a)  $\sigma_T(x)$ and $\sigma_a(x)$ vs. $x$. (b) Training loss vs. epoch for the learning of the operator of RTE with continuous $\sigma_T(x)$, $\sigma_a(x)$. The blue curve is the total risk Eq. (\ref{RTEloss}). The other four curves represent four terms in the total risk Eq. (\ref{RTEloss}), respectively. (c,d,e,f) Comparison between numerical solution and  MOD-Net solution on density function corresponding to each of four boundary conditions $\phi_{L}=\phi_{R}=\cos(v) + a_2\sin(v) + 2$ determined by $a_2=0.02,0.11,0.51,0.99$.}
    \label{fig:bolzmann_continuous_op}
\end{figure}

Only for illustration of data regularization in solving RTE, we  consider a simple case, $x_L=0, x_R=2$, and
\begin{equation}
    \sigma_T(x) =
    \left\{
    \begin{aligned}
        x+1, \quad 0\leq x<1, \\
        2, \quad 1\leq x\leq2,
    \end{aligned}\right.
\end{equation}
and
\begin{equation}
    \sigma_a(x) =
    \left\{
    \begin{aligned}
        x, \quad  0\leq x<1, \\
        1, \quad  1\leq x\leq2.
    \end{aligned}\right.
\end{equation}
the graphs of which are shown in Fig. \ref{fig:bolzmann_continuous_op}(a).
To avoid the multi-scale phenomenon,  we take $\varepsilon(x)=1$.
We set the boundary condition $\phi_{L}=\phi_{R}=a_1\cos(\omega v) + a_2\sin(\omega v) + 2$, which is determined by $a_1,a_2,\omega $. In fact, many boundary conditions can be represented with these basic functions.
For convenience, we set $\omega =1$,$a_1=1$ and only change $a_2$. Then the boundary condition $\phi_{L}=\phi_{R}=\cos(v) + a_2\sin(v) + 2$.

\begin{figure} [H]
    \centering
    \subfigure[Training loss]{
        \includegraphics[scale=0.11]{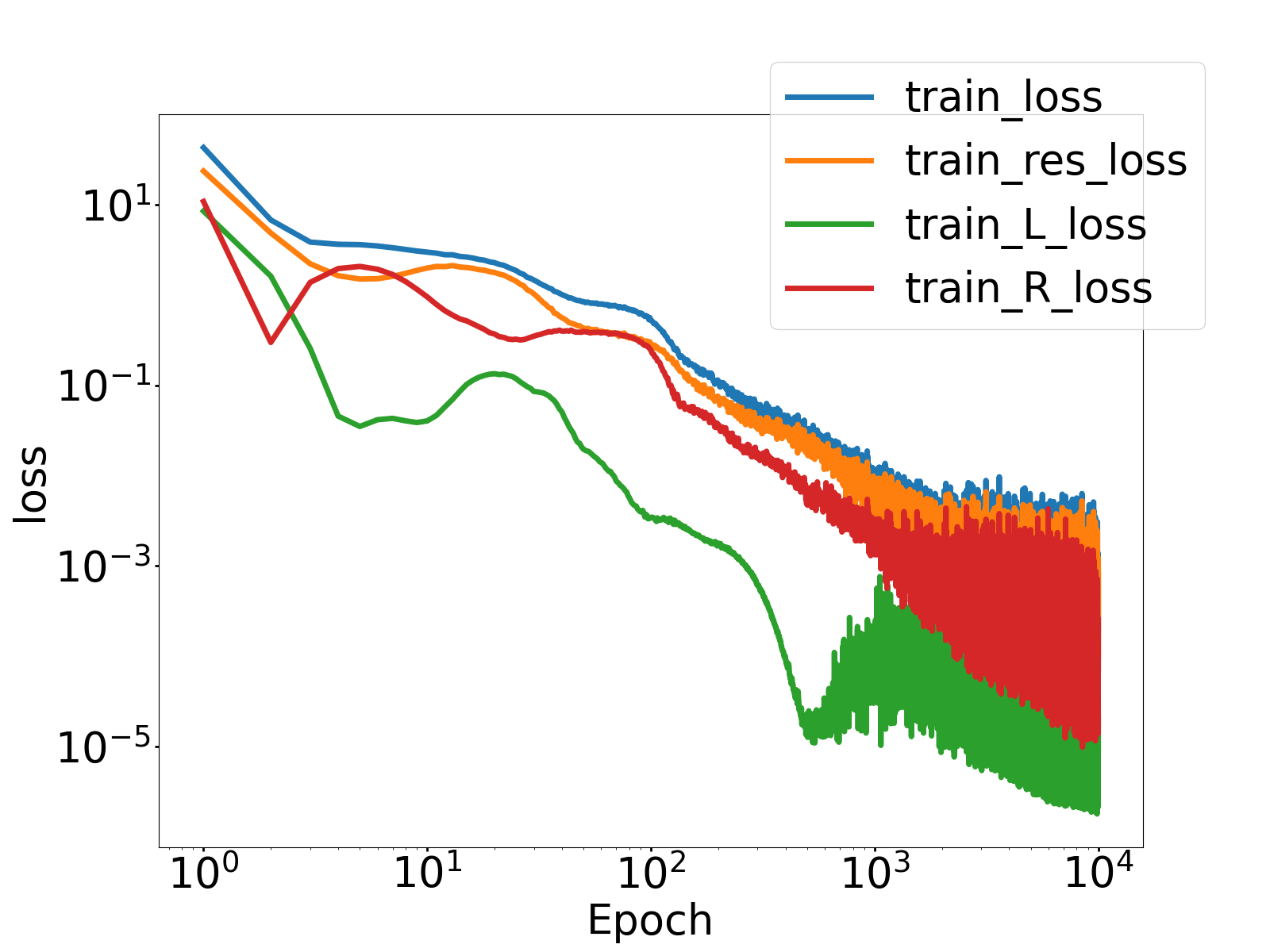}
    }
    \subfigure[Density function]{
        \includegraphics[scale=0.11]{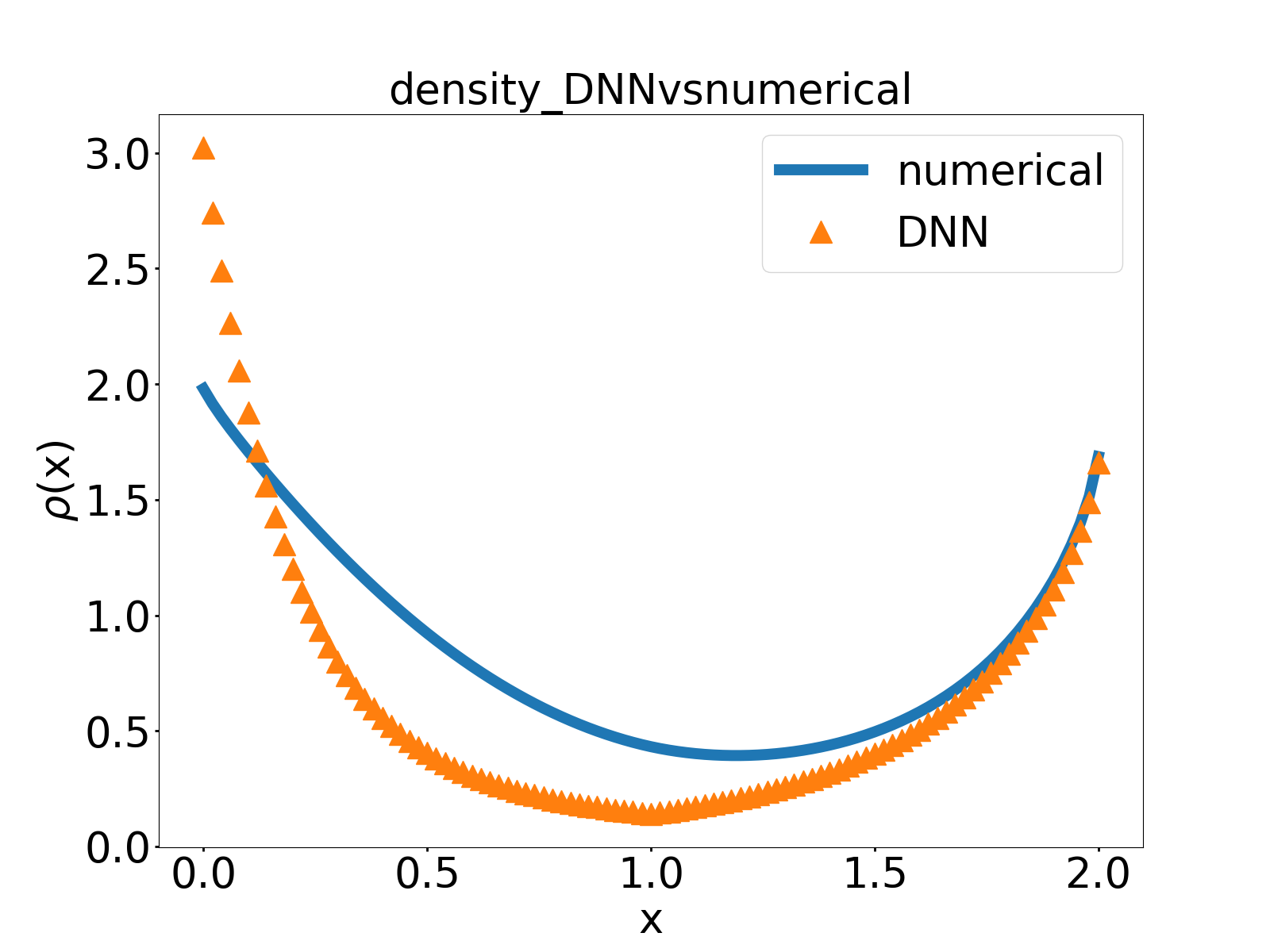}
    }
    \\
    \subfigure[solution on $x=0$]{
        \includegraphics[scale=0.11]{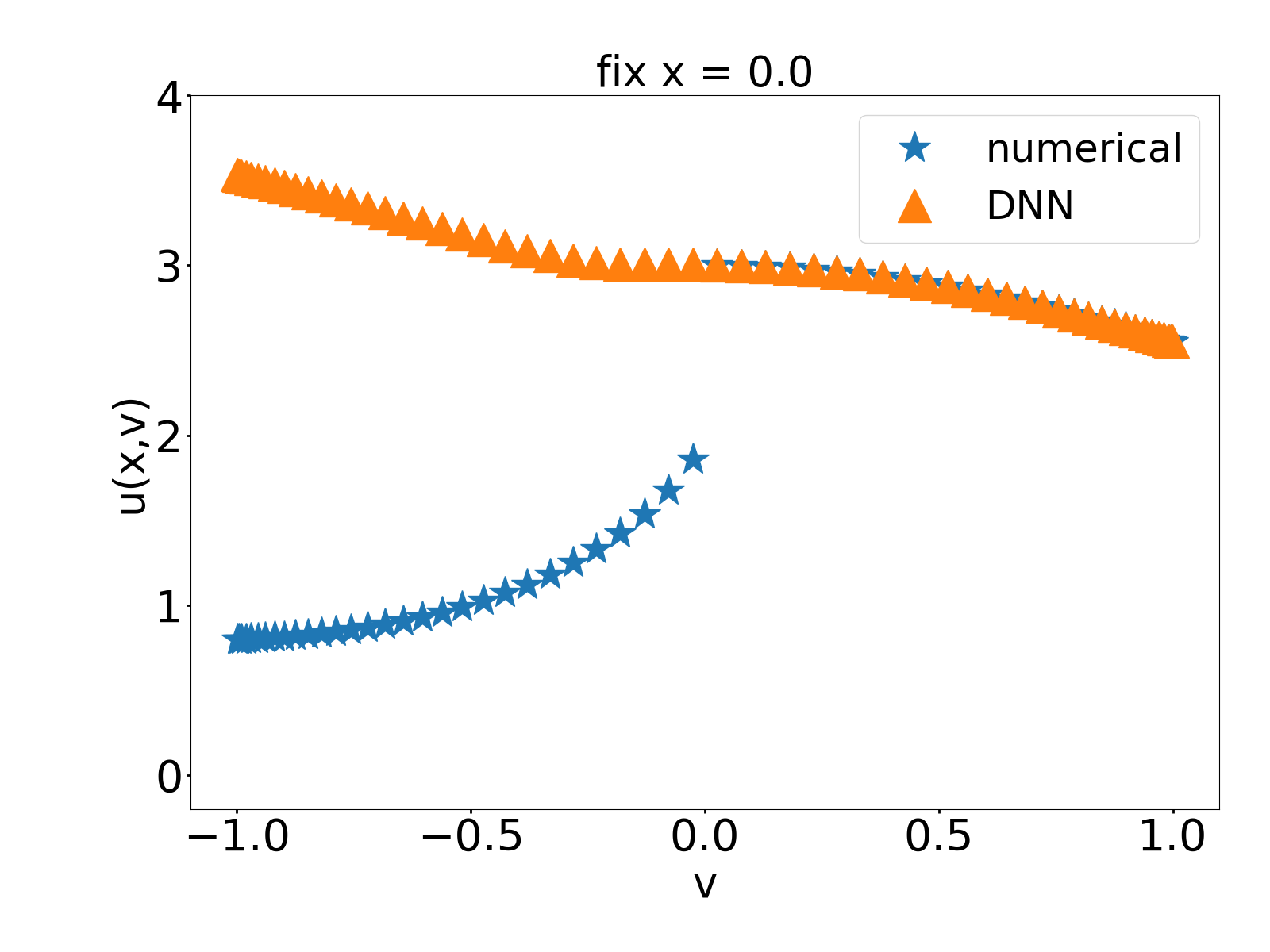}
    }
    \subfigure[ solution on $x=1$]{
        \includegraphics[scale=0.11]{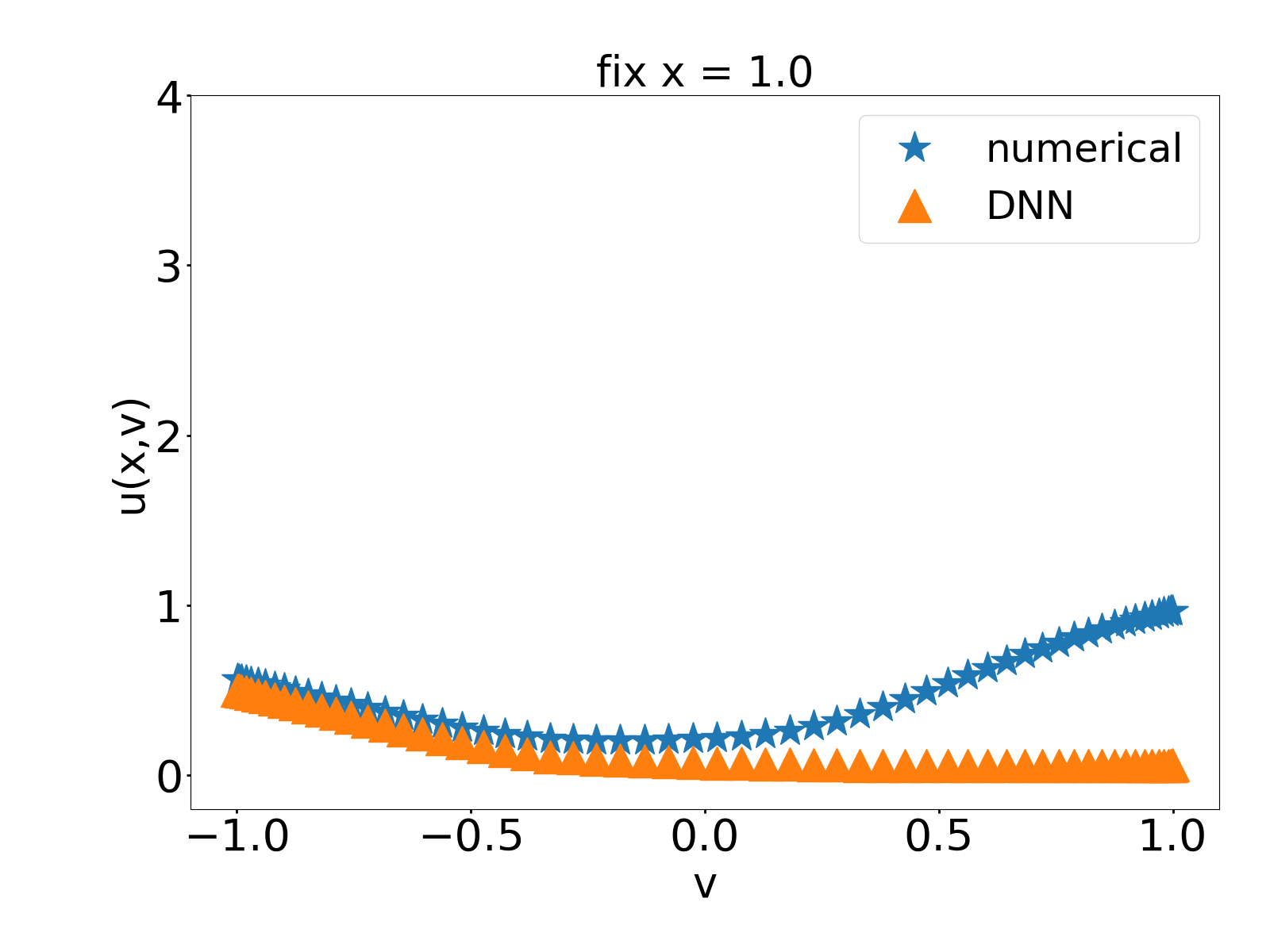}
    }
    \subfigure[ solution on $x=2$]{
        \includegraphics[scale=0.11]{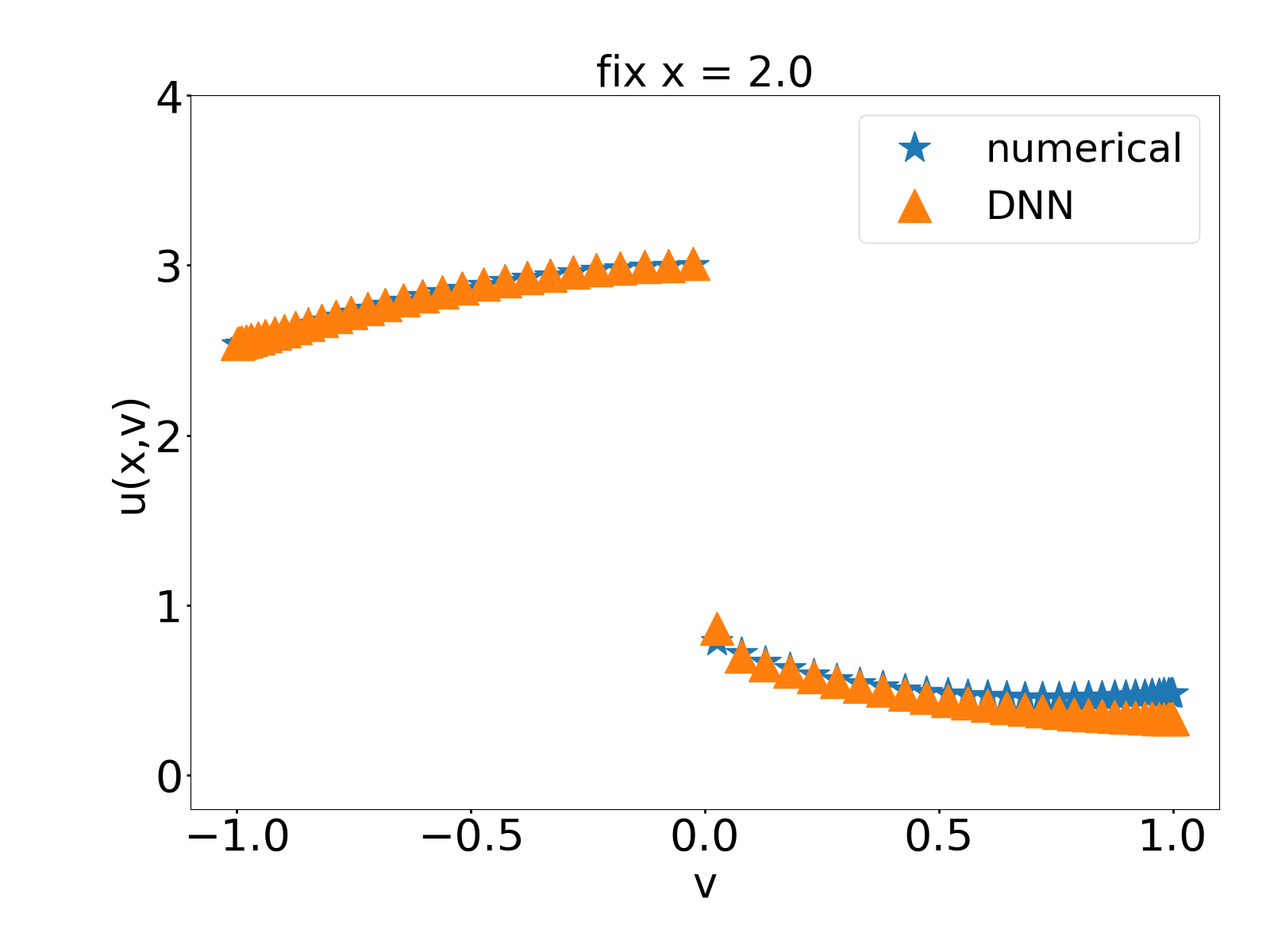}
    }
    \caption{Example 4. Trained by only the PDE,i.e., the governing equation and boundary condition. (a)
        Training loss. The  blue curve is the total risk Eq. (\ref{RTEloss}). The other three curves represent three terms in the total risk Eq.  (\ref{RTEloss}) except the last term, respectively. (b) Comparison between numerical solution and MOD-Net solution on density function. (c, d, e) Comparison between numerical solution and MOD-Net solution  on $x=0,1,2$.}
    \label{fig:bolzmann_continuous_fixed_PDE}
\end{figure}

Fix $a_2=0.01$, by TFPM method, we obtain $20$ labeled data $S^{u,k}=\{x_i,v_i,u^k(x_i,v_i)\}_{i=1}^{20}$ on coarse grids, i.e., $5$ equidistributed isometric points in $x$ direction and $4$ equidistributed isometric points in $v$ direction and $10$ labeled data $S^{\rho,k}=\{x_i,\rho^k(x_i)\}_{i=1}^{10}$, where $\rho^k(x_i) = \frac{1}{2} \sum_{\xi'\in S_v'} \omega_{\xi'} u^k(x_i, \xi')$ and $S_v'$ consists of fixed Gauss-Legendre integration points, $|S_v'|=30$, and $\omega_{\xi'}$'s are corresponding coefficients.
For this example, we train the solver of the RTE with four different loss functions.

First, we set the $\lambda_3$ in (\ref{RTEloss}) be zero, i.e., training the MOD-Net only by the governing  equation and boundary conditions, the empirical risk, i.e., training loss is shown in Fig. \ref{fig:bolzmann_continuous_fixed_PDE}(a).
To test the performance of the obtained MOD-Net solution $u_{\vtheta}$, we calculate the  density function $\rho(x)$ of $u$.
The density function of the MOD-Net solution and numerical TFPM solution with equidistributed isometric grids $(x_i,v_j)_{i\in[101],j\in[60]}$ are significantly different, as shown in Fig. \ref{fig:bolzmann_continuous_fixed_PDE}(b).
For visualization, we plot the MOD-Net solution and the corresponding numerical solution on $x=0,1,2$.
As shown in Fig. \ref{fig:bolzmann_continuous_fixed_PDE}(c, d, e), the MOD-Net solution is far from the numerical solution.

\begin{figure}
    \centering
    \subfigure[Training loss]{
        \includegraphics[scale=0.11]{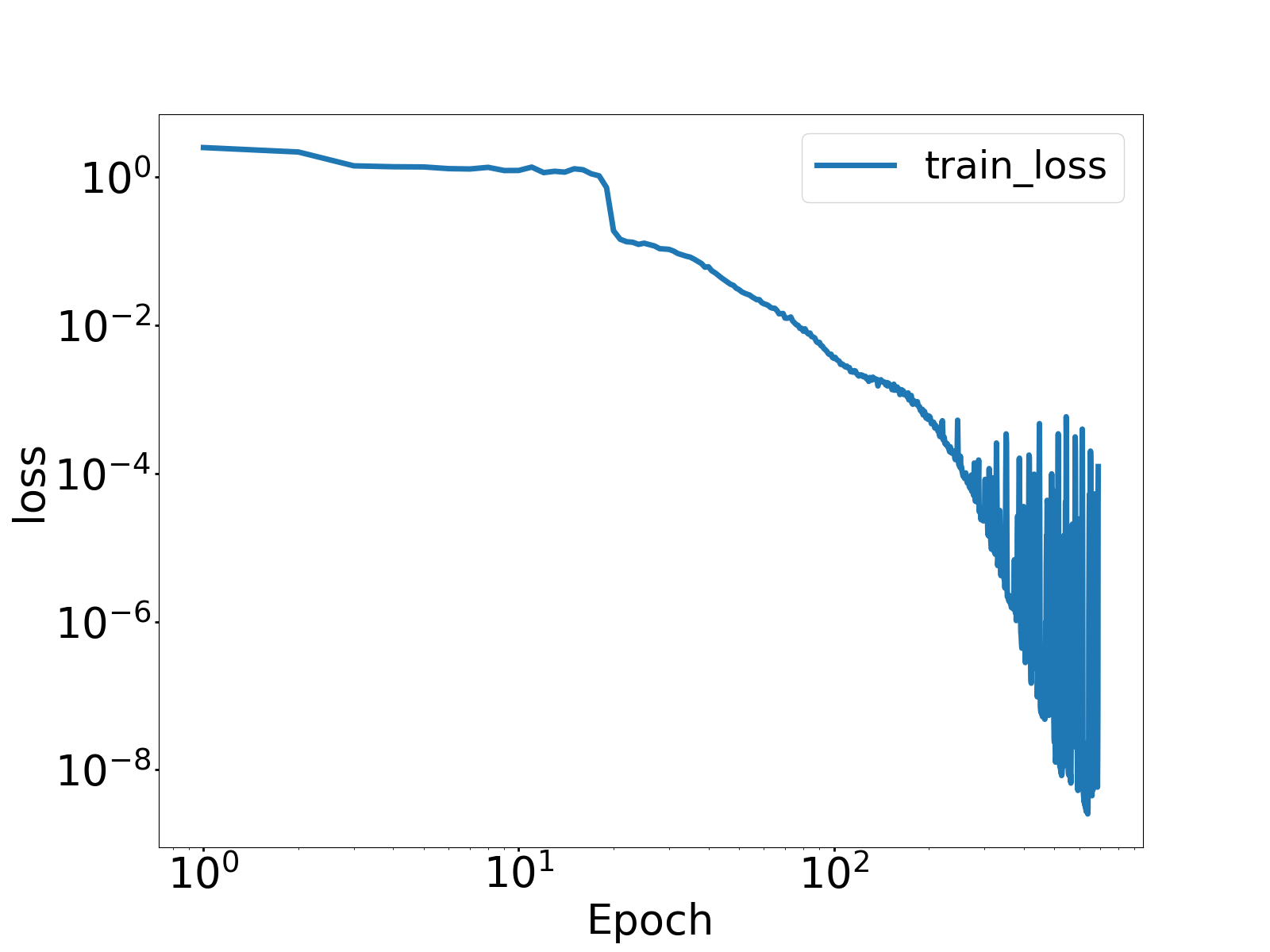}
    }
    \subfigure[Density function]{
        \includegraphics[scale=0.11]{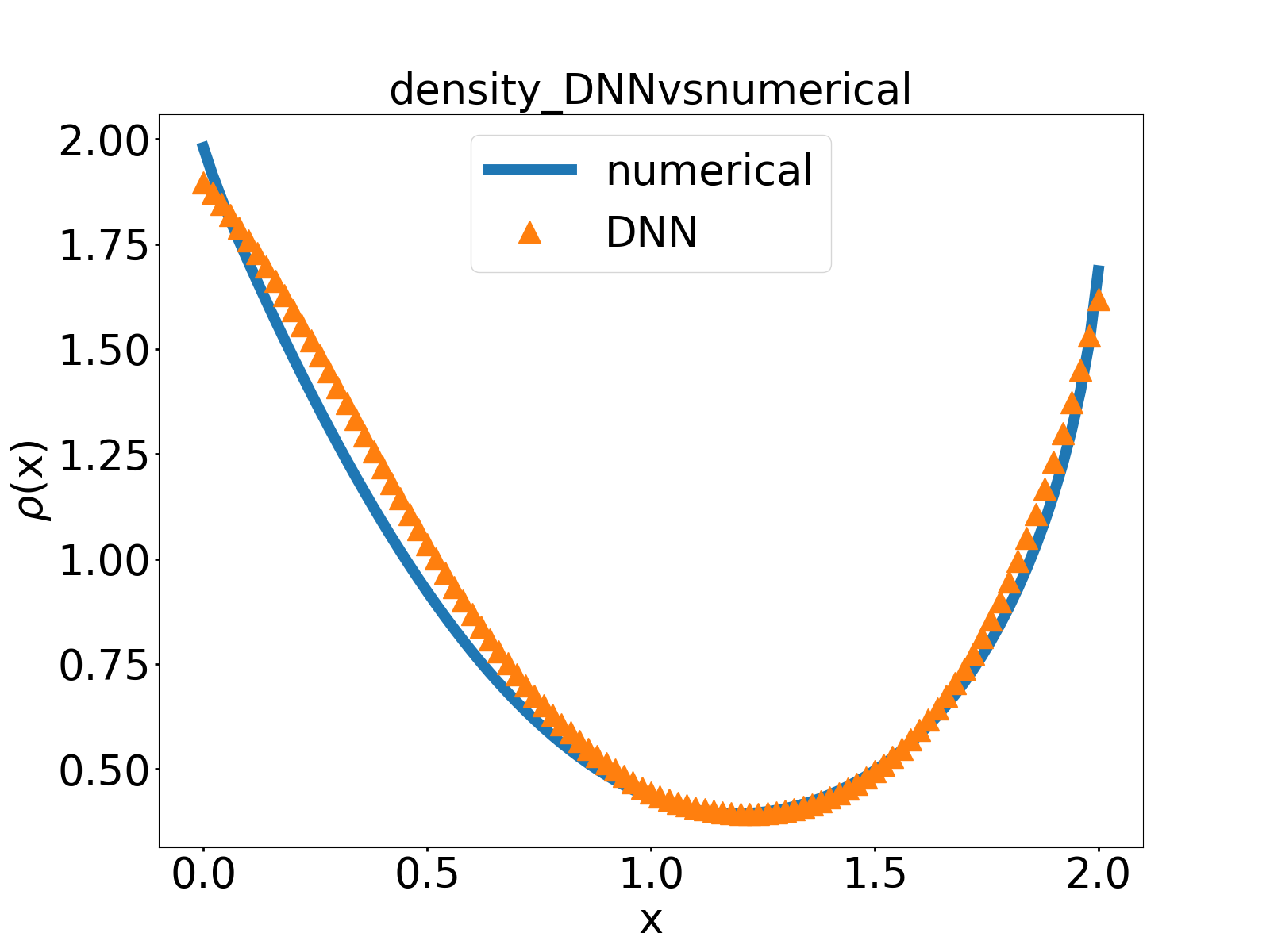}
    }
    \\
    \subfigure[solution on $x=0$]{
        \includegraphics[scale=0.11]{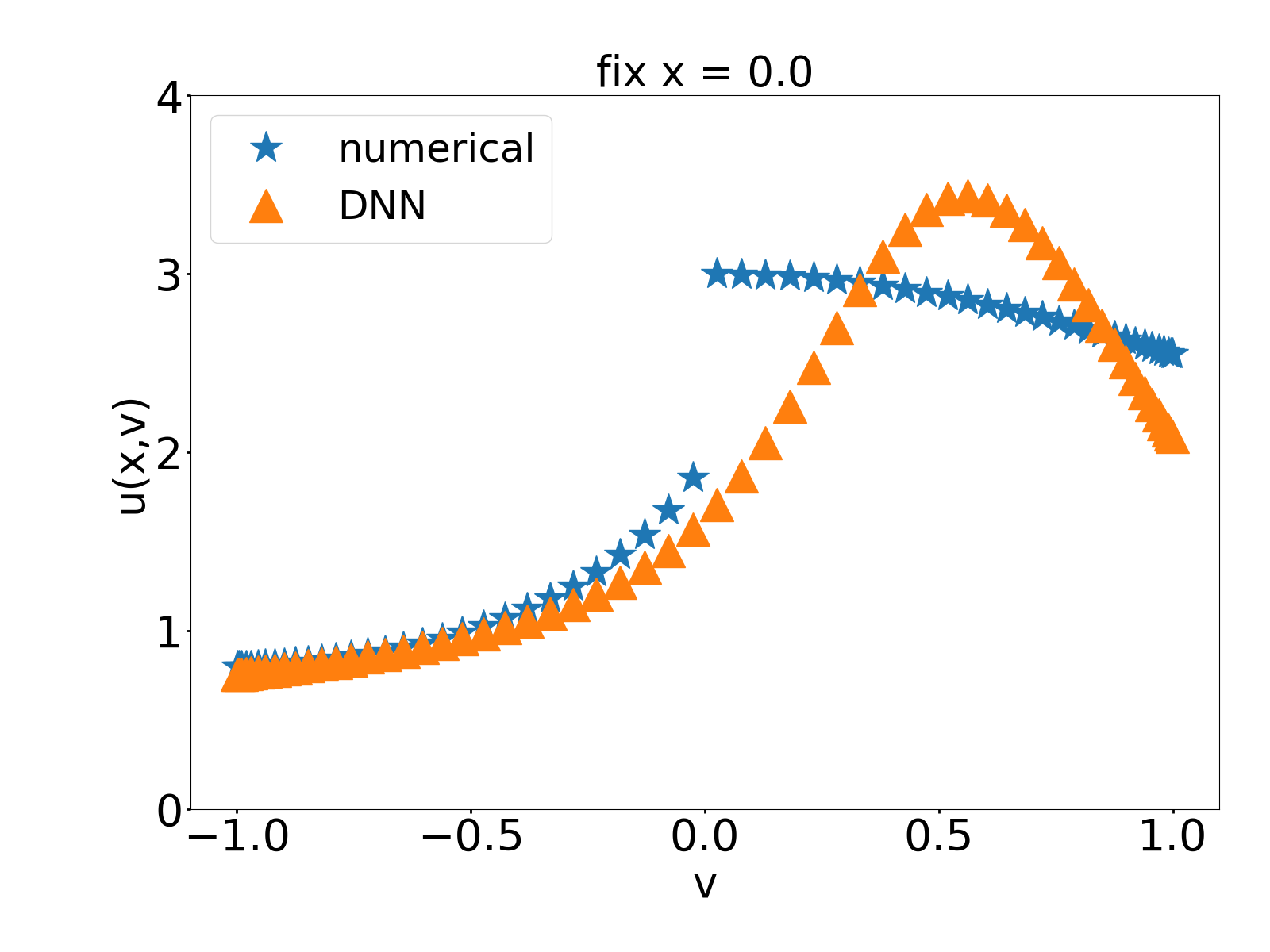}
    }
    \subfigure[solution on $x=1$]{
        \includegraphics[scale=0.11]{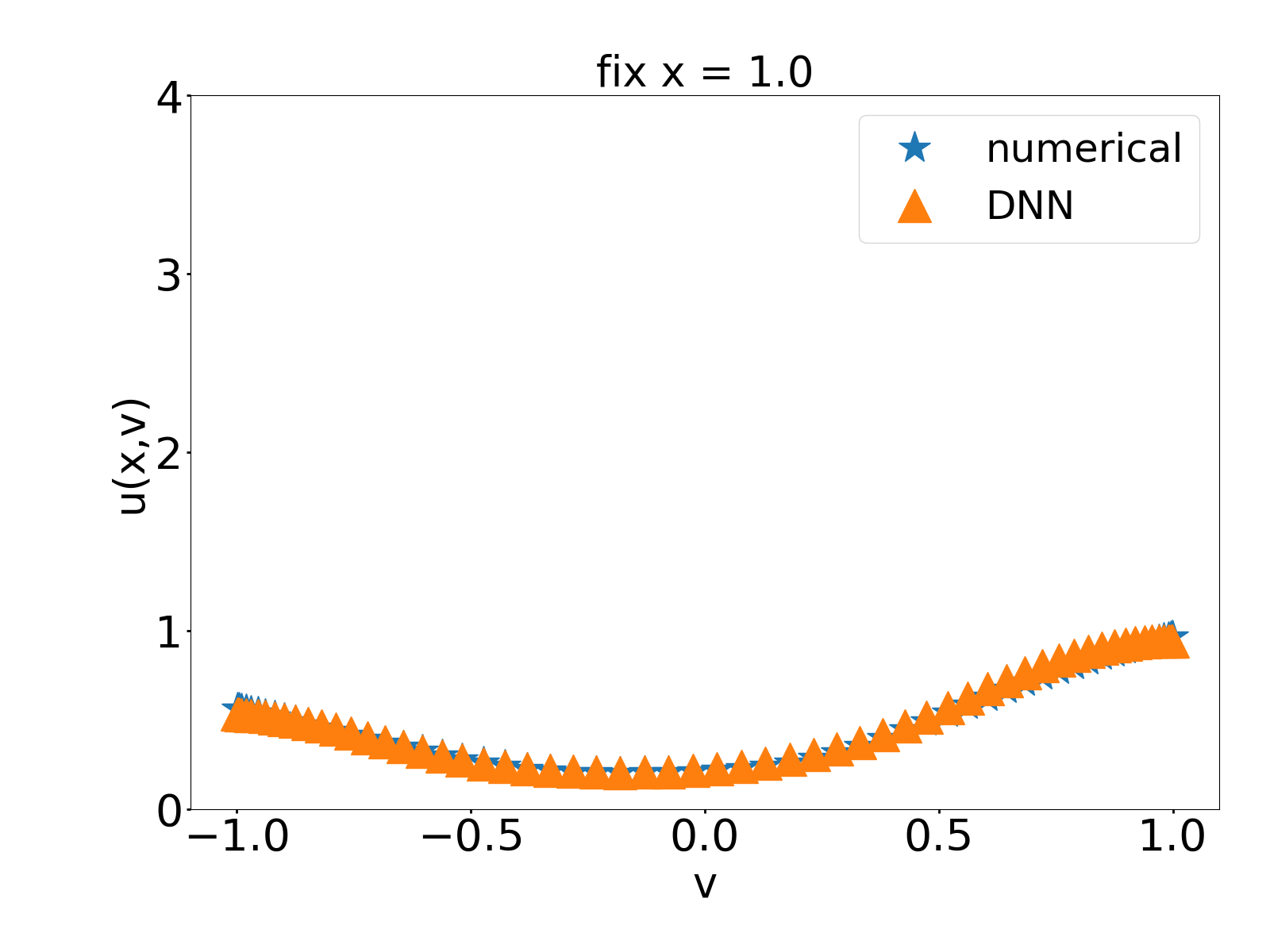}
    }
    \subfigure[solution on $x=2$]{
        \includegraphics[scale=0.11]{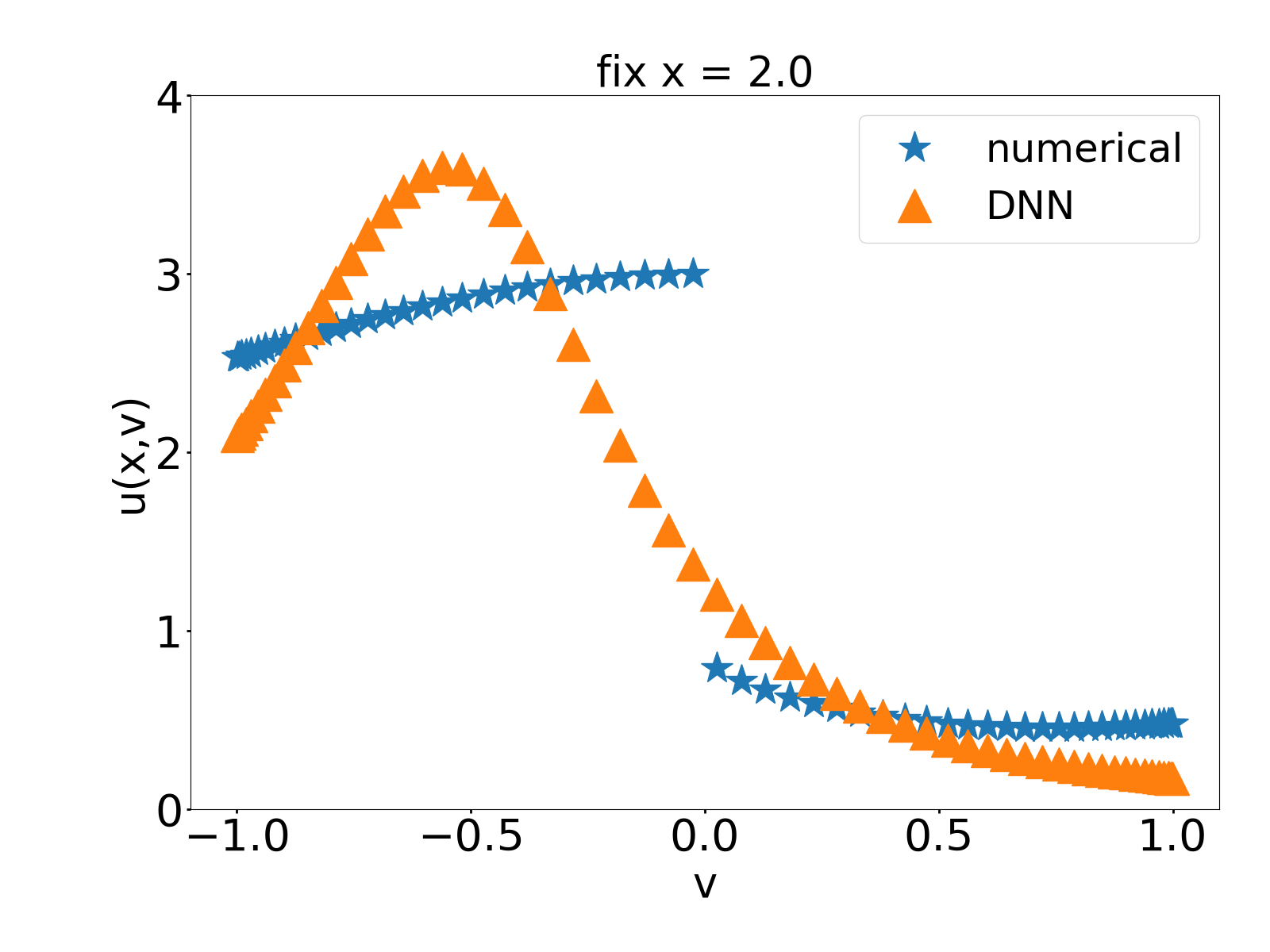}
    }
    \caption{Example 4. Trained by only the labeled data. (a) Training loss. The  blue curve is the last term in total risk Eq. (\ref{RTEloss}). The other three curves represent three terms in the total risk Eq.  (\ref{RTEloss}) except the last term, respectively.  (b) Comparison between numerical solution and MOD-Net solution on density function. (c, d, e) Comparison between numerical solution and MOD-Net solution on $x=0,1,2$.}
    \label{fig:bolzmann_continuous_fixed_data}
\end{figure}
Second,  $\lambda_1,\lambda_{21},\lambda_{22}$ in (\ref{RTEloss}) are set as zero, i.e., training the MOD-Net by only the $20$ labeled data with the training loss  shown in Fig. \ref{fig:bolzmann_continuous_fixed_data}(a).
The training loss (indicated by blue curve) decays with training epoch except oscillation appears in the final stage. We use the trick of early stopping with tolerance $50$ to obtain a DNN with small loss.
Similarly, to test the performance of the obtained MOD-Net solution $u_{\vtheta}$, we  calculate the  density function $\rho(x)$ of $u$.
The density function of the MOD-Net solution and TFPM solution have clear difference as shown in Fig. \ref{fig:bolzmann_continuous_fixed_data}(b).
For visualization, we plot  the MOD-Net solution and the corresponding numerical solution on $x=0$.
As shown in Fig. \ref{fig:bolzmann_continuous_fixed_data}(c), the MOD-Net solution significantly deviates from the numerical solution.

\begin{figure}
    \centering
    \subfigure[Training loss]{
        \includegraphics[scale=0.11]{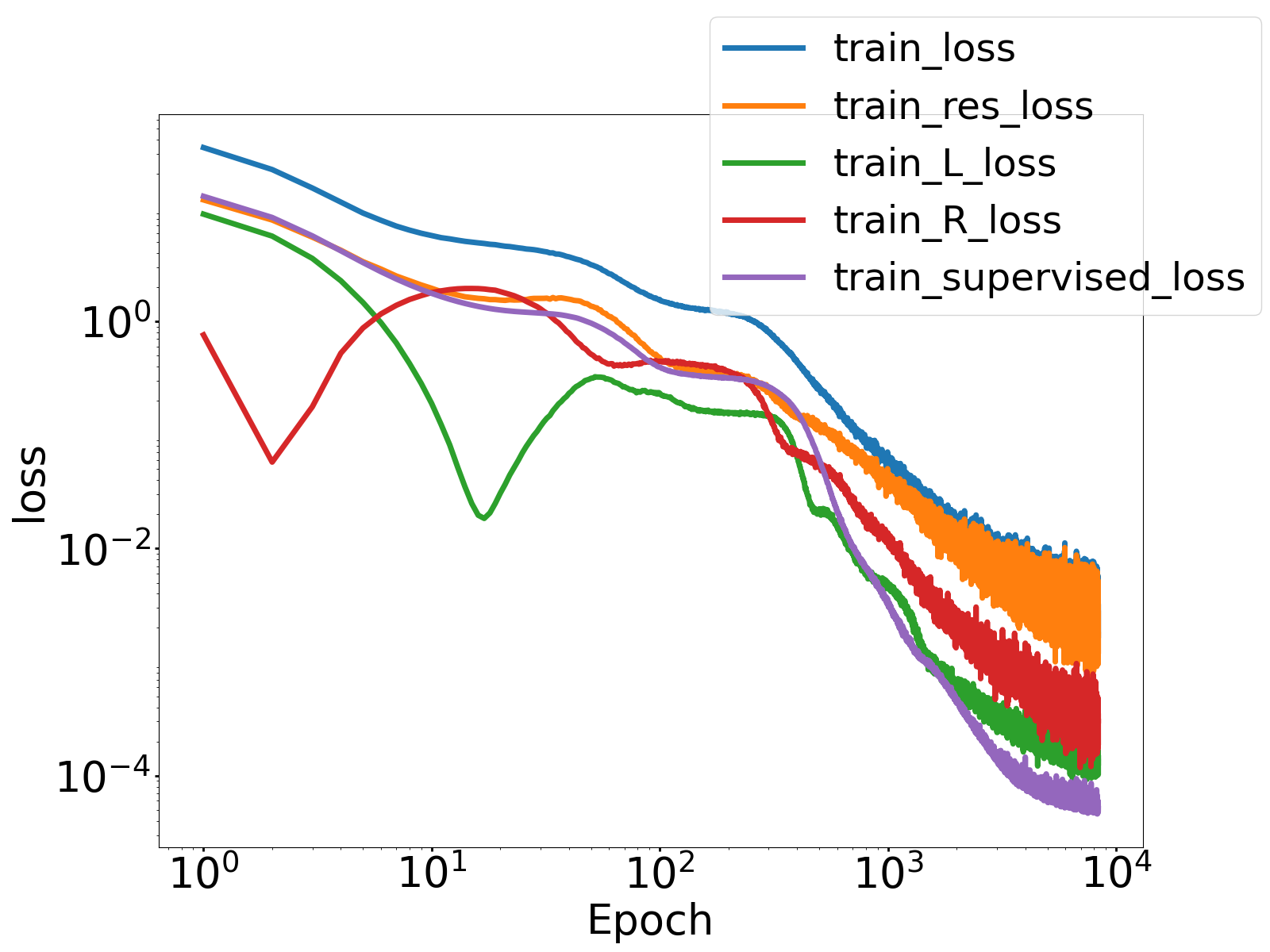}
    }
    \subfigure[Density function]{
        \includegraphics[scale=0.11]{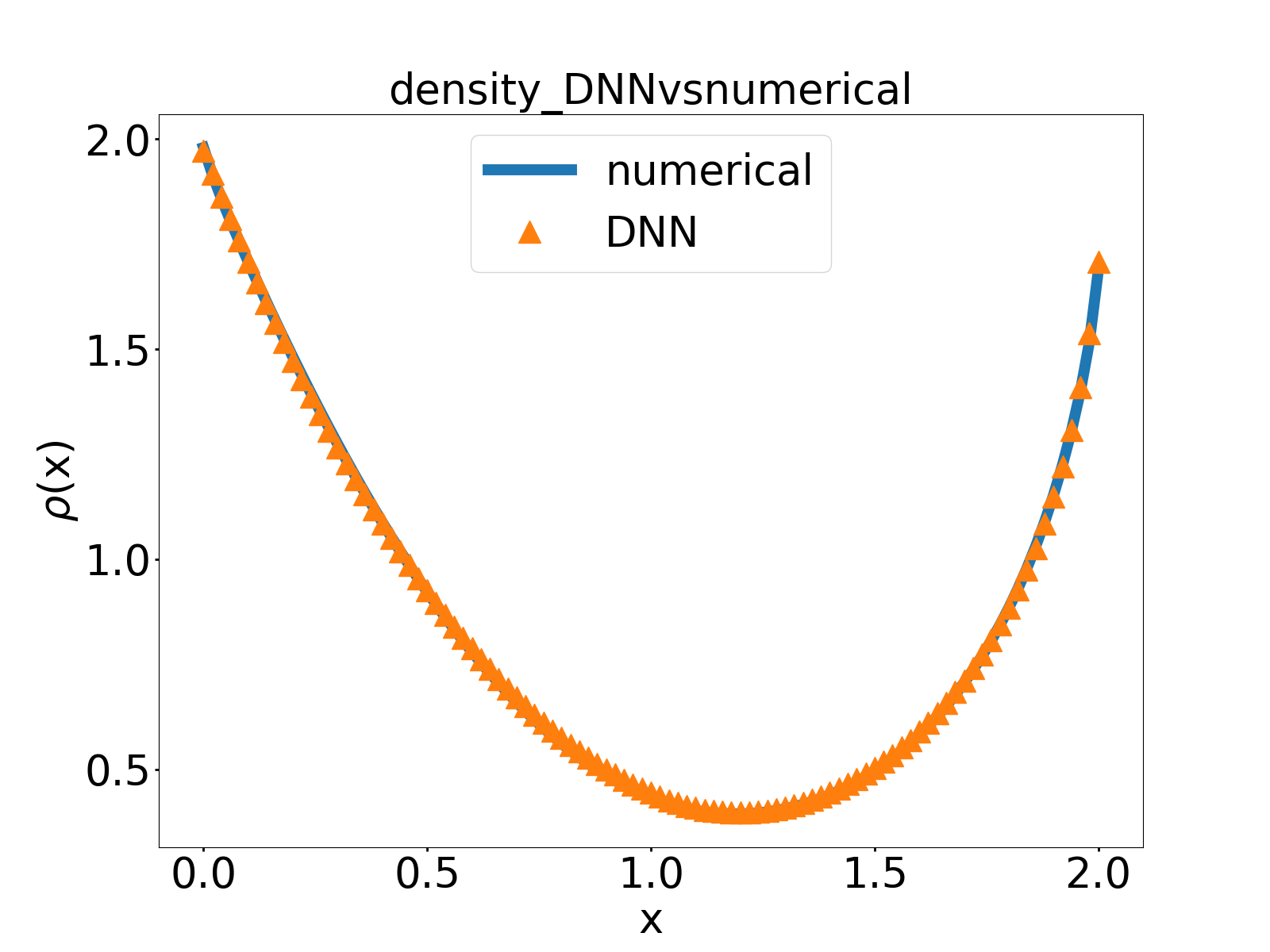}
    }
    \\
    \subfigure[solution on $x=0$]{
        \includegraphics[scale=0.11]{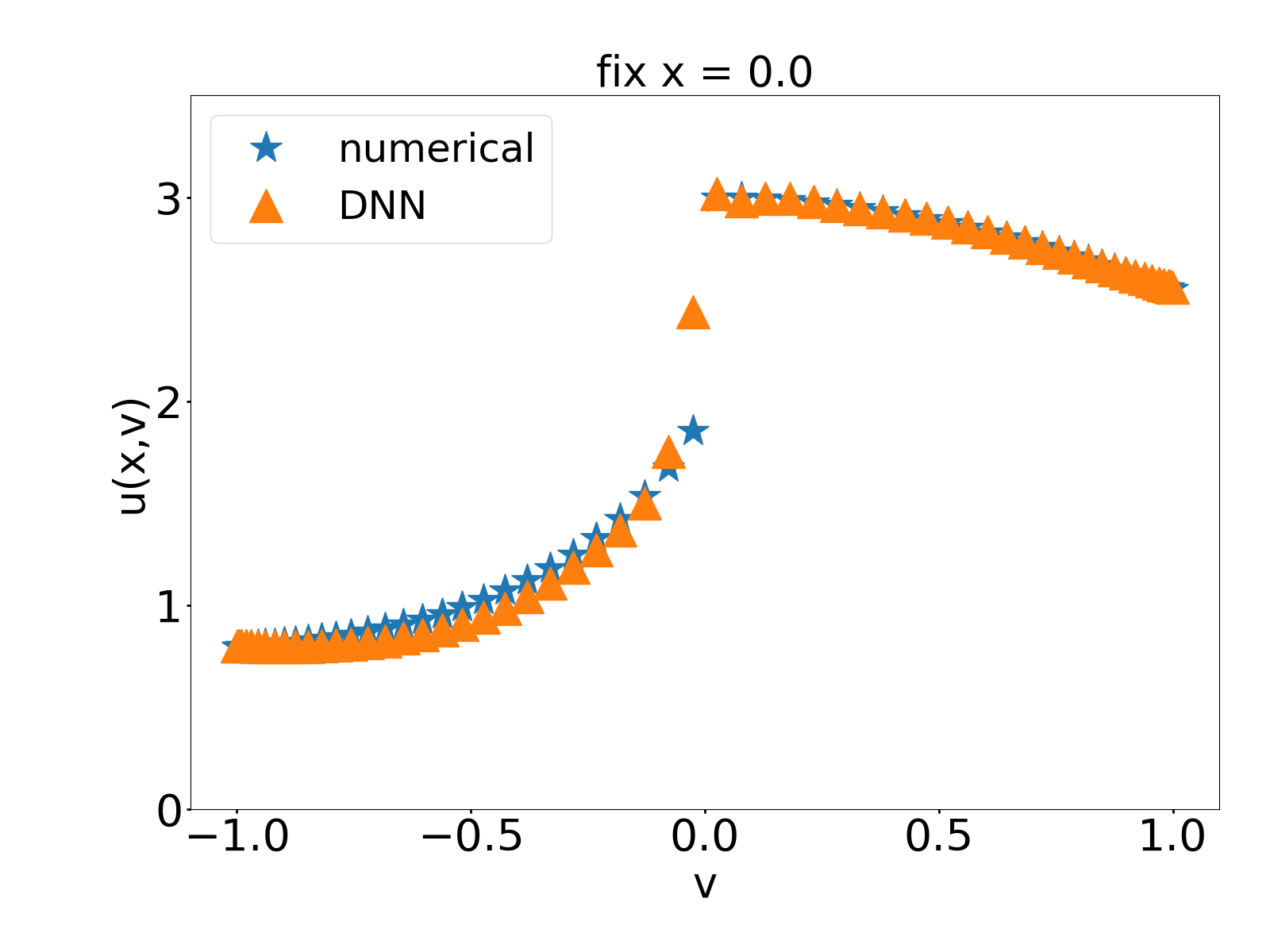}
    }
    \subfigure[solution on $x=1$]{
        \includegraphics[scale=0.11]{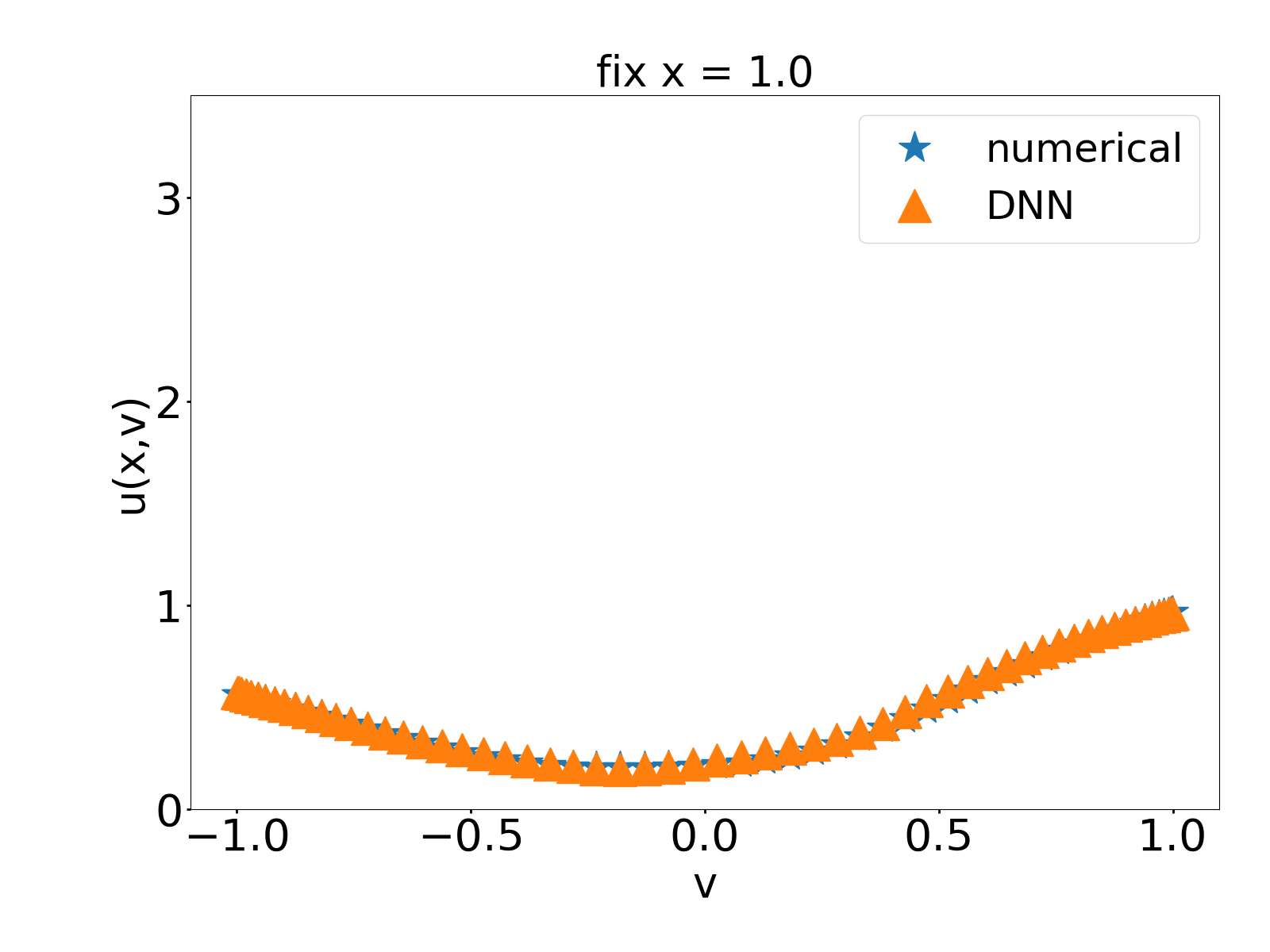}
    }
    \subfigure[solution on $x=2$]{
        \includegraphics[scale=0.11]{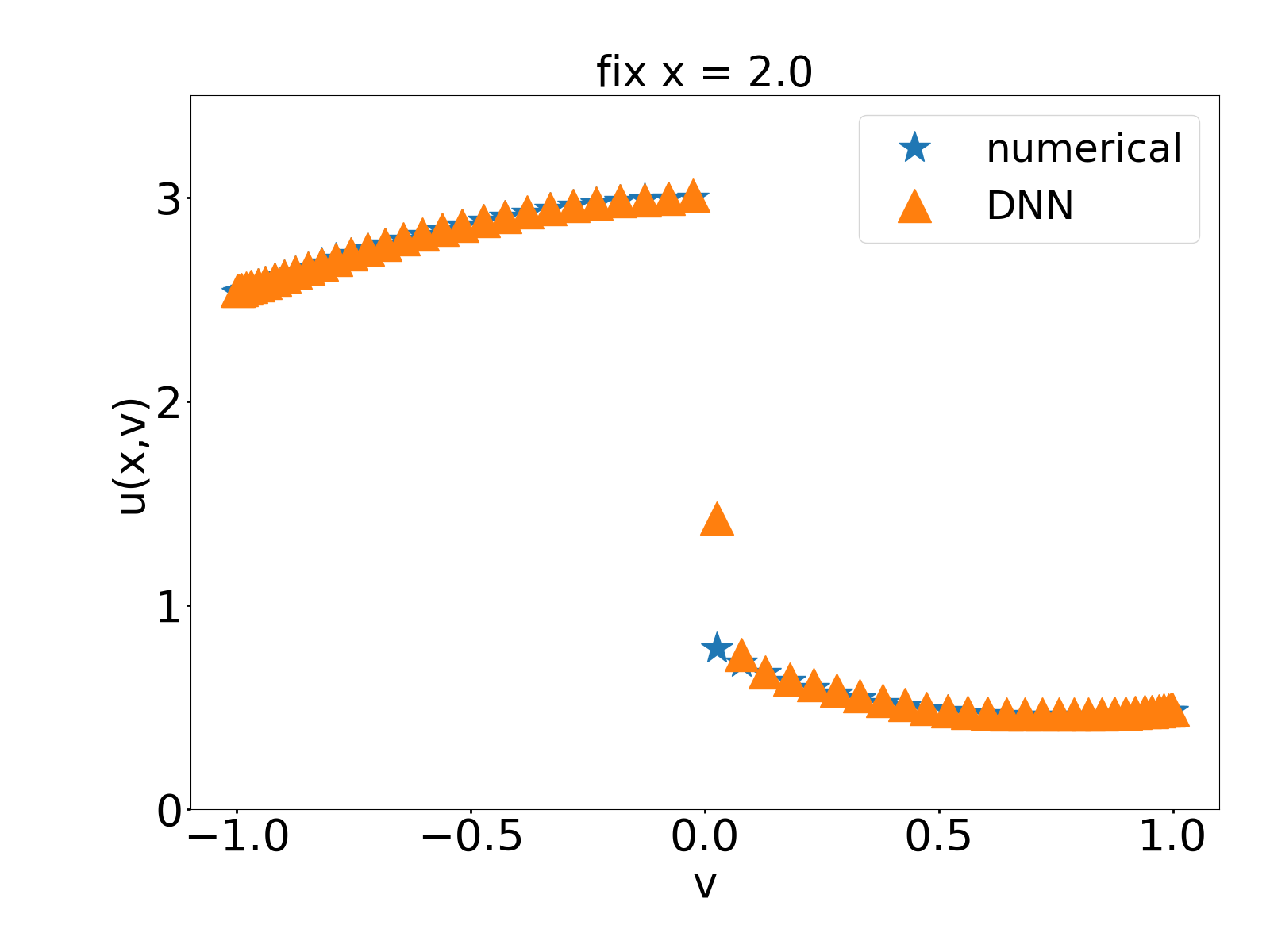}
    }
    \caption{Example 4. Trained by the governing PDE, boundary condition, and coarse grid data together. (a)
        Training loss. The  blue curve is the total risk Eq. (\ref{RTEloss}). The other four curves represent four terms in the total risk Eq.  (\ref{RTEloss}), respectively. (b) Comparison between numerical solution and MOD-Net solution on density function. (c, d, e) Comparison between numerical solution and MOD-Net solution on $x=0,1,2$.}
    \label{fig:bolzmann_continuous_fixed_PDEanddata}
\end{figure}

Third, we train the MOD-Net by the governing equation, boundary condition and $20$ labeled data $S^{u,k}$  simultaneously, and $\lambda_1,\lambda_{21},\lambda_{22},\lambda_{3}$ in (\ref{RTEloss}) are set as one. The training total loss (\ref{RTEloss}) decays with training epoch as shown in Fig. \ref{fig:bolzmann_continuous_fixed_PDEanddata}(b).
Similarly, we  calculate the  density function of  MOD-Net solution and TFPM solution. They are very consistent as shown in Fig. \ref{fig:bolzmann_continuous_fixed_PDEanddata}(c).
As shown in Fig. \ref{fig:bolzmann_continuous_fixed_PDEanddata}(d, e, f), the MOD-Net solution overlap  with the numerical solution  very well except for the discontinuous points.

\begin{figure}[!htb]
    \centering
    \subfigure[Training loss]{
        \includegraphics[scale=0.11]{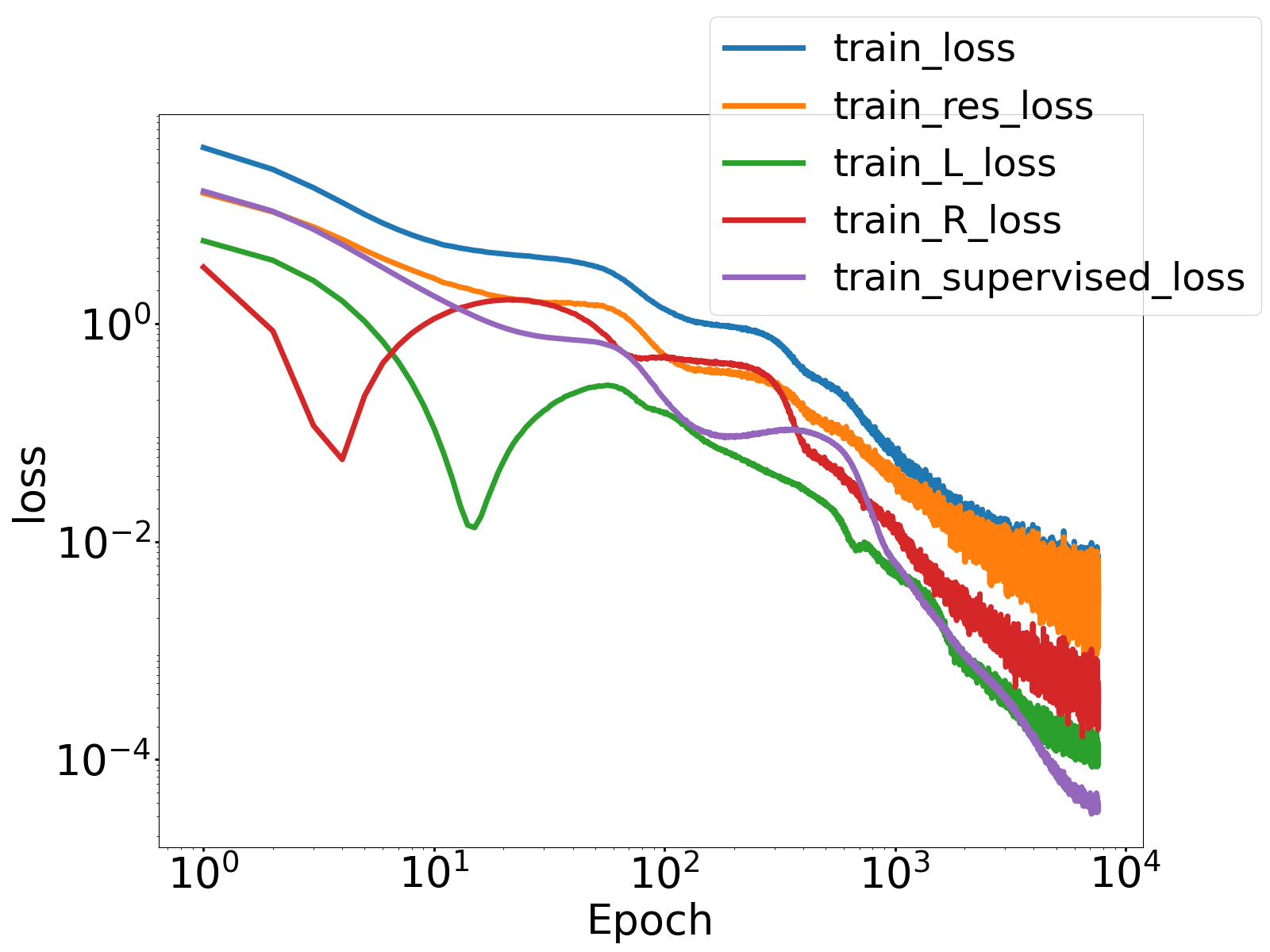}
    }
    \subfigure[Density function]{
        \includegraphics[scale=0.11]{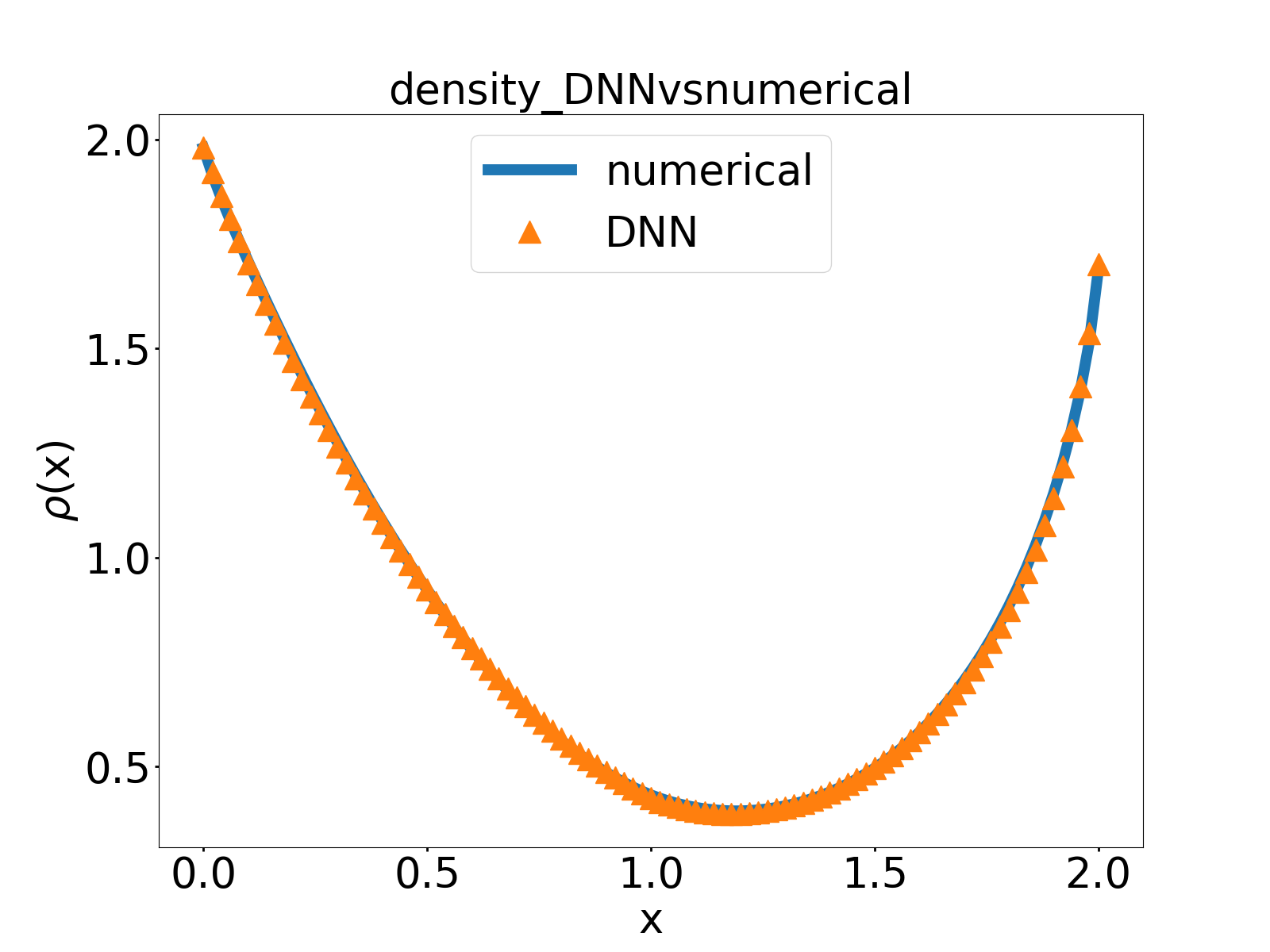}
    }
    \\
    \subfigure[solution on $x=0$]{
        \includegraphics[scale=0.11]{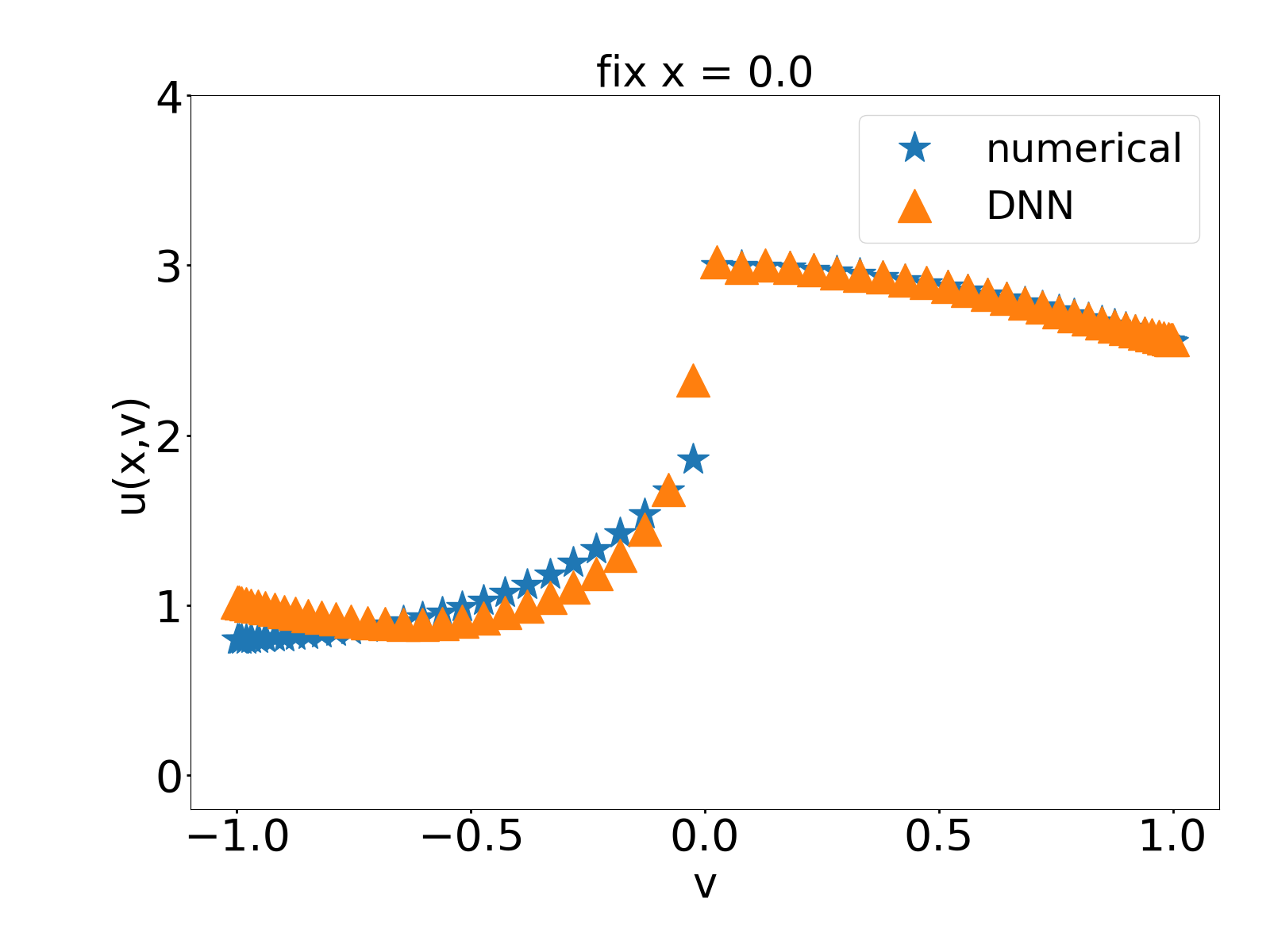}
    }
    \subfigure[solution on $x=1$]{
        \includegraphics[scale=0.11]{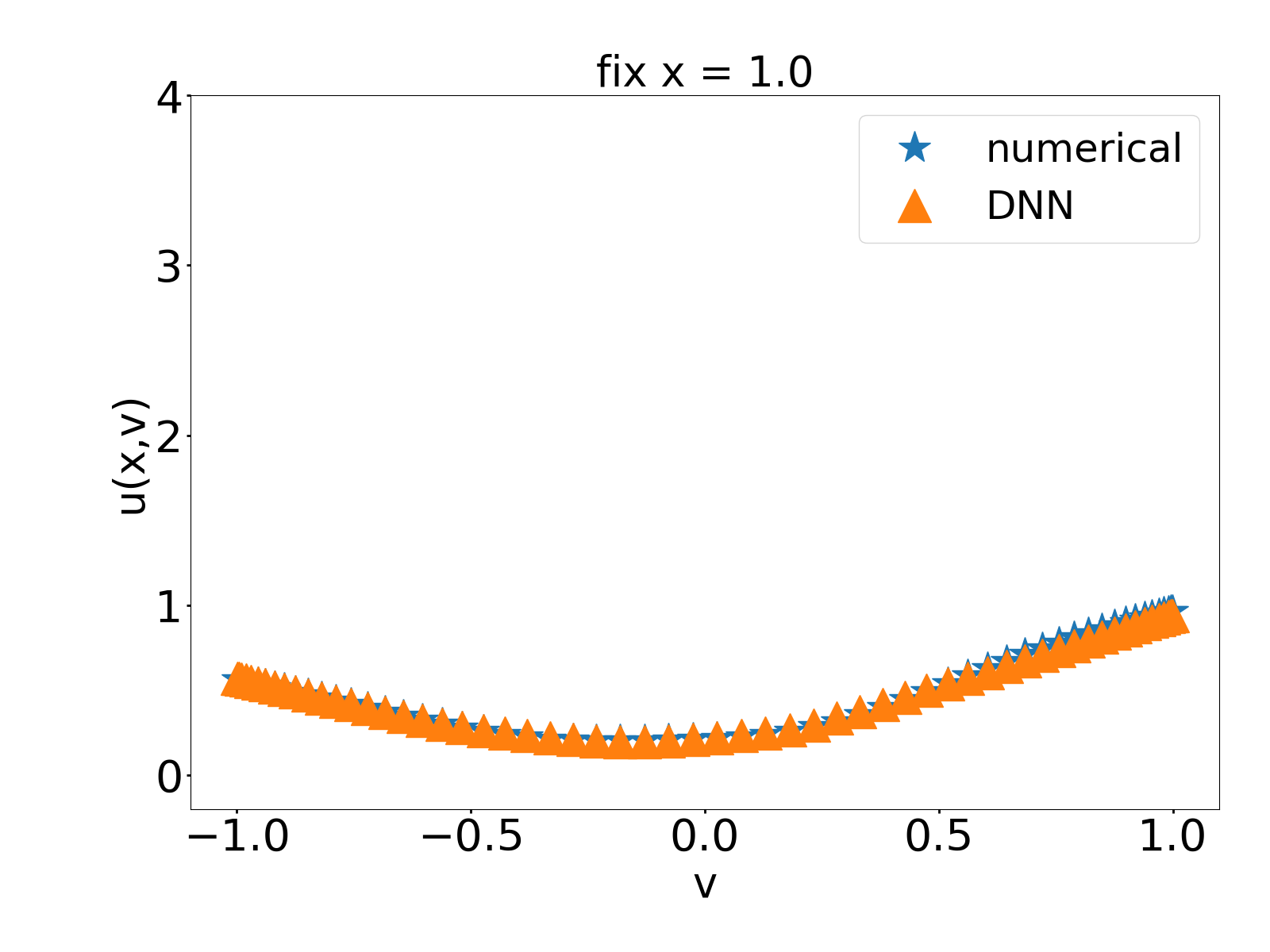}
    }
    \subfigure[solution on $x=2$]{
        \includegraphics[scale=0.11]{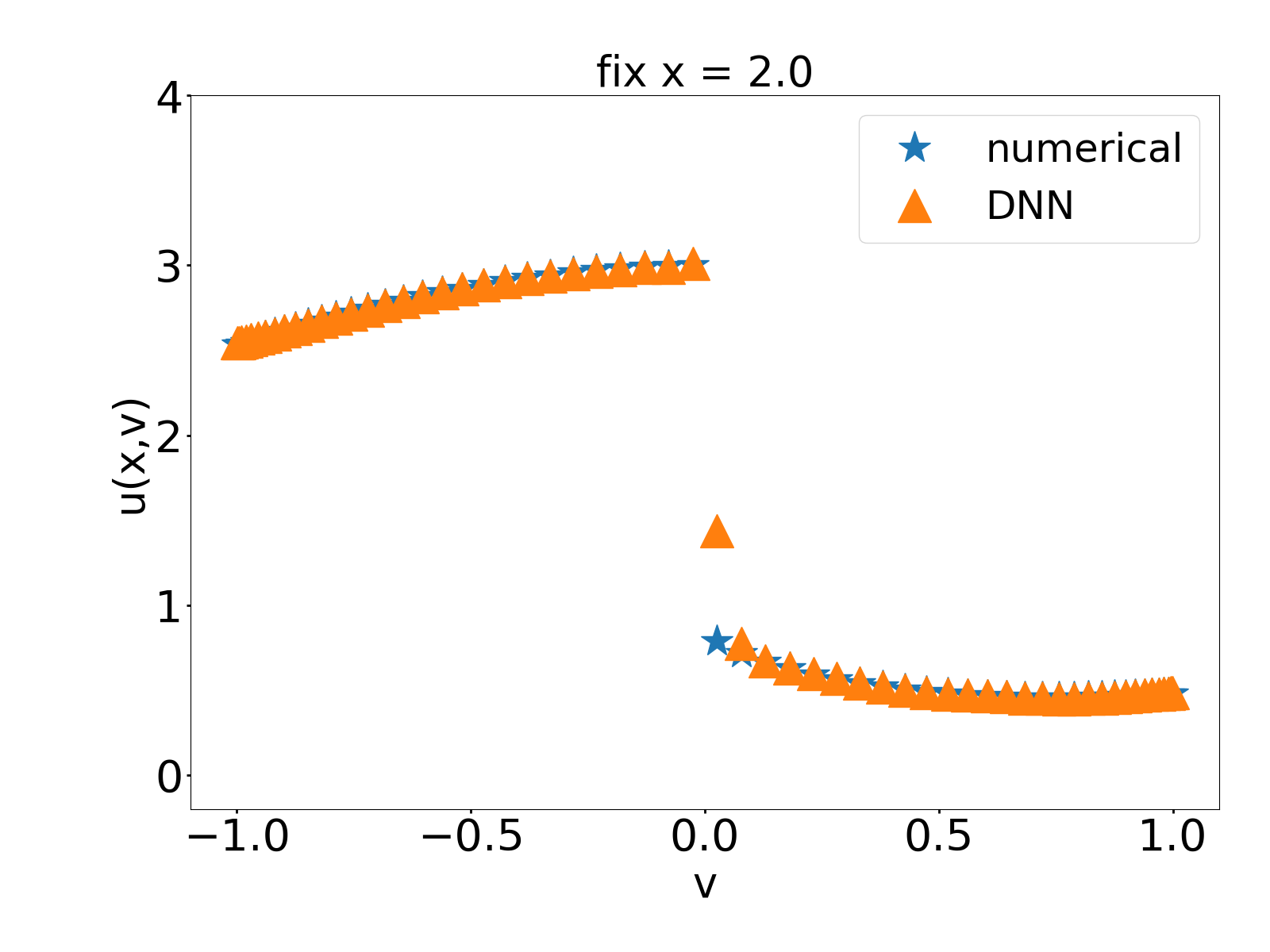}
    }
    \caption{Example 4. Trained by the governing PDE, boundary condition, and coarse density function $\rho$ together. (a)
        Training loss. The  blue curve is the total risk Eq. (\ref{RTElossrho}). The other four curves represent four terms in the total risk Eq.  (\ref{RTElossrho}), respectively. (b) Comparison between numerical solution and MOD-Net solution on density function. (c, d, e) Comparison between numerical solution and MOD-Net solution on $x=0,1,2$.}
    \label{fig:bolzmann_continuous_fixed_PDEandrho}
\end{figure}

Four, we train the MOD-Net by the governing equation, boundary condition and $10$ labeled data $S^{\rho,k}$ simultaneously, and $\lambda_1,\lambda_{21},\lambda_{22},\lambda_{3}$ in (\ref{RTElossrho}) are set as one. The training total loss (\ref{RTElossrho}) decays with training epoch as shown in Fig. \ref{fig:bolzmann_continuous_fixed_PDEandrho}(b).
To see the performance of the trained MOD-Net, we  calculate the  density function of  MOD-Net solution and TFPM solution. They are very consistent as shown in Fig. \ref{fig:bolzmann_continuous_fixed_PDEandrho}(c).
As shown in Fig. \ref{fig:bolzmann_continuous_fixed_PDEandrho}(d, e, f), the MOD-Net solution overlap  with the numerical solution  well except for the discontinuous points.

To sum up, with the information of only PDE or only a few labeled data, the MOD-Net cannot be trained well, however,
combining these two information,
we can train the MOD-Net very well. Note that the labeled data is not restricted in the solution, the information of the density function can also benefit the training of MOD-Net.


To solve the PDE with various $a_2$ accurately and quickly, we use MOD-Net approach to learn a neural operator. For this case, we train MOD-Net with PDE and a few labeled data simultaneously. For each $a_2$ in $\{0.03k\}_{k=1}^{33}$ , by TFPM method, we calculate the corresponding numerical solution on $20$ fixed coarse grids , i.e., $5$ equidistributed isometric points in $x$ direction and $4$ equidistributed isometric points in $v$ direction, as labeled data and obtain $S^{u,k}=\{x_i,v_i,u^k(x_i,v_i)\}_{i=1}^{20}$.

We set $K=20$ and for each epoch, we sample $K$ different $a_2$'s uniformly from  $\{0.03k\}_{k=1}^{33}$.
In training,  we set the number of integration points $|S_{G_v}^+|=30$,$|S_{G_v}^-|=30$,$|S_{v}|=30$, and for each epoch, for each $k$,
we uniformly sample $500$ points in $[0,2]\times[-1,1]$ for $S^{\Omega,S,k}$, $500$ points $v_i$'s in $[0,1]$ for $S^{\partial \Omega_L,k}=\{(0,v_i)|i=1,\cdots,500\}$, another $500$ points $v_i$'s in $[-1,0]$ for $S^{\partial \Omega_R,k}=\{(2,v_i)|i=1,\cdots,500\}$ and use the $20$ labeled data in $S^{u,k}$ .
The training loss is shown in Fig. \ref{fig:bolzmann_continuous_op}(b). The total risk (indicated by blue curve) decays with training epoch. For visualization of each part of the total risk, we also display the residual loss of the RTE (indicated by orange curve), the loss of two boundary lines (indicated by green and red curves) and the loss of supervised data (indicated by remaining curve), which  decay with the training epoch overall.

To test the performance of the trained MOD-Net on $a_2=0.02,0.11,0.51,0.99$.
For each $a_2$, we calculate the  density function $\rho(x)$ of  the solution.
For all $a_2$'s, the density function of MOD-Net solution and numerical solution by TFPM method with equidistributed isometric grids $(x_i,v_j)_{i\in[101],j\in[60]}$ are very consistent as shown in Fig. \ref{fig:bolzmann_continuous_op}.

\section{Numerical experiments: $1$D Burgers equation} \label{sec:Burgers}

To show the performance of the MOD-Net for a nonlinear PDE, first we consider the $1$D Burgers equation in the steady state,
\begin{equation}
    \begin{aligned}
         & \partial_{x}(\frac{1}{2}{u^2(x)})=\nu\partial_{xx}u(x)+g(x), \quad x\in \Omega,
        \\
         & u(x)=  \phi(x), \quad x\in \partial\Omega,
    \end{aligned}
    \label{Burgers}
\end{equation}

For illustration,  we take $\Omega=[-1,1]$,$\nu=1$,$g(x)=0$ and $\phi(x)=c_1\cos(k_1x) + c_2\sin(k_2x) $ for all $x\in \Omega$,
\begin{equation}
    \begin{aligned}
         & \partial_{x}(\frac{1}{2}{u^2(x)})=\partial_{xx}u(x), \quad x\in \Omega,
        \\
         & u(-1)=c_1\cos(-k_1) + c_2\sin(-k_2),
        \\
         & u(1)=c_1\cos(k_1) + c_2\sin(k_2),
    \end{aligned}
    \label{Burgers example}
\end{equation}
in which $c_1, c_2, k_1, k_2$ control the boundary condition.
For this problem , it is difficult to obtain the analytical solution, but a lot of traditional numerical schemes can be used to solve it. In this paper, we use the upwind scheme.

\subsection{Use DNN to fit nonlinear operator}
For a linear PDE, Green's function can help us obtain the solution of PDE due to the superposition principle.
For a nonlinear PDE, similarly we use the following representation of the solution of PDE (\ref{Burgers}),

\begin{equation}
    u(x;\phi,g) =F(\int_{\Omega} G_1(x,x')g(x') \diff{x'}
    +\int_{\partial \Omega} G_2(x,x')\phi(x') \diff{x'}),
\end{equation}

In the considered  example, $\Omega=[-1,1]$ and $g(x)=0$, we have
\begin{equation}
    \begin{aligned}
        u(x;\phi) & =F(\int_{\partial \Omega} G_2(x,x')\phi(x')
        \diff{x'})                                              \\
                  & =F( G_2(x,-1)\phi(-1)+G_2(x,1)\phi(1)
        ).
    \end{aligned}
\end{equation}
We use a  DNN of hidden layer size $256$-$256$-$256$-$256$ equipped with sigmoid function as activation function  $G_{\vtheta_L}(x)$ to fit $G_2(x,-1)$, a DNN with the same setting $G_{\vtheta_R}(x)$ to fit $G_2(x,1)$ and a DNN of one hidden layer, the size of which is $256$, equipped with sigmoid function as activation function  $F_{\vtheta}(x)$ .
Then we can represent neural operator $u_{\vtheta}(x;g)$ with  DNN $G_{\vtheta_L}(x)$, $G_{\vtheta_R}(x)$ and $F_{\vtheta}(x)$, that is,
\begin{equation}
    \begin{aligned}
        u_{\vtheta}(x;\phi)
         & =F_{\vtheta}( G_{\vtheta_L}(x)\phi(-1)+G_{\vtheta_R}(x)\phi(1)).
    \end{aligned}
    \label{burgers_op_rep}
\end{equation}

\subsection{Empirical risk function}
For this simple example, to learn a neural operator, which can approximate the operator $\fG:\phi \mapsto u$, we use no labeled data and only utilize the information of PDE, i.e., governing equation and boundary condition, $\phi^k$, $k=1,2,\cdots,K$.

To utilize the  constraint of governing equation of PDE, we uniformly sample a set of data from $\Omega=[-1,1]$, i.e., $S^{\Omega,k}$.
Since the boundary $\partial\Omega$ consists of two points, i.e. $-1$ and $1$,  then we directly use the information of these two points.

The empirical risk for this example is as follows,
\begin{equation}
    \begin{aligned}
        \RS
         & =  \frac{1}{K}\sum_{k\in[K]}  \bigg(\lambda_1
        \frac{1}{|S^{\Omega,k}|}\sum_{x \in  S^{\Omega,k}}\Big(
        \partial_{x}(\frac{1}{2}{u^2_{\vtheta}(x;\phi^k)})-\partial_{xx}u_{\vtheta}(x;\phi^k)
        \Big)^{2}                                                             \\
         & \quad +  \lambda_2  \Big( u_{\vtheta}(-1;\phi^k)-\phi^k(-1)) \Big)
        \\
         & \quad +  \lambda_3 \Big( u_{\vtheta}(1;\phi^k)-\phi^k(1)) \Big)
        \bigg).
        \label{burgers_leastsquareloss}
    \end{aligned}
\end{equation}

\subsection{Learning process}

For each training epoch, we first randomly choose source functions denoted by $\{\phi^k\}_{k=1}^{K}$. Second, we randomly sample data and obtain data set $S^{\Omega,k}$.
In the following, we feed the data into the neural networks $G_{\vtheta_L}(x)$, $G_{\vtheta_R}(x)$ and calculate the total risk (\ref{burgers_leastsquareloss})  with neural operator $u_{\vtheta}(x;g)$, see Eq. (\ref{burgers_op_rep}).
Train neural networks $G_{\vtheta_L}(x)$, $G_{\vtheta_R}(x)$, and $F_{\vtheta}(x)$ with Adam to minimize the total risk, and finally we  obtain well-trained  DNNs $G_{\vtheta_L}(x)$, $G_{\vtheta_R}(x)$ and $F_{\vtheta}(x)$, furthermore, according to Eq. (\ref{burgers_op_rep}), we obtain a neural operator $u_{\vtheta}(x;\phi)$.

\subsection{Results}
The boundary condition $\phi(x)=c_1\cos(k_1x) + c_2\sin(k_2x)$ is determined by $c_1,c_2,k_1,k_2$. To illustrate our approach,
we  train a neural operator mapping from $\phi(x)$ to the solution $u(x)$ of the Burgers equation (\ref{Burgers example}).

\paragraph{Example 5: MOD-Net for Burgers equation}

To solve the PDE with various $c_1,c_2,k_1,k_2$ accurately and quickly, we use MOD-Net approach to learn a neural operator.
For simplicity, we fixed $c_1=4,k_1=2,k_2=10$ and only change $c_2$.

\begin{figure}[!htb]
    \centering
    \includegraphics[scale=0.11]{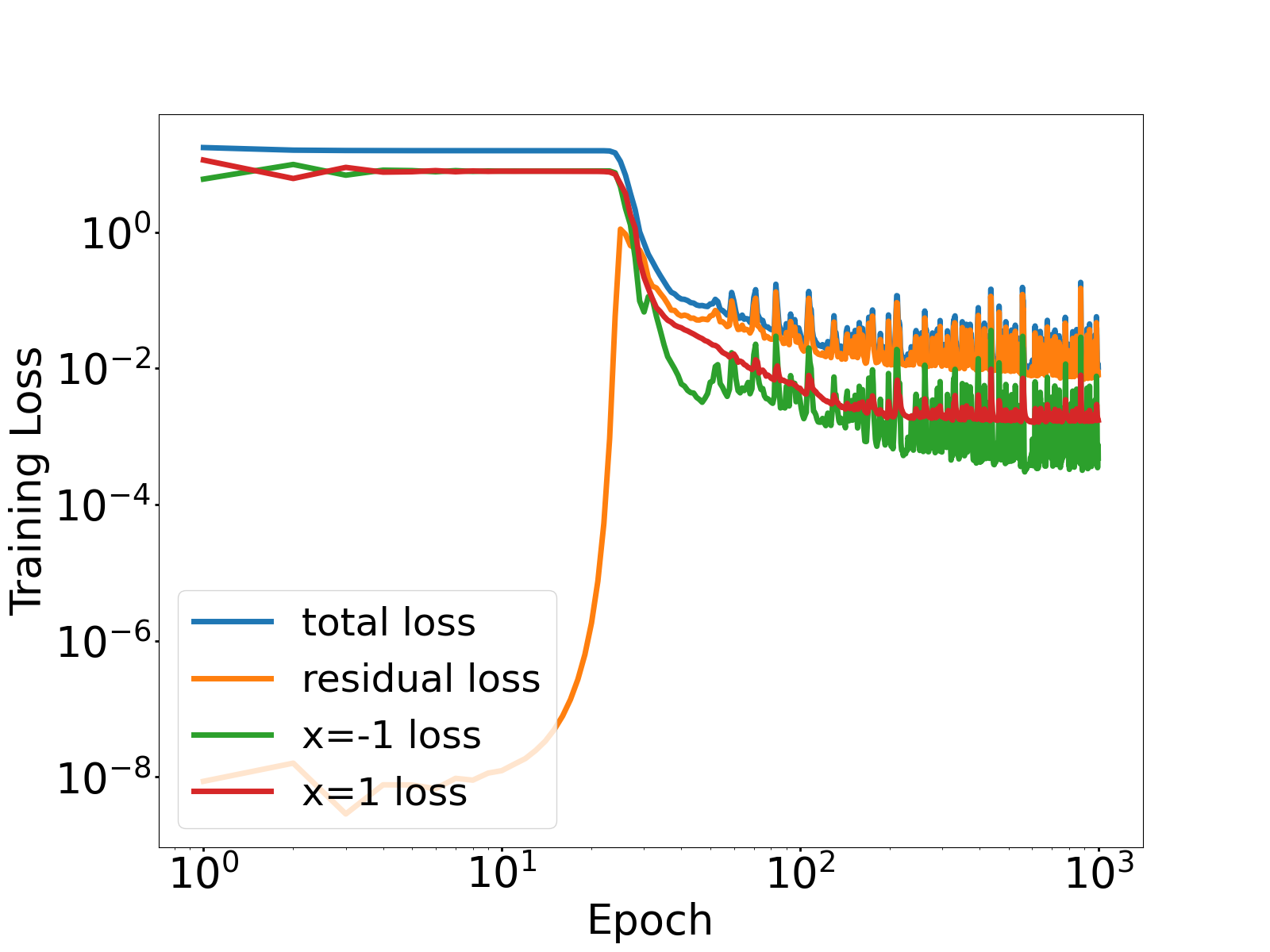}
    \caption{Example 5. Training loss vs. epoch for the learning of the operator of Burgers equation. The blue curve is the total risk Eq. (\ref{burgers_leastsquareloss}). The other three curves represent three terms in the total risk Eq. (\ref{burgers_leastsquareloss}), respectively.}
    \label{fig:burgers_op_loss}
\end{figure}

\begin{figure}[!htb]
    \centering
    \subfigure[$c_2$=3.1]{
        \includegraphics[scale=0.11]{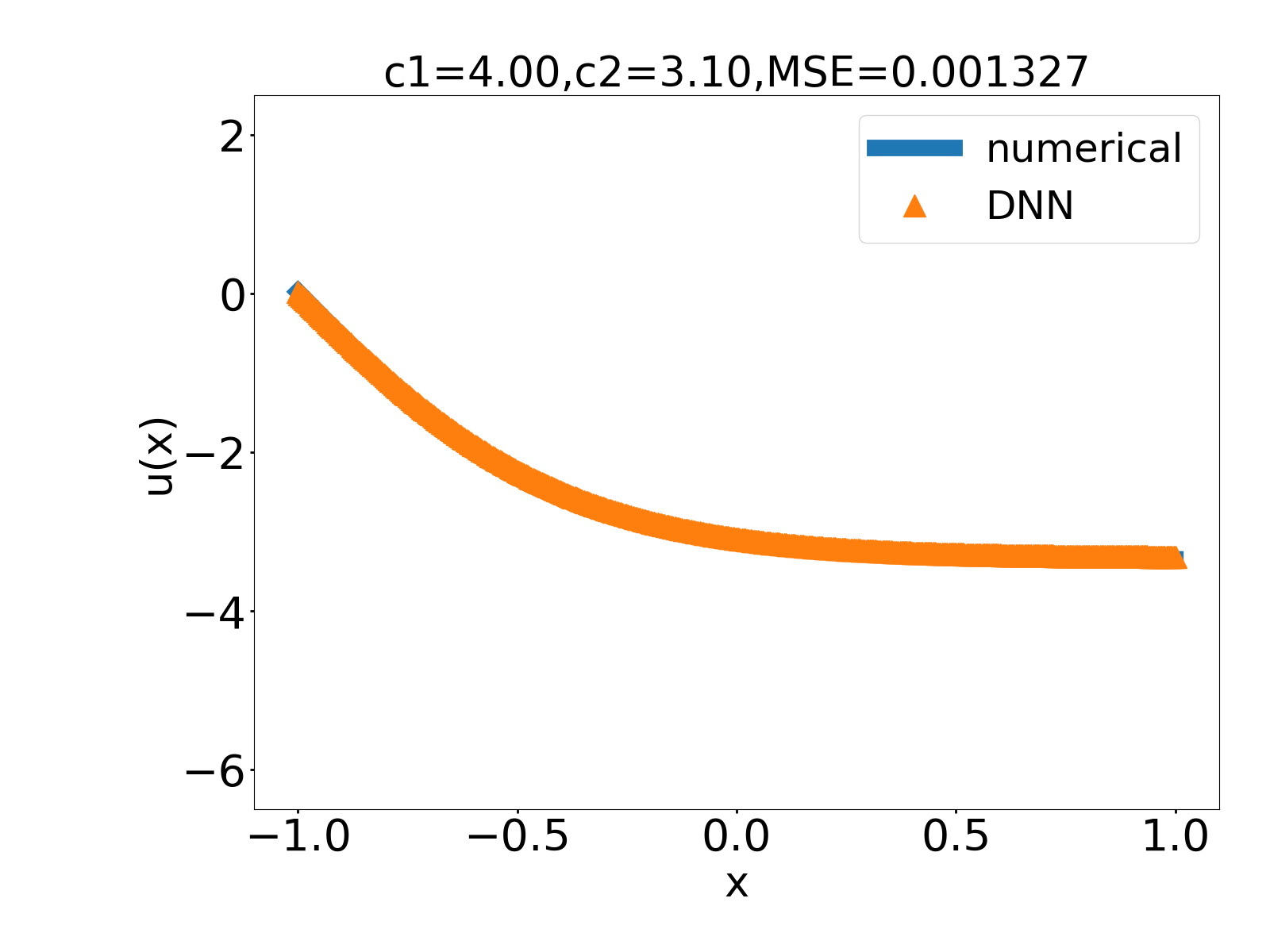}
    }
    \subfigure[$c_2$=4.1]{
        \includegraphics[scale=0.11]{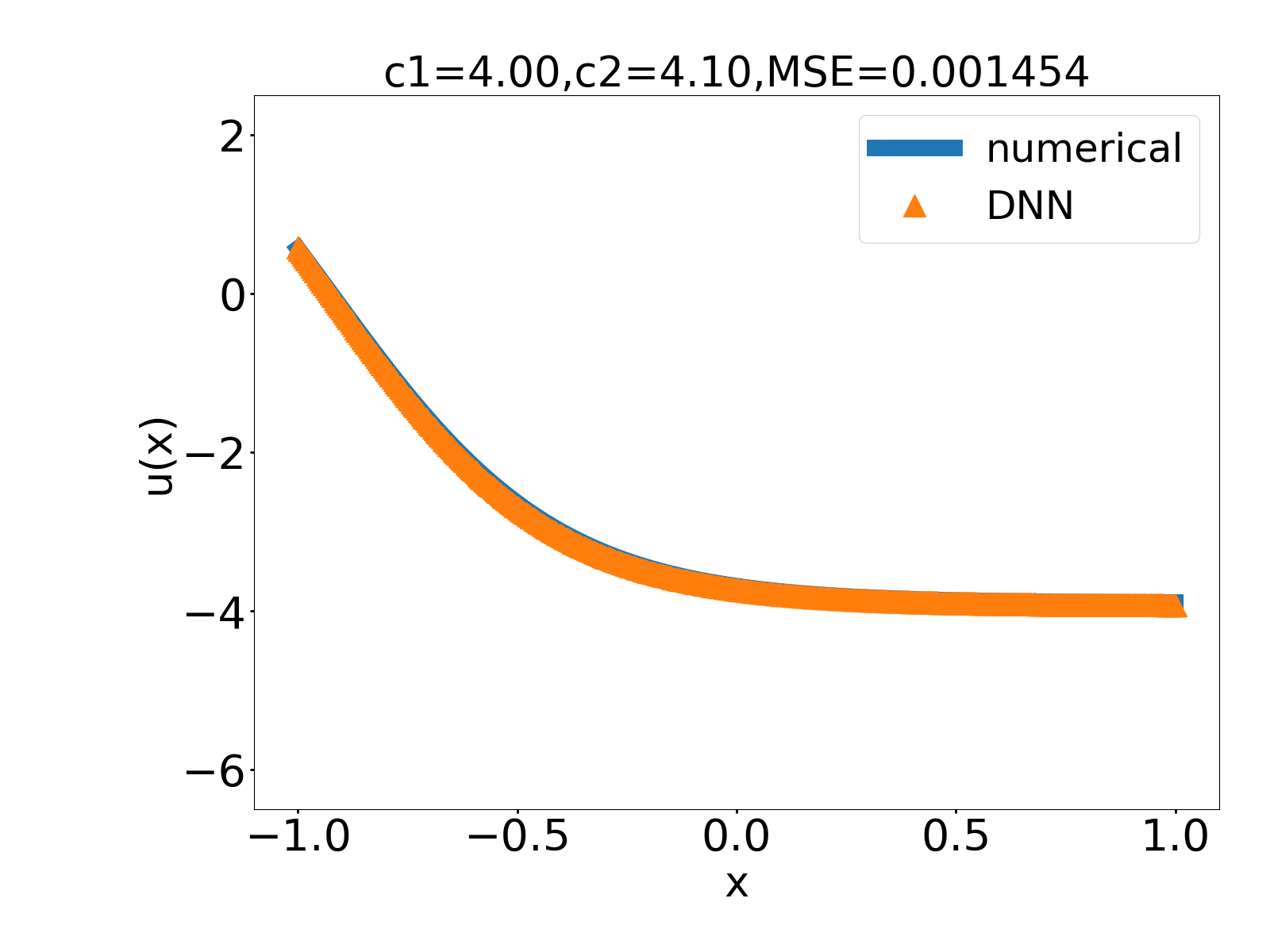}
    }
    \subfigure[$c_2$=5.1]{
        \includegraphics[scale=0.11]{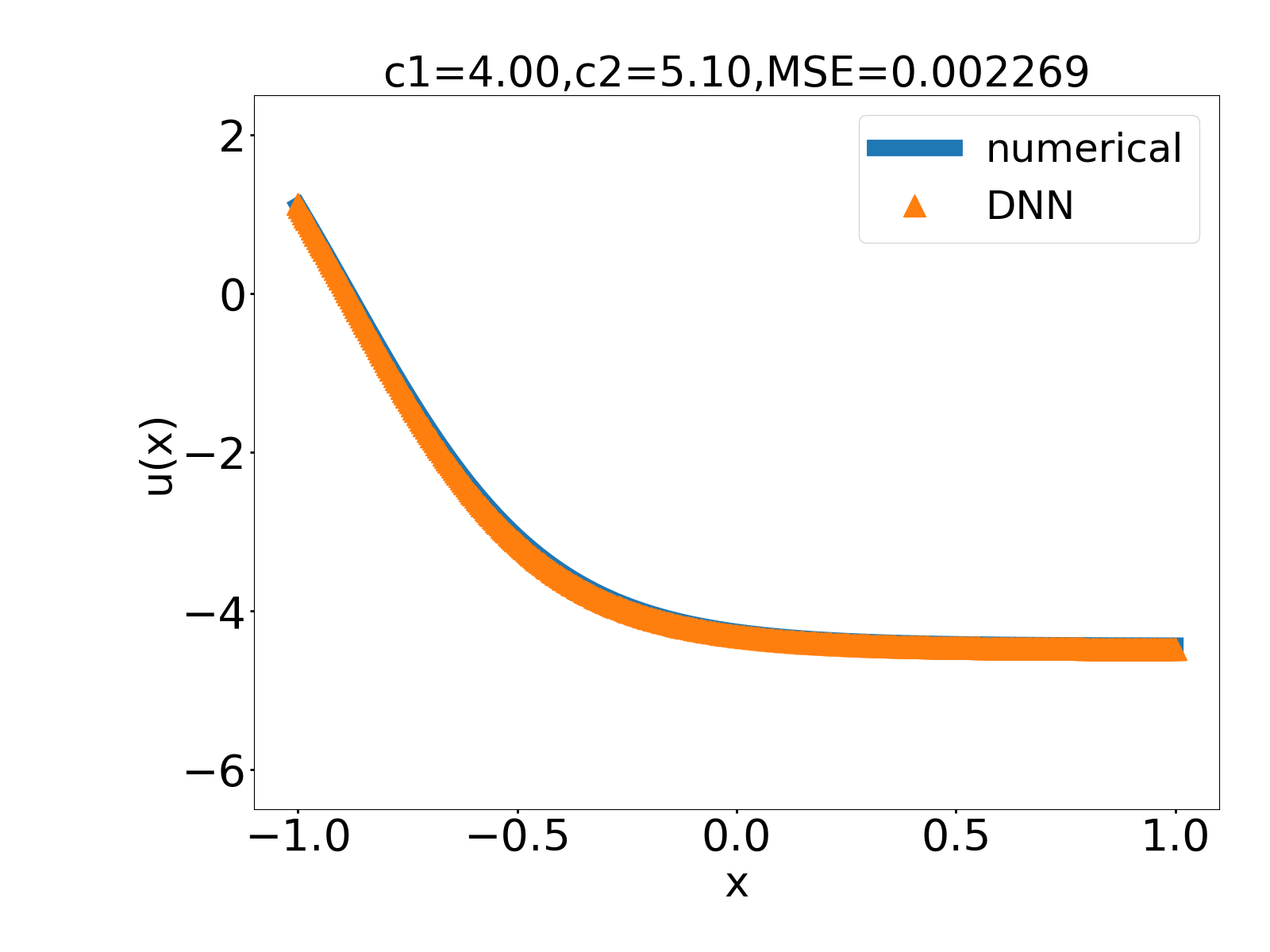}
    }
    \subfigure[$c_2$=6.1]{
        \includegraphics[scale=0.11]{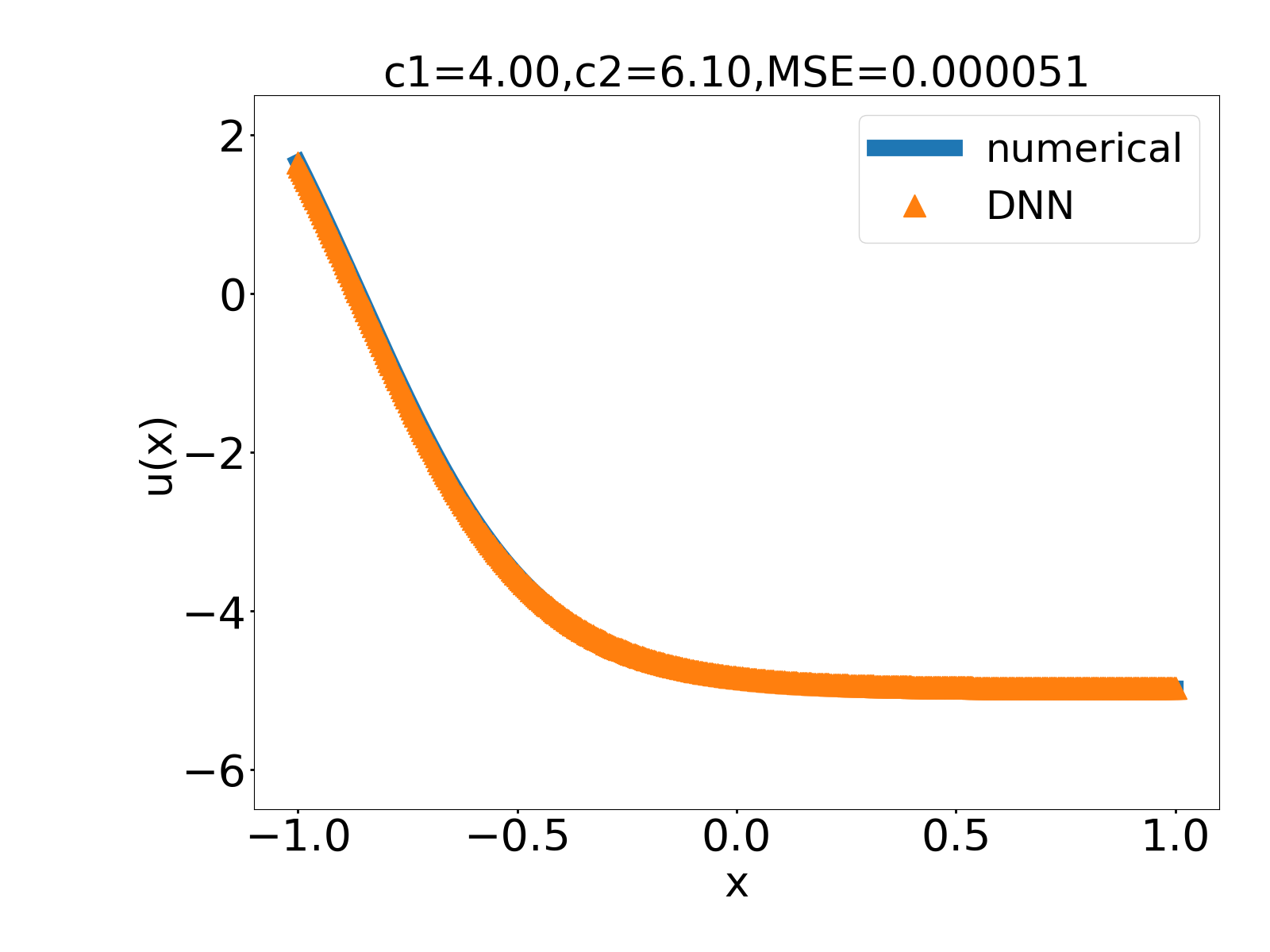}
    }
    \subfigure[$c_2$=6.9]{
        \includegraphics[scale=0.11]{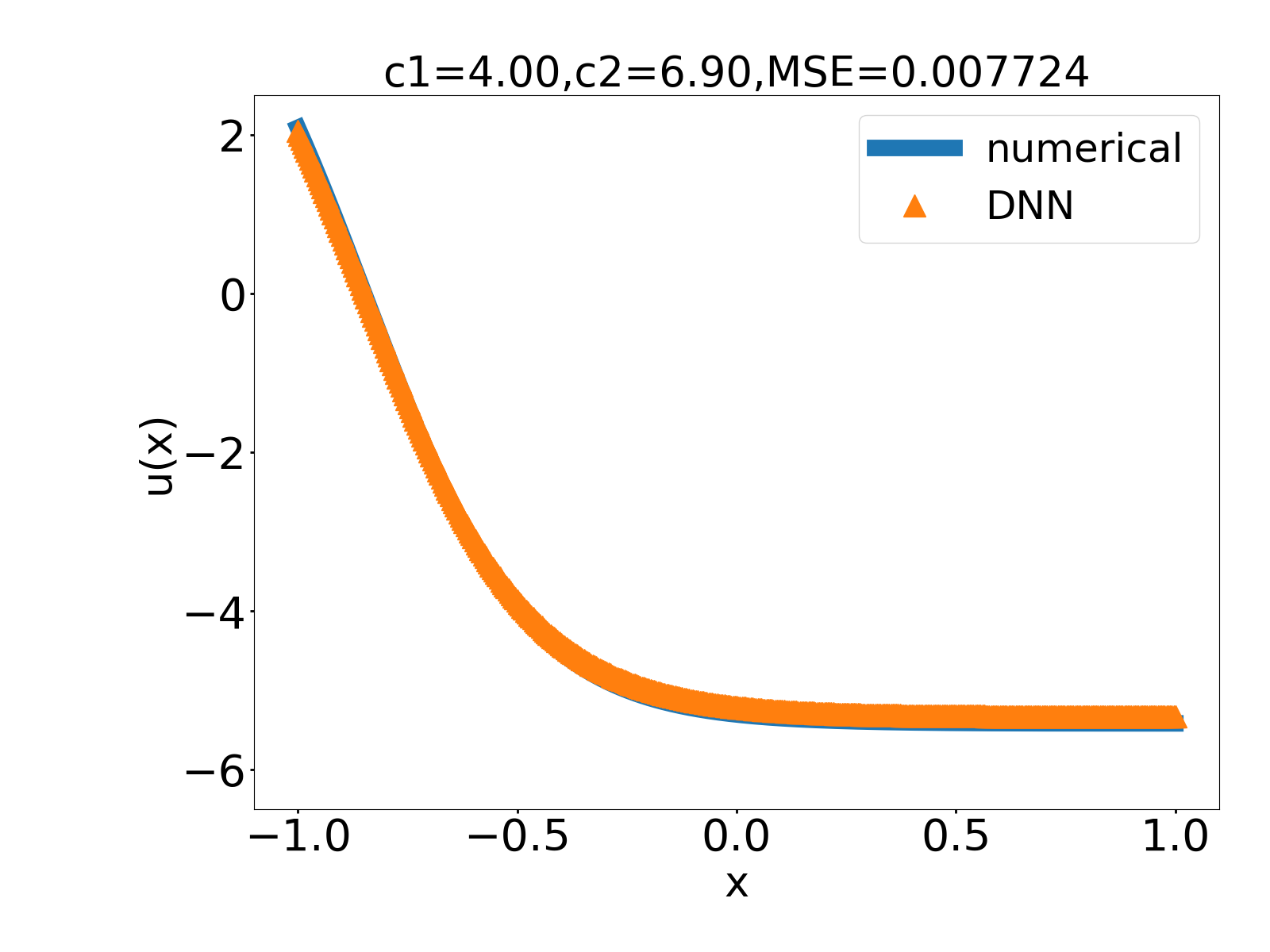}
    }
    \caption{Example 5. Comparison between numerical solution and  MOD-Net solution  on  equidistributed isometric points in $\Omega=[-1,1]$ corresponding to each of boundary conditions  determined by $\phi(x)=4\cos(2x) + c_2\sin(10x)$.}
    \label{fig:burgers_op}
\end{figure}

In training process, we set $K=21$ and for each epoch, we choose  $c_2$'s, i.e., $\{3+0.2(k-1)\}_{k=1}^{21}$.
For each $k$, we uniformly sample $500$ points in $[-1,1]$ for $S^{\Omega,S,k}$.
The training loss is shown in Fig. \ref{fig:burgers_op_loss}. The total risk (indicated by blue curve) roughly decays  with training epoch. For visualization of each part of the total risk, we also display the residual loss of the Burgers equation (indicated by orange curve) and the loss of two boundary loss (indicated by green and red curves), which  roughly decay  with the training epoch at the final stage.

To test the performance of the trained MOD-Net on $c_2=3.1,4.1,5.1,6.1,6.9$.
For each $c_2$,  we  plot the MOD-Net solution and the corresponding numerical solution on $500$ equidistributed isometric points sampled in $\Omega=[-1,1]$.
As shown in Fig. \ref{fig:burgers_op}, the MOD-Net solutions overlap well with the numerical solutions  for all cases.

\section{Numerical experiments: $2$D nonlinear equation}
\label{sec:nonlinear}


We consider the following nonlinear equation,
\begin{equation}\label{Poisson-nonlinear}
    \begin{aligned}
         & -\Delta u(\vx)+0.01u^3(\vx)=g(\vx), \quad \vx\in \Omega,
        \\
         & u(\vx)=0, \quad \vx\in \partial \Omega.
    \end{aligned}
\end{equation}
Take a $2$D case as example, in which source function $g(x,y)=-a(x^2-x+y^2-y)+0.01(\frac{a}{2}x(x-1)y(y-1))^3$,
\begin{equation}
    \begin{aligned}
         & -(\partial_{xx}u+\partial_{yy}u)+0.01u^3=-a(x^2-x+y^2-y)+0.01(\frac{a}{2}x(x-1)y(y-1))^3, \quad (x,y)\in \Omega,
        \\
         & u=0, \quad (x,y)\in \partial \Omega,
    \end{aligned}
    \label{2D_nonlinear_case}
\end{equation}
where $\Omega=[0,1]^2$ and $a$ controls the source term $g(x,y)$.
Obviously, the analytical solution is  $u(x,y;g)=\frac{a}{2}x(x-1)y(y-1)$.

\subsection{Use DNN to fit nonlinear operator}

In the considered $2$D case, $g(x,y)=-a(x^2-x+y^2-y)+0.01(\frac{a}{2}x(x-1)y(y-1))^3$, we have
\begin{equation}
    u(x,y;g) = F(\int_0^1\int_0^1 G(x,y,x',y')g(x',y') \diff{x'}\diff{y'}).
\end{equation}
We use a  DNN of hidden layer size $128$-$128$-$128$-$128$ $G_{\vtheta}(x,y,x',y')$ to fit the function $G(x,y,x',y')$ and a DNN of one hidden layer, the size of which is $256$, equipped with sigmoid function as activation function  $F_{\vtheta}(x,y)$.

When we calculate the integral, We use the Gauss-Legendre numerical integral. Then we can represent neural operator $u_{\vtheta}(x,y;g)$ with  DNN $G_{\vtheta}(x,y,x',y')$ and $F_{\vtheta}(x,y)$, that is,
\begin{equation}
    \begin{aligned}
        u_{\vtheta}(x,y;g) = F_{\vtheta}(\sum_{x^\prime \in S_{G_x}}\sum_{y^\prime \in S_{G_y}} \omega_{x^\prime} \omega_{y^\prime} G_{\vtheta}(x,y,x^\prime,y^\prime)g(x^\prime,y^\prime)),
    \end{aligned}
    \label{nonlinear_Gausianint}
\end{equation}
where  $S_{G_x}\subset[0,1]$,$S_{G_y}\subset[0,1]$ consist fixed integration points, determined by $1$D Gauss-Legendre quadrature and $\omega_{x^\prime}$, $\omega_{y^\prime}$ are corresponding coefficients.

\subsection{Empirical risk function}
For this  example, to train the neural networks, we use no labeled data and only utilize the information of PDE, i.e., governing equation and boundary condition, for each $g^k,k=1,2,\cdots,K$.
To utilize the  constraint of governing equation of PDE, we uniformly sample a set of data from $\Omega=[0,1] \times [0,1]$, i.e., $S^{\Omega,k}$. To utilize the information of boundary constraint, for each $k$, we uniformly sample a set of data from $\partial\Omega$.
Since the boundary $\partial\Omega$ consists of four line segments, i.e., $\partial\Omega=\bigcup_{i=1}^{4}\partial\Omega_i$, where $\partial\Omega_1=\{(0,y)\}_{y\in[0,1]}$,$\partial\Omega_2= \{(1,y)\}_{y\in[0,1]}$,$\partial\Omega_3=\{(x,0)\}_{x\in[0,1]}$,$\partial\Omega_4=\{(x,1)\}_{x\in[0,1]}$, we uniformly sample a set of data from $\partial\Omega_i$ respectively, i.e., $S^{\partial\Omega_i,k}$, i=1,2,3,4.

Since we use no labeled data, $\lambda_3$ in the general definition  is set as zero.
The empirical risk for this example is as follows,
\begin{equation}
    \begin{aligned}
        \RS
         & =  \frac{1}{K}\sum_{k\in[K]}  \bigg(\lambda_1
        \frac{1}{|S^{\Omega,k}|}\sum_{(x,y) \in S^{\Omega,k}}\Big(
        \partial_{xx}u_{\vtheta}(x,y;g^k)+\partial_{yy}u_{\vtheta}(x,y;g^k)+u_{\vtheta}(x,y;g^k)^3
        -g^k(x,y)
        \Big)^{2}
        \\
         & \quad +  \lambda_2 \frac{1}{|S^{\partial\Omega_1,k}|}\sum_{(x,y) \in S^{\partial\Omega_1,k}} {u_{\vtheta}(x,y;g^k)}^{2}
        \\
         & \quad +  \lambda_2 \frac{1}{|S^{\partial\Omega_2,k}|}\sum_{(x,y) \in S^{\partial\Omega_2,k}} {u_{\vtheta}(x,y;g^k)}^{2}
        \\
         & \quad +  \lambda_2\frac{1}{|S^{\partial\Omega_3,k}|}\sum_{(x,y) \in S^{\partial\Omega_3,k}} {u_{\vtheta}(x,y;g^k)}^{2}
        \\
         & \quad + \lambda_2\frac{1}{|S^{\partial\Omega_4,k}|}\sum_{(x,y) \in S^{\partial\Omega_4,k}} {u_{\vtheta}(x,y;g^k)}^{2}
        \bigg).
        \label{nonlinearcase_leastsquareloss}
    \end{aligned}
\end{equation}

\subsection{Learning process}

For each training epoch, we first randomly choose source functions $\{g^k\}_{k=1}^{K}$ and calculate their values on fixed  integration points$(x',y')$, where $x' \in S_{G_x}$ and $ y' \in S_{G_y}$. Second, we randomly sample data and obtain data set $S^{\Omega,k},S^{\partial \Omega_i,k},i=1,2,3,4$.
We obtain the data set $D= \{(x,y,x',y',g^k(x',y'))| (x,y)\in S^{\Omega,k} \cup \bigcup_{i=1}^{4} S^{\partial\Omega_i,k},x'\in S_{G_x},y' \in S_{G_y}\}$. In the following, we feed the data into the neural network $G_{\vtheta}(x,y,x',y')$ and calculate the total risk (\ref{nonlinearcase_leastsquareloss}) )  with neural operator $u_{\vtheta}(x,y;g)$, see Eq. (\ref{nonlinear_Gausianint}).
Train neural networks $G_{\vtheta}$ and $F_{\vtheta}$ with Adam to minimize the total risk, and finally we  obtain  well-trained  DNNs $G_{\vtheta}$ and $F_{\vtheta}$, furthermore, according to Eq. (\ref{nonlinear_Gausianint}), we obtain a neural operator $u_{\vtheta}(x,y;g)$.

\subsection{Results}
The source function $g(x,y)=-a(x^2-x+y^2-y)+0.01(\frac{a}{2}x(x-1)y(y-1))^3$ is determined by $a$. We then denote source function $g(x,y)$ by $g_a$.
To illustrate our approach,
we  train a neural operator mapping from $g$ to the solution of the  equation (\ref{2D_nonlinear_case}).
\subsubsection{Example 6: MOD-Net for a family of nonlinear equations}
For illustration, we would train the MOD-Net by various source functions and test the MOD-Net with several source functions that are not used for training.

\begin{figure}[!htb]
    \centering

    \subfigure[$a=15$]{
        \includegraphics[scale=0.11]{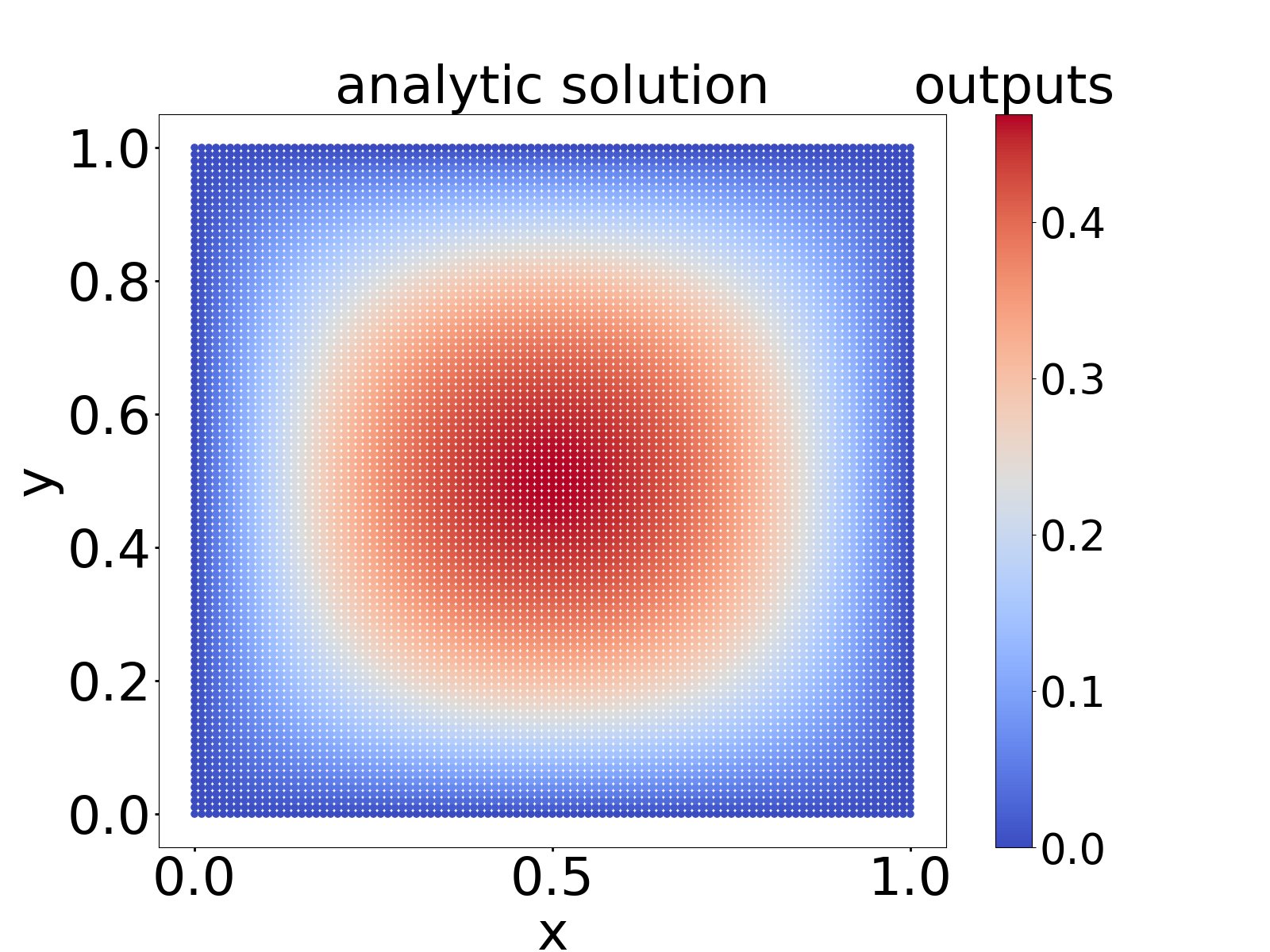}
    }
    \subfigure[$a=15$]{
        \includegraphics[scale=0.11]{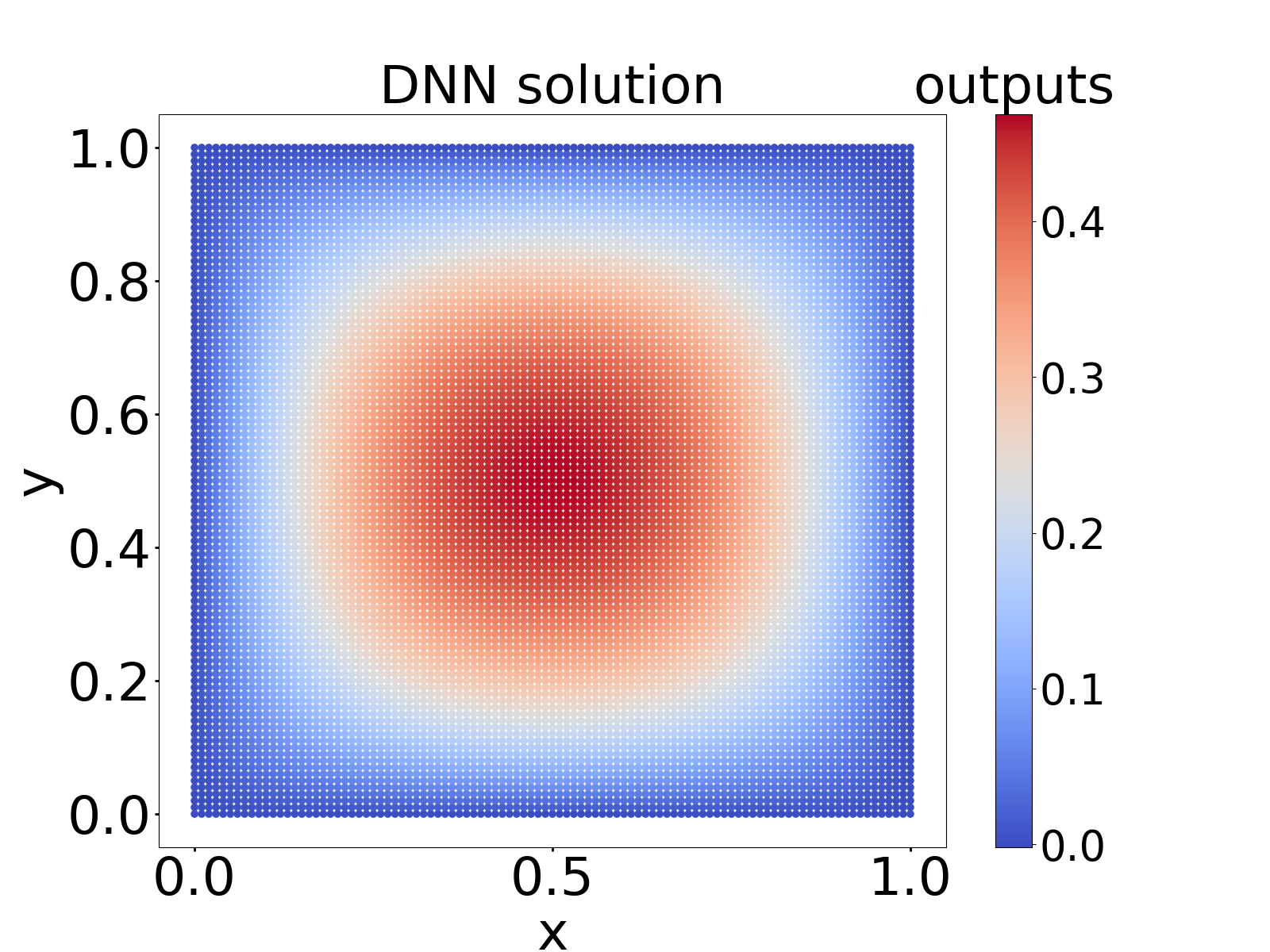}
    }
    \subfigure[$a=15$]{
        \includegraphics[scale=0.11]{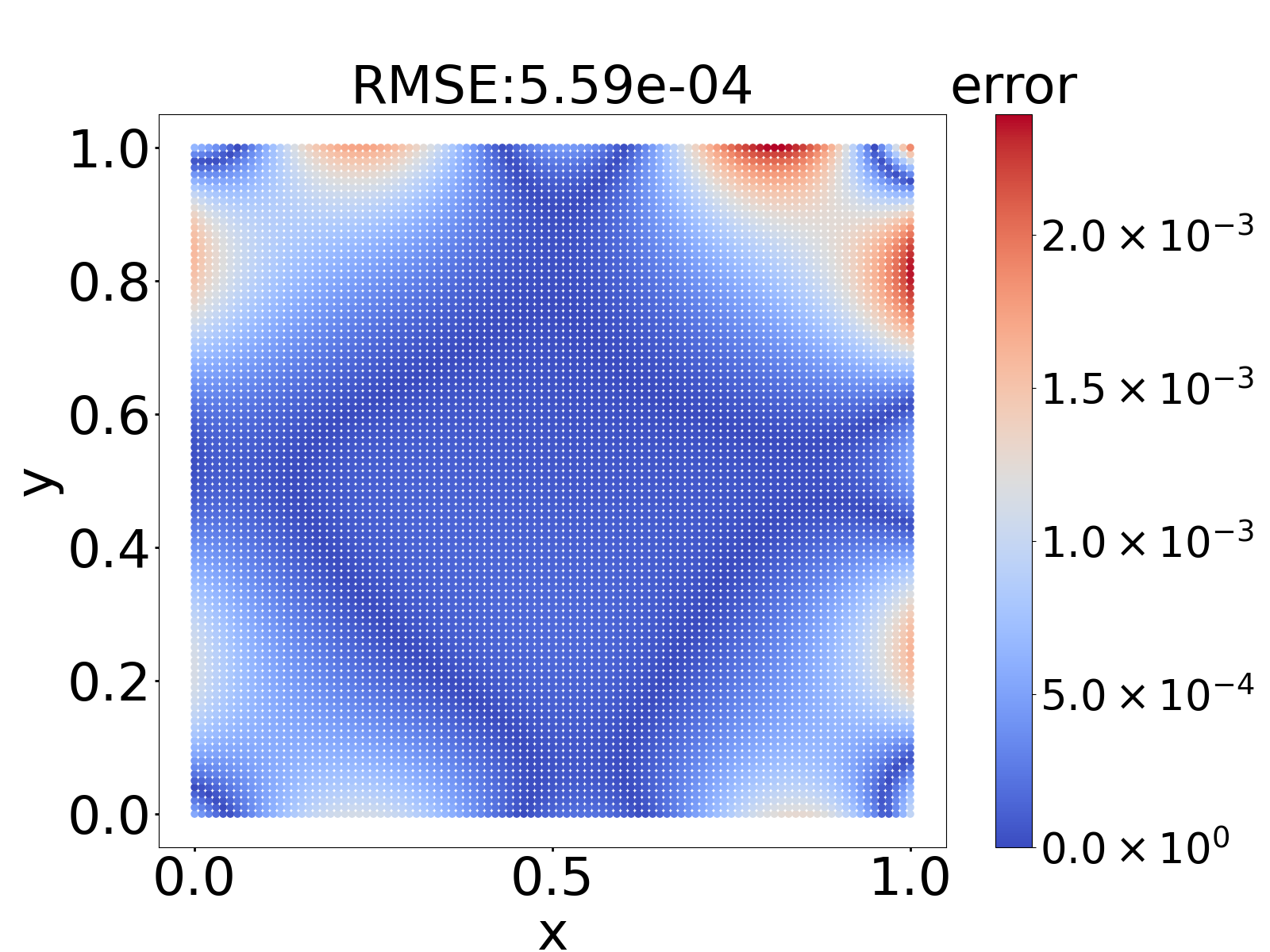}
    }
    \caption{Example 6. Comparison between analytic solution and  MOD-Net solution on $101\times101$ grid points corresponding to source terms determined by $a=15$. (a) Analytic solution. (b) MOD-Net solution. (c) The difference between the analytic solution and MOD-Net solution. Since the error at most points are very small, here we reduce the upper limit of the Color bar to see more information. }
    \label{fig:nonlinear_op_contrast_leastsquare_error}
\end{figure}

\begin{figure}[!htb]
    \centering

    \subfigure[$a=15,x=0.0$]{
        \includegraphics[scale=0.11]{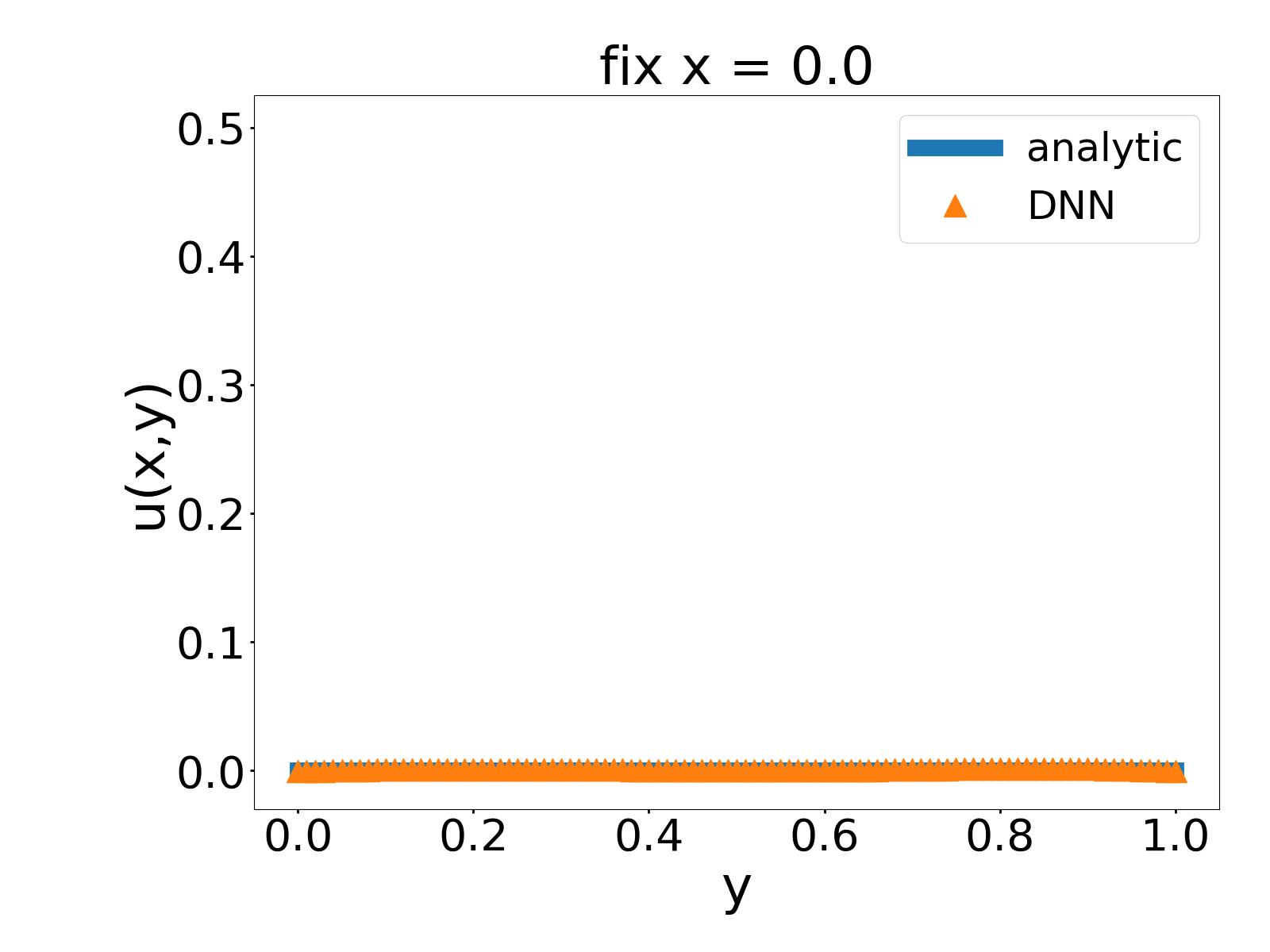}
    }
    \subfigure[$a=15,x=0.5$]{
        \includegraphics[scale=0.11]{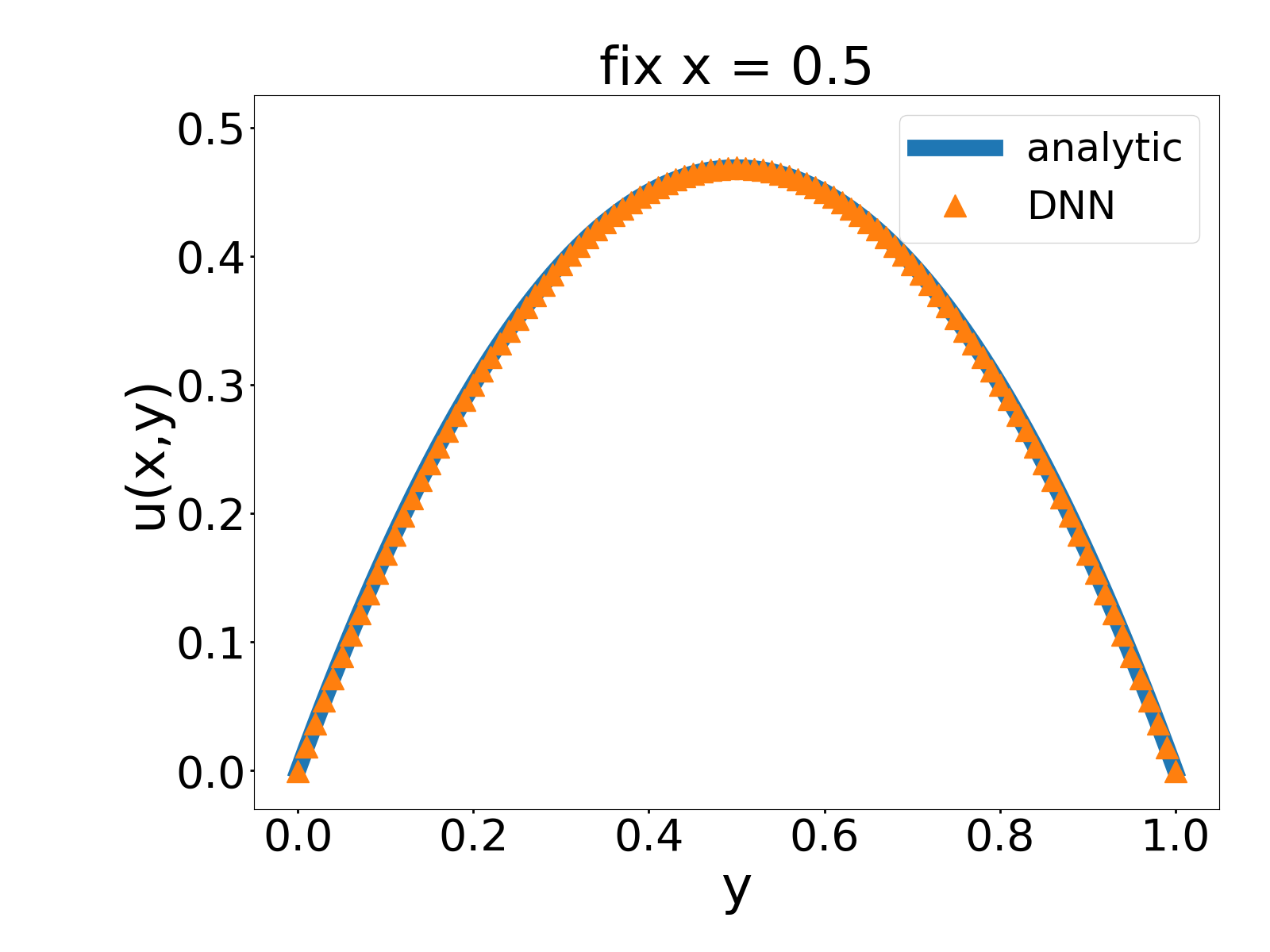}
    }
    \subfigure[$a=15,x=1.0$]{
        \includegraphics[scale=0.11]{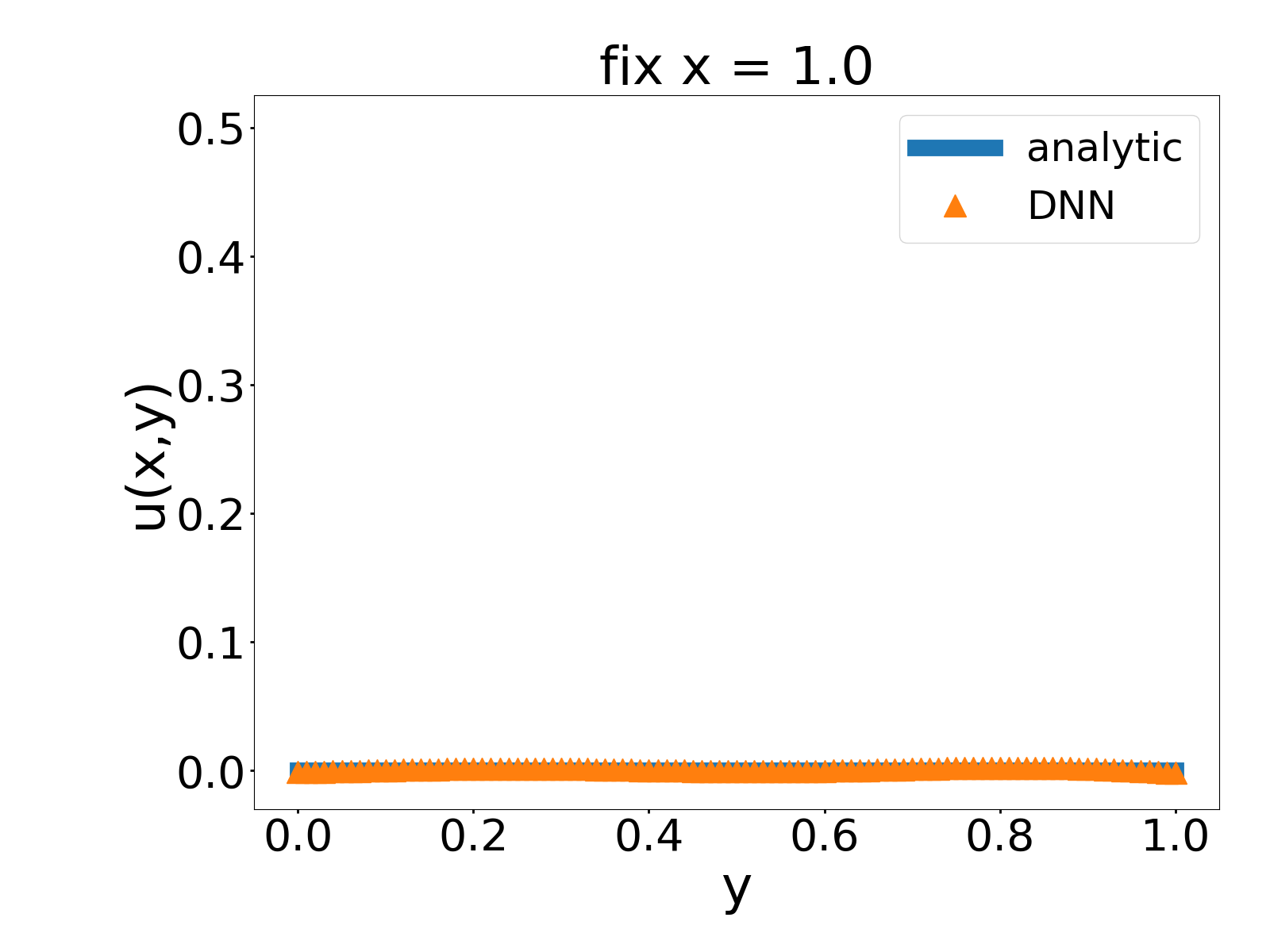}
    }
    \caption{Example 6. Comparison between analytic solution and  MOD-Net solution   on  $x=0,0.5,1$ corresponding to  source function determined by test $a=15$.}
    \label{fig:nonlinear_op_contrast_leastsquare_section}
\end{figure}

In training, we set $K=10$ and for each epoch, we choose a family of source functions $g_a(x,y)=-a(x^2-x+y^2-y)+0.01(\frac{a}{2}x(x-1)y(y-1))^3$
, where the control parameter $a \in \{10k\}_{k=1}^{10}$.
We set the number of integration points $|S_{G_x}|=10$ and $|S_{G_y}|=10$, and for each epoch, we randomly sample $200$ points in $[0,1]^2$ and sample  $200$ points in each line of boundary, respectively. We set $\lambda_1=1,\lambda_2=1$ in the empirical risk using least square loss (\ref{nonlinearcase_leastsquareloss}).

We test the performance of this well-trained MOD-Net on $a=15$.
For each fixed $a$, the corresponding exact true solution is $u(x,y)=\frac{a}{2}x(x-1)y(y-1)$.
For visualizing the performance of our well-trained MOD-Net, we show the analytic solution and  MOD-Net solution on $101\times101$ equidistributed isometric grid points by color. To compare these two solutions more intuitively, we calculate the difference of the two solutions at each point and show the error by color ,i.e., as the error increases, the color changes from blue to red, in Fig. \ref{fig:nonlinear_op_contrast_leastsquare_error}. The root mean square error (RMSE) is $\sim 5.6 \times 10^{-4}$. We also show the solution obtained by MOD-Net and the corresponding analytical solution on fixed $x=0,0.5,1$ for  considered test source function. As shown in  Fig. \ref{fig:nonlinear_op_contrast_leastsquare_section}, the MOD-Net can well predict the analytic solutions.

\section{Conclusion and discussion}
\label{sec:con}
In this work, we propose a model-operator-data network (MOD-Net) for solving PDEs.
The advantage of the MOD-Net is that it solves a family of PDEs but not a specific one in a meshless way without expensive labeled data.
For illustration, we use MOD-Net approach solving Poisson equation and a family of nonlinear PDEs, in which the empirical risk of MOD-Net only requires the physical information of the PDE and does not require any labeled data. And in some examples, such as constructed equation, which can be regarded as a simplification  of RTE, and the real one-dimensional RTE,  with only physical information and operator representation built into the empirical risk is not enough to guarantee the learning of correct solutions despite the small training loss.  By adding \emph{few labeled data as regularization} which are calculated analytically or computed by traditional numerical schemes on the coarse grid points with cheap computational cost, we show that the MOD-Net can be well trained to approximate the correct solution.

Though the good results are shown above, there are many important problems of efficiently using DNN to solve PDE for future works.
One is that it remains unclear what are the appropriate weights of different components to use in  the total loss. Second is to understand why the learned solver can be far away from the true solution even when residual loss and boundary loss are almost zero and why data regularization can rescue this phenomenon.
And Third, the MOD-Net is a framework to guide us how to take advantage of DNN to solve the physical PDEs. The aforementioned three key components should be adapted to specific PDEs.
For nonlinear operator it is still not clear what representations we should use, although many existing work propose their own representations~\cite{lu2019deeponet,wang2021learning,li2020fourier,li2020neural}
. For the data part, how many labeled data should we use and whether noisy data will have the same good regularization effect are all very important work to be done.
Other works need to be done including how we can use DNN to represent a Green's function related to variable coefficient, such as the $\sigma$ in RTE. Typically traditional DNN is not possible to parameterize an operator that is imposed on functions, so a careful design and build of the architecture is needed to feed information into the network of operator. We leave all this for future discussion.

\section*{Acknowledgments}
This work is sponsored by the National Key R\&D Program of China  Grant No. 2019YFA0709503 (Z. X.) and No. 2020YFA0712000 (Z. M.), the Shanghai Sailing Program (Z. X.), the Natural Science Foundation of Shanghai Grant No. 20ZR1429000  (Z. X.), the National Natural Science Foundation of China Grant No. 62002221 (Z. X.), the National Natural Science Foundation of China Grant No. 12101401 (T. L.), the National Natural Science Foundation of China Grant No. 12101402 (Y. Z.), Shanghai Municipal of Science and Technology Project Grant No. 20JC1419500 (Y. Z.), the National Natural Science Foundation of China Grant No. 12031013 (Z. M.), Shanghai Municipal of Science and Technology Major Project No. 2021SHZDZX0102, and the HPC of School of Mathematical Sciences and the Student Innovation Center at Shanghai Jiao Tong University.


\bibliographystyle{plain}
\bibliography{DLref}

\begin{thebibliography}{10}

\bibitem{cai2020phase}
Wei Cai, Xiaoguang Li, and Lizuo Liu.
\newblock A phase shift deep neural network for high frequency approximation
  and wave problems.
\newblock {\em SIAM Journal on Scientific Computing}, 42(5):A3285--A3312, 2020.

\bibitem{chandrasekhar2013radiative}
Subrahmanyan Chandrasekhar.
\newblock {\em Radiative transfer}.
\newblock Courier Corporation, 2013.

\bibitem{dissanayake1994neural}
MWMG Dissanayake and Nhan Phan-Thien.
\newblock Neural-network-based approximations for solving partial differential
  equations.
\newblock {\em communications in Numerical Methods in Engineering},
  10(3):195--201, 1994.

\bibitem{weinan2020machine}
Weinan E.
\newblock Machine learning and computational mathematics.
\newblock {\em arXiv preprint arXiv:2009.14596}, 2020.

\bibitem{weinan2017deep}
Weinan E, Jiequn Han, and Arnulf Jentzen.
\newblock Deep learning-based numerical methods for high-dimensional parabolic
  partial differential equations and backward stochastic differential
  equations.
\newblock {\em Communications in Mathematics and Statistics}, 5(4):349--380,
  2017.

\bibitem{E2020algorithms}
Weinan E, Jiequn Han, Arnulf Jentzen, et~al.
\newblock Algorithms for solving high dimensional pdes: From nonlinear monte
  carlo to machine learning.
\newblock {\em arXiv preprint arXiv:2008.13333}, 2020.

\bibitem{weinan2018deep}
Weinan E and Bing Yu.
\newblock The deep ritz method: A deep learning-based numerical algorithm for
  solving variational problems.
\newblock {\em Communications in Mathematics and Statistics}, 6(1):1--12, 2018.

\bibitem{fan2019multiscale}
Yuwei Fan, Lin Lin, Lexing Ying, and Leonardo Zepeda-N{\'u}nez.
\newblock A multiscale neural network based on hierarchical matrices.
\newblock {\em Multiscale Modeling \& Simulation}, 17(4):1189--1213, 2019.

\bibitem{hamilton2019dnn}
A~Hamilton, T~Tran, MB~Mckay, B~Quiring, and PS~Vassilevski.
\newblock Dnn approximation of nonlinear finite element equations.
\newblock Technical report, Lawrence Livermore National Lab.(LLNL), Livermore,
  CA (United States), 2019.

\bibitem{he2020relu}
Juncai He, Lin Li, Jinchao Xu, and Chunyue Zheng.
\newblock Relu deep neural networks and linear finite elements.
\newblock {\em Journal of Computational Mathematics}, 38(3):502--527, 2020.

\bibitem{lenoble1985radiative}
Jacqueline Lenoble.
\newblock {\em Radiative transfer in scattering and absorbing atmospheres:
  standard computational procedures}, volume 300.
\newblock A. Deepak Hampton, Va., 1985.

\bibitem{li2020multi}
Xi-An Li, Zhi-Qin~John Xu, and Lei Zhang.
\newblock A multi-scale dnn algorithm for nonlinear elliptic equations with
  multiple scales.
\newblock {\em Communications in Computational Physics}, 28(5):1886--1906,
  2020.

\bibitem{li2020fourier}
Zongyi Li, Nikola Kovachki, Kamyar Azizzadenesheli, Burigede Liu, Kaushik
  Bhattacharya, Andrew Stuart, and Anima Anandkumar.
\newblock Fourier neural operator for parametric partial differential
  equations.
\newblock {\em arXiv preprint arXiv:2010.08895}, 2020.

\bibitem{li2020neural}
Zongyi Li, Nikola Kovachki, Kamyar Azizzadenesheli, Burigede Liu, Kaushik
  Bhattacharya, Andrew Stuart, and Anima Anandkumar.
\newblock Neural operator: Graph kernel network for partial differential
  equations.
\newblock {\em arXiv preprint arXiv:2003.03485}, 2020.

\bibitem{liao2021deep}
Yulei Liao and Pingbing Ming.
\newblock Deep nitsche method: Deep ritz method with essential boundary
  conditions.
\newblock {\em Communications in Computational Physics}, 29(5):1365--1384,
  2021.

\bibitem{liu2020multi}
Ziqi Liu, Wei Cai, and Zhi-Qin~John Xu.
\newblock Multi-scale deep neural network (mscalednn) for solving
  poisson-boltzmann equation in complex domains.
\newblock {\em Communications in Computational Physics}, 28(5):1970--2001,
  2020.

\bibitem{lu2019deeponet}
Lu~Lu, Pengzhan Jin, and George~Em Karniadakis.
\newblock Deeponet: Learning nonlinear operators for identifying differential
  equations based on the universal approximation theorem of operators.
\newblock {\em arXiv preprint arXiv:1910.03193}, 2019.

\bibitem{markidis2021physics}
Stefano Markidis.
\newblock Physics-informed deep-learning for scientific computing.
\newblock {\em arXiv preprint arXiv:2103.09655}, 2021.

\bibitem{raissi2019physics}
Maziar Raissi, Paris Perdikaris, and George~E Karniadakis.
\newblock Physics-informed neural networks: A deep learning framework for
  solving forward and inverse problems involving nonlinear partial differential
  equations.
\newblock {\em Journal of Computational Physics}, 378:686--707, 2019.

\bibitem{strofer2019data}
Carlos~Michelen Strofer, Jin-Long Wu, Heng Xiao, and Eric Paterson.
\newblock Data-driven, physics-based feature extraction from fluid flow fields
  using convolutional neural networks.
\newblock {\em Communications in Computational Physics}, 25(3):625--650, 2019.

\bibitem{Matthew2020Fourier}
Matthew {Tancik}, Pratul~P. {Srinivasan}, Ben {Mildenhall}, Sara
  {Fridovich-Keil}, Nithin {Raghavan}, Utkarsh {Singhal}, Ravi {Ramamoorthi},
  Jonathan~T. {Barron}, and Ren {Ng}.
\newblock Fourier features let networks learn high frequency functions in low
  dimensional domains.
\newblock {\em arXiv preprint arXiv:2006.10739}, 2020.

\bibitem{wang2021learning}
Sifan Wang, Hanwen Wang, and Paris Perdikaris.
\newblock Learning the solution operator of parametric partial differential
  equations with physics-informed deeponets.
\newblock {\em arXiv preprint arXiv:2103.10974}, 2021.

\bibitem{wang2021eigenvector}
Sifan {Wang}, Hanwen {Wang}, and Paris {Perdikaris}.
\newblock {On the eigenvector bias of Fourier feature networks: From regression
  to solving multi-scale PDEs with physics-informed neural networks}.
\newblock {\em Computer Methods in Applied Mechanics and Engineering}, 384(C),
  2021.

\bibitem{xu2019frequency}
Zhi-Qin~John Xu, Yaoyu Zhang, Tao Luo, Yanyang Xiao, and Zheng Ma.
\newblock Frequency principle: Fourier analysis sheds light on deep neural
  networks.
\newblock {\em Communications in Computational Physics}, 28(5):1746--1767,
  2020.

\bibitem{xu_training_2018}
Zhi-Qin~John Xu, Yaoyu Zhang, and Yanyang Xiao.
\newblock Training {Behavior} of {Deep} {Neural} {Network} in {Frequency}
  {Domain}.
\newblock In {\em Neural {Information} {Processing}}, Lecture {Notes} in
  {Computer} {Science}, pages 264--274, 2019.

\end{thebibliography}

\end{document}